\documentclass[10pt]{amsart}
\usepackage{amsthm,amsmath,amssymb}
\usepackage{enumitem}
\usepackage{bm}
\usepackage[
pdfencoding=auto,
colorlinks = true,
citecolor = blue,
linkcolor = blue,
anchorcolor = blue,
urlcolor = blue, 
filecolor=blue,
pdfkeywords={}
]{hyperref}
\usepackage{color}

\newtheorem{teor}{Theorem}[section]

\newtheorem{lemm}[teor]{Lemma}
\newtheorem{osse}[teor]{Remark}
\newtheorem{prop}[teor]{Proposition}
\newtheorem{defi}[teor]{Definition}
\newtheorem{coro}[teor]{Corollary}

\newcommand{\xe}{\bm{\xi}}
\newcommand{\ve}{\mathbf{v}}
\newcommand{\He}{\mathbf{H}}
\newcommand{\qe}{\mathbf{q}}
\newcommand{\no}{\mathbf{n}_{\partial \Omega}}
\newcommand{\opdiv}{\operatorname{div}}

\newcommand{\dx}{\mathrm{d}x}
\newcommand{\dt}{\mathrm{d}t}

\def\d{{\rm d}}

\def\H{\mathbf{H}}
\def\W{\mathbf{W}}

\def\L{\mathbf{L}}

\def\ddt{\frac{\d}{\d t}}

\def\n{\mathbf{n}}

\newcommand{\numberset}{\mathbb}
\newcommand{\N}{\numberset{N}}
\newcommand{\R}{\numberset{R}}

\def\P{\mathbf{P}}
\def\v{\mathbf{v}}

\def\vphi{\varphi}

\def\eps{\varepsilon}

\def\phit_{\phi_t}

\def\tT{{T}}
\def\H{\mathbf{H}}
\def\Zut{Z^1_{\tT}}
\def\Zdt{Z^2_{\tT}}
\def\Ztt{Z^3_{\tT}}

\def\XT{X_{{T}}}
\def\YT{Y_{{T}}}
\def\norm#1{\left\Vert#1\right\Vert }
\def\norma#1{\left\vert#1\right\vert }

\def\f{\mathbf f}

\def\ddt{\frac{d}{dt}}

\def\LLQ(#1,#2){L^{#1}(0,\widetilde{T};\L^{#2}(\Gamma_0))}

\def\div#1{\text{div}_{#1}}

\def\non{\nonumber}

\def\Ya{L^p(0,T;\L^p_\sigma(\Omega))}
\def\Yan{L^p(0,T;\L^p(\Omega))}
\def\Yb{L^q(0,T;L^q(\Omega))}
\def\Yc{L^r(0,T;L^r(\Omega))}
\def\th{\theta}
\def\Vs{\mathbf{V}_\sigma}
\def\Hs{\mathbf H_\sigma}

\def\LL{\mathbf L}
\def\bv{\mathbf v}
\def\bV{\mathbf V}
\def\vdn{\bv_\delta^n}
\def\vphidn{\vphi_\delta^n}
\def\thetadn{\theta_\delta^n}
\def\qdn{q_\delta^n}

\numberwithin{equation}{section}
\begin{document}

	\title[A non-isothermal Navier--Stokes/Allen--Cahn model]{Well-posedness and sharp interface limit of a non-isothermal Navier--Stokes/Allen--Cahn model}

\author[H. Abels]{Helmut Abels} 
\address[Helmut Abels]{Faculty of Mathematics, Universit\"at Regensburg,
	Universitätsstraße 31, 93053 Regensburg, 
	Germany}
\email{helmut.abels@mathematik.uni-regensburg.de}

\author[A. Marveggio]{Alice Marveggio}
\address[Alice Marveggio]{Hausdorff Center for Mathematics, Universit\"at Bonn, Endenicher Allee 62, 53115 Bonn,
	Germany}
\email{\tt alice.marveggio@hcm.uni-bonn.de}

\author[A. Poiatti]{Andrea Poiatti}
\address[Andrea Poiatti]{Faculty of Mathematics, University of Vienna,	Oskar-Morgenstern-Platz 1, 1090 Wien, Austria}
\email{\tt andrea.poiatti@univie.ac.at}


\begin{abstract}
	We propose a thermodynamically consistent phase-field model for the flow of a mixture of two different viscous incompressible fluids of equal density in a bounded domain.
	We prove the well-posedness of local-in-time strong solutions by means of maximal regularity and contraction mapping arguments. 
	We introduce a suitable entropic weak formulation of the problem, replacing the heat equation by the total energy inequality and an entropy production inequality, and we rigorously prove global-in-time existence of such weak solutions, developing a novel approximation scheme.
    We also show that an entropic weak solution to this non-isothermal phase-field model converges to a distributional (or $BV$) solution to a non-isothermal Navier--Stokes/mean curvature flow, under an energy convergence assumption.
\end{abstract}

\maketitle

\tableofcontents

\section{Introduction}
\label{sec:intro}
In this manuscript we propose the study of a non-isothermal diffuse interface model describing the evolution of a mixture of two viscous incompressible Newtonian fluids of equal density and viscosity in a bounded domain $\Omega \subset \mathbb{R}^d$, $d=2,3$, and in a time interval $[0,T)$, $T>0$.  
Our model couples the incompressible Navier-Stokes system for the velocity $\v_\eps=\v_\eps(x,t)$, as a single velocity field for the mixture, with an Allen-Cahn system with convection, which also accounts for the effects of the (absolute) temperature \(\theta_\eps=\theta_\eps(x,t)\). 
 We introduce the scalar field $\varphi_\eps=\varphi_\eps(x,t)$ as the phase-field variable or (non-conserved) order parameter (e.g., representing the rescaled density difference of the two involved phases), which takes value $-1$ and resp. 1 in the pure states, while it varies from $-1$ to 1 in diffused transition layers of typical width $\varepsilon > 0$.

Denoting by $ \frac{D  }{Dt} $ the material derivative $\partial_t + \ve_\varepsilon \cdot \nabla  $, we consider the following system of PDEs: 
 \begin{align}
	\partial_t \ve_\varepsilon + \ve_\varepsilon \cdot \nabla \ve_\varepsilon &= - \nabla p + \opdiv (\mathbb{S}_\varepsilon) - \opdiv(\varepsilon \nabla \varphi_{\varepsilon}  \otimes \nabla \varphi_{\varepsilon} ), \label{eq:intro1}\\
	\opdiv(\ve_\varepsilon)&=0, \label{eq:intro2}\\
	\frac{D \varphi_{\varepsilon} }{Dt} &= \Delta \varphi_\varepsilon - \frac{1}{\varepsilon^2} \partial_\varphi W(\varphi_\varepsilon) + \frac{1}{\varepsilon}\ell(\theta_\varepsilon) , \label{eq:intro3}\\
	\frac{D Q(\theta_\varepsilon )  }{Dt}+ \ell(\theta_\varepsilon) \frac{D \varphi_{\varepsilon}  }{Dt} &= \opdiv(\kappa(\theta_\varepsilon) \nabla \theta_\varepsilon) + \nu(\theta_\varepsilon) |\nabla \ve_\varepsilon + \nabla \ve_\varepsilon^\mathsf{T}|^2 + \varepsilon \left| \frac{D \varphi_{\varepsilon} }{Dt}\right|^2, \label{eq:intro4}
\end{align}
on $\Omega \times (0,T)$, where the scalar function \(p=p(x,t)\) is the local pressure, $\mathbb{S}(\theta_\varepsilon)\mathbf v_\varepsilon= {\nu(\theta_\eps)} (\nabla \ve_\eps + \nabla \ve_\eps^\mathsf{T})$ represents the viscous (i.e., dissipative) part of the stress tensor, where$\nu(\theta_\eps)>0$ is the viscosity of the mixture, $W$ is a standard polynomial double well potential with zeros at $\pm 1$, i.e., $W(s)= \frac14(1-s^2)^2$,  $Q(\theta_\eps)$ denotes the heat capacity,
 $\kappa(\theta_\eps)$ indicates the heat conductivity,
and $\ell(\theta_\eps)$ is the latent heat. 
Suitable assumptions on the functions $\nu, Q, \kappa$, and $\ell$ of the temperature $\theta_\eps$ will be introduced in Section \ref{sec:main}.
We endow the system of PDEs  \eqref{eq:intro1}-\eqref{eq:intro4} with suitable initial and boundary conditions, assuming that the system is insulated from the exterior, i.e.,
complete slip conditions for $\ve_\varepsilon$, homogeneous Neumann condition for $\varphi_{\varepsilon} $, and no (heat) flux condition for $\theta_\varepsilon$ (for more details, see Section \ref{sec:derivation}).

 A non-isothermal Navier-Stokes/Allen-Cahn system was first systematically presented by T. Blesgen \cite{Blesgen}. The resulting model is a first step to describe cavitation in a flowing liquid: gas phases grow or shrink due to changes of temperature and density in the fluid and they are transported with the current. A generalization of this model to a broader class of fluids has been derived by means of the first principles of thermodynamics by H. Freist\"uhler and M. Kotschote \cite{FK}.
 Here we obtain the system of PDEs \eqref{eq:intro1}-\eqref{eq:intro1} by following the general approach devised Frémond \cite{Fremond} for the physical derivation of phase-field models in a non-isothermal setting (see Section \ref{sec:derivation},  also \cite{ERS3D}).

Even though there are several works on Allen-Cahn models coupled with the Navier-Stokes equation in the constant temperature case (e.g., \cite{GGP,GigaLiu,Li, numerics}, see also \cite{GP} for a multi-component version of the Allen-Cahn equation), fewer contributions in the literature concern the non-isothermal setting.
The well-posedness of a system consisting in the Navier–Stokes equation, a convective Allen–Cahn equation, and an energy transport equation for the temperature was studied in \cite{Lopes}.
A model describing the evolution of a liquid crystal substance in the nematic phase was investigated in \cite{FeirRocSchi}.
The proposed model \eqref{eq:intro1}-\eqref{eq:intro4} couples the incompressible Navier-Stokes system with an Allen-Cahn type equation, also supplemented by a heat equation. 
Since an uncontrolled dissipation of the kinetic energy may not be captured by any term appearing in the classical formulation of the problem and since the loss of kinetic energy in any energetically closed system must be compensated by a source term in the internal energy balance, mathematical difficulties arise whenever the equations for the kinetic and internal (heat) energy are separated. To avoid this difficulty, existence of global-in-time weak solutions for the resulting PDE system is thus proved for a suitable weak formulation of the problem replacing the heat equation with the total energy balance and an entropy production inequality. 
This notion of weak solution was originally developed in \cite{Feireisl1} in the framework of heat conduction phenomena in fluids. It embodies the basic principles of thermomechanics and ensures the thermodynamical consistency of the model.  
Subsequently, this ``entropic" notion of weak solution has been used to prove global-in-time existence results for non-isothermal systems, including models for liquid crystals \cite{FFRS,FRSZ,FRSZ2}, phase transitions and damage phenomena in thermoviscoelastic materials \cite{RoccaRossi}, and binary incompressible immiscible fluids \cite{ERS3D}. 
In particular, in \cite{ERS3D} the authors study a model coupling the incompressible Navier-Stokes equation with a Cahn-Hilliard system with convection, which also accounts for the effects of the absolute temperature (see also \cite{ERS2D,Lasarzik}). A global-in-time weak existence result was obtained in \cite{ERS3D}, under specific growth conditions for the heat capacity $Q$ and the heat conductivity $\kappa$, excluding some physically relevant cases. In \cite{Lasarzik}, an extension to these cases is provided by considering measure-valued solutions and, more generally, dissipative solutions, relying on an appropriate (relative) energy inequality. 
In case $d=2$ and for linear $Q$, global-in-time existence for a stronger formulation of the model was proved in \cite{ERS2D}.
Similar techniques were also adopted in \cite{LRS}, where an Allen-Cahn type equation with internal energy balance was considered. 

Inspired by the works mentioned above, we investigate the well-posedness and the existence of solutions for the thermodynamically consistent phase-field model \eqref{eq:intro1}-\eqref{eq:intro4} under suitable assumptions on the functions $\nu, Q, \kappa$, and $\ell$ of the temperature $\theta_\eps$ (see Section \ref{sec:main}).
To this aim, we first prove the existence of a unique local-in-time strong solution to \eqref{eq:intro1}-\eqref{eq:intro4} by means of maximal regularity and a contraction mapping argument (see Theorem \ref{thm:strong} and Section \ref{sec::strong}). 
Then, we introduce a suitable weak solution concept for \eqref{eq:intro1}-\eqref{eq:intro4} (see Section \ref{sec:main}) and, introducing a suitable approximation scheme, we establish their global-in-time existence (see Theorem \ref{theo:weaksol} and Section \ref{sec:weaksol}), as long as we consider a class of sublinear and bounded latent heat $\ell$.
We remark that our notion of weak solution is ``entropic", meaning that the heat equation \eqref{eq:intro4} is replaced by the total energy inequality and an entropy production inequality.
Note that we use the boundedness of the latent heat $\ell$ to deduce further key properties of the (absolute) temperature, in particular its strict positivity (almost everywhere).
To prove the existence of such weak  solutions we rigorously develop a novel approximation scheme for \eqref{eq:intro1}-\eqref{eq:intro4}. Similar schemes are often omitted (or simply sketched) for similar systems in the literature. 
Finally, we prove the uniform strict positivity of the temperature (whenever the initial temperature is uniformly strictly positive) as a key property, adopting a De Giorgi iteration scheme (see Section \ref{Degiorgi}). To the best of the authors' knowledge, this approach, albeit very flexible (it has been successfully applied to many different problems, see, e.g., \cite{EPPS,GalPoia,GP, P}) has never been proposed in the literature for similar models. The closest result is available only in \cite{LRS}, obtained by means of a standard maximum principle in case of linear $Q$. Here we extend the validity of such a maximum principle to a more general class of functions $Q$ (see assumptions in Section \ref{sec:main}). 

The second part of our work is dedicated to the study of the sharp interface limit of the thermodynamically consistent phase-field model \eqref{eq:intro1}-\eqref{eq:intro4} by letting the diffuse interface width $\eps$ tend to zero (see Theorem \ref{theo:convergenceBV} and Section \ref{sec:sharp}). 
Our aim is to pass to the limit $\eps \to 0 $ in the entropic weak formulation of \eqref{eq:intro1}-\eqref{eq:intro4}, under a (time-integrated) energy converence hypothesis, in the spirit of the work by Luckhhaus and Sturzenhecker \cite{LuckhausStur}. Following this approach, we obtain in the limit $\eps \to 0$ (entropic) distributional (or $BV$) solutions to a non-isothermal Navier--Stokes/mean curvature flow, whose strong formulation reads as follows:
For a time dependent family of sets $ A(t)\subset \Omega$, $t \in [0,T)$, a velocity field $\ve= \ve(x,t)$ and a temperature field $\theta=\theta (x,t)$, 
 \begin{align}
	\partial_t \ve + \ve \cdot \nabla \ve &= - \nabla p + \opdiv (\mathbb{S}) - \mathbf{H}_A |\nabla \chi_A| && \text{ on } \Omega \times (0,T), \label{eq:intros1}\\
	\opdiv(\ve)&=0 && \text{ on } \Omega \times (0,T), \label{eq:intros2}\\
	V &= H_A -\ell(\theta)  &&\text{ on }  (\partial A \cap \Omega )\times (0,T), \label{eq:intros3}\\
	\frac{D Q(\theta)  }{Dt}+ \ell(\theta) V |\nabla \chi_A| &= \opdiv(\kappa(\theta) \nabla \theta) + \nu(\theta) |\nabla \ve + \nabla \ve^\mathsf{T}|^2 + |V|^2|\nabla \chi_A| && \text{ on } \Omega \times (0,T), \label{eq:intros4}
\end{align}
where $H_A = -\nabla \cdot{\n_{A}}$ ($\mathbf{H}_A = H_A  \n_A)$ denote the scalar (vectorial) mean curvature, $\n_A$ is the inner normal of $\partial A$, $V$ is the surface normal velocity (in the direction $\n_A$), and $|\nabla \chi_A|$ represents the surface measure of $A$. Similarly to \eqref{eq:intro1}-\eqref{eq:intro4}, this system has to be endowed with complete slip conditions for $\ve$,  ninety degree contact angle between $\partial A$ and $\partial \Omega$, and no (heat) flux condition for $\theta$.

Before commenting on our analysis, we recall some of the current available literature on the study of the sharp interface limit for evolutionary phase-field models.
Several results on the convergence of the Allen-Cahn type models to (mean) curvature driven flows are available
, ranging from formal asymptotic analysis to rigorous results for desparate notions of solution. 
For instance, for the convergence of the Allen-Cahn equation towards a smooth evolution we refer to the works \cite{DeMottoniSchatzman,Chen,BronsardR,AbelsMoser} and, more recently, to \cite{FischerLauxSimon,HenselMoser,FischerMarveggio,LauxKroemer} by means of the relative energy method, whereas in the framework of Brakke solutions we refer to \cite{Ilmanen}, in that of viscosity solutions to \cite{EvansSS}, and in that of De Giorgi varifold solution to \cite{HenselLaux-varifold,P1}. In the framework of the Navier--Stokes/Allen--Cahn system, convergence as $\eps \to 0$ to the corresponding smooth evolution was shown in the case of non-vanishing mobility, both via asymptotic analysis \cite{AbelsFei} and via relative energy method \cite{HenselLiu}, but also is the physical relevant case of vanishing mobility, both via asymptotic analysis \cite{AbelsFeiMoser} and via relative energy method \cite{AbelsFischerMoser}.
More importantly, in order to situate our result in the context of current research, we recall the (energy conditional) convergence result by Laux and Simon \cite{LauxSimon} in the setting of distributional (or BV) solutions to (multiphase) mean curvature flow. 
Similar results were then established to include contact angle conditions \cite{HenselLaux}, anisotropy \cite{LauxStinsonUllrich,EP}, volume-preserving contraint \cite{Laux}, and heterogeneity of the double well potential \cite{GanediMS} into the analysis.
To the best of our knowledge, only a few results are available on the analysis of the sharp interface limit of non-isothermal phase-field models.
We recalll the Caginalp model \cite{Caginalp1,Caginalp2}, consisting in the Allen-Cahn equation together with a simplier version of our heat equation \eqref{eq:intro4} (i.e., without quadratic terms on the right-hand side). The convergence of the Caginalp model to the Stefan problem as $\eps \to 0$ was rigorously proven via asymptotic expansions in \cite{CaginalpChen}. In the framework of weak solutions, convergence was shown for the quasi-stationary Caginalp model \cite{Schaetzle} and for the Caginalp model with kinetic undercooling \cite{Soner}.

A key originality of our work is therefore the study of the sharp interface limit of a novel non-isothermal phase-field model \eqref{eq:intro1}-\eqref{eq:intro4}, characterized by a thermodynamically consistent heat equation with (quadratic) dissipation terms \eqref{eq:intro4}. This motivates the entropic weak framework while also addressing the mathematical difficulties encountered during the analysis.
In order to pass to the sharp interface limit ($\eps \to 0$) in the (entropic) weak formulation of \eqref{eq:intro1}-\eqref{eq:intro4}, we consider a class of sublinear and bounded latent heat $\ell$, as we need to rely on the properties of (entropic) weak solutions. 
Indeed, the boundedness of $\ell$ is crucial to obtain the compactness of the phase indicator function
$\psi_\eps:= \int_{-1}^{\varphi_{ \eps}} \sqrt{2W(s)}\, \d s$
by means of a Modica-Mortola type trick for our entropic formulation.
Moreover, we need to show the boundedness of the phase-field variable $\varphi_{ \eps}$ and the strict separation of temperature $\theta_{\eps}$ from zero, holding uniformly in $\eps>0$. This latter property is nontrivial to establish as, in our entropic formulation with the use of the bounded latent heat $\ell$, the heat content $Q$ and the caloric part of the entropy have the same space integrability property. 
We obtain distributional (or $BV$) solutions to the entropic formulation of \eqref{eq:intros1}-\eqref{eq:intros4}, relying on the assumption that the time-integrated energies of the phase-field approximation converge to the perimeter (as in \cite{LauxSimon}), as well as on the strong convergence of the heat content, as $\eps \to 0$.
Furthermore, we observe that we are able to pass to the limit $\eps \to 0$ in the entropic counterpart of the quadratic terms in \eqref{eq:intro4} by using a lower-semicontinuity result by Ioffe \cite{Ioffe}. However in the passage to the limit we loose some information for the last term in \eqref{eq:intro4}. Indeed, this passes to its limit only as a (nonnegative) measure with a support containing the perimeter of $A$.
\linebreak

\paragraph{\bf Outline.} The rest of the paper is organized as follows.
We derive the non-isothermal phase-field model \eqref{eq:intro1}-\eqref{eq:intro4} in Section \ref{sec:derivation}.
We introduce the notation and the functional setting in Section \ref{sec:notation}.
We present the precise mathematical statement of our results in Section \ref{sec:main}, including the assumptions exploited throughout the paper.
The well-posedness of local strong solutions is shown in Section \ref{sec::strong}. We establish the global-in-time existence of (entropic) weak solutions in Section \ref{sec:weaksol}.
The convergence of weak solutions to a distributional (or $BV$) solution to a non-isothermal Navier--Stokes/mean curvature flow is proven in Section \ref{sec:sharp}.


 \section{Derivation of the phase-field model} \label{sec:derivation}
 In this section we derive the system of PDEs \eqref{eq:intro1}-\eqref{eq:intro4} by following the approach devised by Frémond \cite{Fremond} in the case of linear latent heat $\ell$.
 In case of more general forms of latent heat $\ell$, our systems of PDEs  \eqref{eq:intro1}-\eqref{eq:intro4} remains a valid model in the regime of temperatures at which $\ell$ is approximately linear (e.g., for temperatures below a critical one). Furthermore, we observe that the system \eqref{eq:intro1}-\eqref{eq:intro4} corresponds to a subclass of the general model derived in the work by H. Freist\"uhler and M. Kotschote \cite{FK}.

 \subsubsection*{Derivation of the model}
We provide a physical derivation of the model by following a variant of the general approach devised by Fr\'emond.
The starting point is the free energy density, which is given by 
\begin{equation}
	f_\eps  = \frac{\varepsilon}{2} |\nabla \varphi_{\varepsilon}|^2 + \frac1{\varepsilon}W( \varphi_{\varepsilon}) + g(\theta_\eps)-   \ell(\theta_\eps) \varphi_{\varepsilon}.
\end{equation}
Here, $g$ refers to pure heat conduction, is the main concave term of the free energy, and is linked to the specific heat $c_V (\theta) = Q^\prime(\theta)$ by the relation $Q(\theta) = g(\theta) -  \theta g^\prime(\theta)$, whereas $\ell $ represents the latent heat function, namely the difference in energy at fixed temperature between two pure phases (i.e., $\varphi_{\varepsilon} = \pm 1$). 
In particular, we consider a latent heat of the form $\ell(s)= \lambda  s$, for some constant $\lambda >0$.
The entropy density of the system is given by
\begin{align} \label{eq:defentropy}
	s_\eps=- \frac{\partial f_\eps}{ \partial\theta_\eps} =  - g'(\theta_\eps) +  \ell^\prime(\theta_\eps) \varphi_{\varepsilon},
\end{align}
and the internal energy density $e_\eps$ is given by the Gibbs relations, which reads as
\begin{align}\label{eq:internalenergy}
e_\eps = f_\eps + \theta_\eps s_\eps= \frac{\varepsilon}{2} |\nabla \varphi_{\varepsilon}|^2 + \frac1{\varepsilon}W(\varphi_{\varepsilon}) + Q(\theta_\eps) 
.
\end{align}

Following Frémond's approach (for more details we refer to Frémond's monograph \cite{Fremond}), we can write a system of the field equations representing the balance of momentum in absence of external forces according to Newton’s second law, the principle of virtual powers ruling the evolution of the order parameter $\varphi_{\varepsilon}$, and the balance of internal energy. More precisely, this system of 
 field equations reads as
\begin{align}
	&\partial_t \ve_\eps + \operatorname{div}(\ve_\eps \otimes \ve_\eps) = \operatorname{div}\sigma_\eps, \quad  \operatorname{div}\ve_\eps=0,  \label{eq:F1}\\
	& \operatorname{div}\He_\eps - B_\eps =0, \label{eq:F2} \\
	& \frac{De_\eps}{Dt} + \operatorname{div}\qe_\eps = \sigma_\eps : (\nabla \ve_\eps + \nabla \ve_\eps^\mathsf{T}) + B_\eps  \frac{D\varphi_{\varepsilon}}{Dt} + \He_\eps \cdot \nabla \frac{D\varphi_{\varepsilon}}{Dt} 
	,\label{eq:F3}
\end{align}
where $\sigma_\eps$ denotes the stress tensor, $\He_\eps$ is the energy flux, $B_\eps$ is the energy density, and $\qe_\eps$ is the heat flux. These are assumed to be the sum of their non-dissipative and dissipative components, namely $\sigma_\eps= \sigma_\eps^{nd}  + \sigma_\eps^d $, $\He_\eps= \He^{nd}_\eps + \He^d_\eps$, $B_\eps = B^{nd}_\eps + B^d_\eps$,  and  $\qe_\eps= \qe^{nd}_\eps+ \qe^{nd}_\eps$. In particular, we have
\begin{align}
	&B^{nd}_\eps= \frac{\partial f_\eps}{\partial \varphi_{\varepsilon}}= \frac1{\varepsilon} \partial_\varphi W(\varphi_{\varepsilon}) - \ell(\theta_\eps ) , \quad B^{d}_\eps= \varepsilon \frac{D\varphi_{\varepsilon} }{Dt},\label{eq:defB} \\
	&\He^{nd}_\eps=  \varepsilon \nabla \varphi_{\varepsilon}, \quad \He^d_\eps=0,\label{eq:defH}
\end{align}
and we set 
\begin{align} \label{eq:dissipdef}
	\qe_\eps^d= - \kappa(\theta_\eps) \nabla \theta_\eps, \quad 
	\sigma_\eps^{d} = \mathbb{S}(\theta_\eps) \ve_\eps - p \mathbb{I}, 
\end{align}
where $\mathbb{S}(\theta_\eps) \ve_\eps= \nu(\theta_\eps)(\nabla \ve_\eps + \nabla \ve_\eps^\mathsf{T}),$
whereas the nondissipative components $\qe^{nd}_\eps$ and $\sigma^{nd}_\eps$ are determined below so that the second law of thermodynamics is satisfied.
To this aim, we impose the validity of the second law of thermodynamics, namely
\begin{align*}
\theta_\eps \left( \frac{Ds_\eps }{Dt} + \operatorname{div}\Big(\frac{\qe_\eps }{\theta_\eps}\Big) \right) \geq 0.
\end{align*}
Then, we compute
\begin{align*}
	\frac{\partial f_\eps}{\partial \nabla \varphi_\eps} \cdot \frac{D\left(\nabla \varphi_\eps\right)}{D t}=\mathbf{H}^{n d}_\eps \cdot\left(\nabla \frac{D \varphi_\eps }{D t}-\nabla \mathbf{v}_\eps  \cdot \nabla \varphi_\eps  \right). 
\end{align*}
Having the equations above at our disposal, we obtain
\begin{align*}
	\theta_\eps \left( \frac{Ds_\eps }{Dt} + \operatorname{div}\Big(\frac{\qe_\eps }{\theta_\eps}\Big) \right) 
	&=
	\frac{D e_\eps }{D t}+\operatorname{div} \mathbf{q_\eps }-\frac{D f_\eps }{D t}-\frac{D \theta_\eps }{D t} s_\eps - \frac{\qe_\eps }{\theta_\eps } \cdot \nabla \theta_\eps 
	\\
&=	\frac{D e_\eps }{D t} +\operatorname{div} \mathbf{q_\eps }-\frac{\partial f_\eps }{\partial \varphi} \frac{D \varphi_\eps }{D t}-\frac{\partial f_\eps }{\partial \nabla \varphi_\eps } \cdot \frac{D \nabla \varphi_\eps }{D t}-\frac{\qe_\eps }{\theta_\eps } \cdot \nabla \theta_\eps 
	\\
& =  \sigma_\eps   : D\mathbf{v_\eps } +B_\eps  \frac{D \varphi_\eps }{D t}+\mathbf{H_\eps } \cdot \frac{D \nabla \varphi_\eps }{D t}
-\frac{\partial f_\eps }{\partial \nabla \varphi_\eps } \cdot \frac{\nabla D   \varphi_\eps }{D t} \\ 
&\quad  -B^{n d}_\eps  \frac{D \varphi_\eps }{D t}+\frac{\kappa(\theta_\eps )}{\theta_\eps }\left|\nabla \theta_\eps \right|^{2}- \frac{\mathbf{q}^{n d}_\eps  }{\theta_\eps }\cdot \nabla \theta_\eps 
	\\ 
& = \left(\sigma^{n d}_\eps +\varepsilon\left(\nabla \varphi_\eps  \otimes \nabla \varphi_\eps \right)\right): D\mathbf{v}_\eps +\nu(\theta)|D\mathbf{v}_\eps |^{2}
\\
 & \quad  + \eps \left| \frac{D\varphi_\eps }{Dt} \right|^2  
 + \frac{\kappa(\theta_\eps )}{\theta_\eps } \left|\nabla \theta_\eps \right|^{2} - \frac{\mathbf{q}^{n d}_\eps  }{\theta_\eps }\cdot \nabla \theta_\eps .
\end{align*}
To obtain non-positivity, we choose:
\begin{align}
	\sigma^{nd}_\eps  = - \varepsilon \nabla \varphi_{\varepsilon} \otimes \nabla \varphi_{\varepsilon} , \quad \qe^{nd}_\eps =0. \label{eq:nondissipdef}
\end{align}

 \subsubsection*{Strong formulation of the model}
From the derivation sketched above, we obtain the system of PDEs representing the
strong formulation of our model.
In particular, from \eqref{eq:F1}, \eqref{eq:dissipdef}, and \eqref{eq:nondissipdef} we deduce that the evolution of the velocity $\ve_\varepsilon$ is ruled by the Navier-Stokes system, namely
 \begin{align}
\text{conservation of momentum: }&&	\partial_t \ve_\varepsilon + \ve_\varepsilon \cdot \nabla \ve_\varepsilon &= - \nabla p + \opdiv (\mathbb{S}(\theta_\varepsilon)\mathbf v_\varepsilon) - \opdiv(\varepsilon \nabla \varphi_{\varepsilon}  \otimes \nabla \varphi_{\varepsilon} )  ,\label{eq:strongmom}\\
\text{incompressibility: }	&&	\opdiv(\ve_\varepsilon)&=0,\label{eq:intro2b}
\end{align}
 where $\mathbb{S}(\theta_\eps) \ve_\eps= \nu(\theta_\eps)(\nabla \ve_\eps + \nabla \ve^\mathsf{T}_\eps)$.
The evolution of the order parameter $\varphi_{\varepsilon}$ is prescribed by \eqref{eq:F2}, which using \eqref{eq:defB}-\eqref{eq:defH} reads as
	 \begin{align}
\text{Allen-Cahn equation: } \quad 	\frac{D \varphi_{\varepsilon} }{Dt} &= \Delta \varphi_\varepsilon - \frac{1}{\varepsilon^2} \partial_\varphi W(\varphi_\varepsilon) + \frac{1}{\varepsilon}\ell(\theta_\varepsilon). \label{eq:strongAC}
\end{align}
Using \eqref{eq:dissipdef},\eqref{eq:defB}, \eqref{eq:defH}, and \eqref{eq:nondissipdef}, as well as the explicit form \eqref{eq:internalenergy}, the internal energy balance \eqref{eq:F3} can be rewritten as 
\begin{align}
\text{heat equation: } \quad 	\frac{D Q(\theta_\varepsilon )  }{Dt}+ \ell(\theta_\varepsilon) \frac{D \varphi_{\varepsilon}  }{Dt} &= \opdiv(\kappa(\theta_\varepsilon) \nabla \theta_\varepsilon) + \nu(\theta_\varepsilon) |\nabla \ve_\varepsilon + \nabla \ve_\varepsilon^\mathsf{T}|^2 + \varepsilon \left| \frac{D \varphi_{\varepsilon} }{Dt}\right|^2. \label{eq:strongheat}
\end{align}
 
 
 \subsubsection*{Balances of total energy and entropy}
 A key point in the weak formulation of our model consists in replacing the heat equation \eqref{eq:strongheat} with the balances of total energy and of entropy,
which are mathematically more tractable, but still encode the main features of the model, i.e., the first and second laws of thermodynamics.
 
 From the equations above, one can derive the total energy balance
 \begin{align}
 	&\partial_t \Big(e_\varepsilon + \frac12 |\ve_\eps|^2  \Big) + \ve_\eps \cdot \nabla \Big(e_\varepsilon + \frac12 |\ve_\eps|^2  \Big) \notag \\
 	&= \operatorname{div}(\kappa(\theta_\varepsilon) \nabla \theta_{\varepsilon}) +\opdiv ((\mathbb{S}(\theta_\eps) \ve_\eps )\ve_\eps - p \ve_\eps )  - \opdiv(\varepsilon(\nabla \varphi_\eps \otimes \nabla \varphi_\eps ) \ve_\eps ) 
 	.\label{eq:energybal}
\end{align} 	
 By multiplying \eqref{eq:strongheat} by ${1}/{\ell(\theta_\eps)}$, one obtains the (modified) entropy balance
 \begin{align}
 	&\partial_t(\Lambda (\theta_{\varepsilon}) + \varphi_\varepsilon) + \ve_\eps \cdot \nabla (\Lambda (\theta_{\varepsilon})  + \varphi_\varepsilon) \notag \\
 	&= \operatorname{div}
 	\Big( \frac{\kappa(\theta_\varepsilon)}{\ell(\theta_\varepsilon)}  \nabla \theta_{\varepsilon} \Big) + \frac{\ell^\prime (\theta_\varepsilon)}{\ell^2(\theta_\varepsilon)} \kappa(\theta_\varepsilon) |\nabla \theta_{\varepsilon}|^2  + \frac{\nu(\theta_\varepsilon)}{\ell(\theta_\varepsilon)} |\nabla \ve_\varepsilon + \nabla \ve_\varepsilon^\mathsf{T}|^2+ \frac{\varepsilon }{\ell(\theta_\varepsilon)} |(\partial_t + \ve _\eps\cdot \nabla)\varphi_\varepsilon|^2,  \label{eq:entropybal}
 \end{align}
where
 \begin{equation}
 	 \Lambda(\theta_\varepsilon):= \int_0^{\theta_\varepsilon} \frac{Q^\prime(s)}{\ell(s)} \, ds.  \label{eq:defmodentropy}
 \end{equation}
We note that, under assumption of linear latent heat $\ell(s)=\lambda s$, the function $\lambda$ defined in \eqref{eq:defentropy} coincides with the purely caloric part of the entropy density (i.e., first term in \eqref{eq:defentropy}). 
Moreover, at this level, the (modified) entropy balance \eqref{eq:entropybal} is an equality, however it will turn into an inequality in the framework
weak solutions. This is due to the behaviour of the quadratic terms on the right hand side of \eqref{eq:entropybal} with respect to weak limits.
In case we could prove that there exist a smooth solution to the weak formulation of the model, then we would recover \eqref{eq:entropybal} as an equality and
\eqref{eq:entropybal} would be equivalent to \eqref{eq:strongheat}.


 \subsubsection*{Boundary and initial conditions}  \label{subsec:bicond}
In accordance with the physical derivation, we assume that the system is insulated from the exterior.
Consistently, we impose complete slip boundary conditions:
 \begin{align} \label{eq:BCv}
 	\ve_\eps \cdot \no= 0 \quad  \text{ and } \quad { (\operatorname{Id} - \no \otimes \no) ((\mathbb{S}(\theta_\eps) \ve_\eps)\no )  = \mathbf{0}}  \quad \text{ along } \partial \Omega. 
 \end{align}
The first condition prescribes that the normal component of the velocity is zero along the boundary, i.e., the fluid cannot exit from the domain $\Omega$ and it can move only tangentially to the boundary.
The second condition excludes external contribution to the viscous stress, hence there are no friction effects with the boundary. 
Note that one could impose the Navier boundary conditions with general friction parameter; however, for simplicity, we consider the complete slip case.
 For the order parameter $\varphi_{\varepsilon}$ we assume a contact angle of $90^\circ$ between the diffuse interface and the boundary of the domain, namely 
 \begin{align}\label{eq:BCvarphi}
 	\nabla \varphi_{\varepsilon} \cdot \no =0  \quad \text{ along } \partial \Omega. 
 \end{align}
Moreover, we consider no heat flux (through the boundary) condition
 \begin{align}\label{eq:BCtheta}
 	\qe_\eps \cdot \no = - \kappa(\theta_\eps) 	 \nabla \theta_\eps \cdot \no =0 \quad \text{ along } \partial \Omega. 
 \end{align}
 At last, the system of PDEs is endowed with the initial conditions
 \begin{equation} 
 	\ve_{\eps}(\cdot, 0) = \ve_{ 0 }, \quad  \varphi_{\varepsilon} (\cdot, 0) = \varphi_{ \eps, 0 }, \quad \theta_{\varepsilon} (\cdot, 0) = \theta_{0 } \quad \text{ a.e. on } \Omega. 
 \end{equation}
 


 \section{Notation and functional setting} \label{sec:notation}
 We start by introducing some notation.
 
 \begin{enumerate}[label=\textnormal{(N\arabic*)},leftmargin=*]
 	
 	\item \textbf{Notation for general Banach spaces.} 
 	For any normed space $X$
    , we denote its norm by $\|\cdot\|_X$,
 	its {dual space by $X'$}.
 	Moreover, if $X$ is a Hilbert space, we write $(\cdot,\cdot)_X$ to denote the corresponding inner product.
 	Furthermore, for any vector space $X$, corresponding spaces of vector-valued or matrix-valued functions with each component in $X$ are denoted by $\mathbf{X}$.
 	
 	\item \textbf{Lebesgue and Sobolev spaces.} 
 	Assume $\Omega$ to be a sufficiently smooth bounded domain in $\R^d$, $d=2,3$.
 	For $1 \leq p \leq \infty$ and $k \in \N$, the classical Lebesgue and Sobolev spaces defined on $\Omega$ are denoted by $L^p(\Omega)$ and $W^{k,p}(\Omega)$, and their standard norms are denoted by $\|\cdot\|_{L^p(\Omega)}$ and $\|\cdot\|_{W^{k,p}(\Omega)}$, respectively.
 	In the case $p = 2$, the notation $H^k(\Omega) = W^{k,2}(\Omega)$ is used. We point out that $H^0(\Omega)$ coincides with $L^2(\Omega)$. We also set $W^{2,s}_N(\Omega):=\{f\in W^{2,s}(\Omega): \partial_{\n_{\partial\Omega}}f=0\quad\text{ on }\partial\Omega\}$, $s\geq2$. 
 	
 	Also, for any interval $I\subset\R$, any Banach space $X$, $1 \leq p \leq \infty$ and $k \in \N$, we write $L^p(I;X)$, $W^{k,p}(I;X)$ and $H^{k}(I;X) = W^{k,2}(I;X)$ to denote the Lebesgue and Sobolev spaces of functions with values in $X$. The canonical norms are denoted by $\|\cdot\|_{L^p(I;X)}$, $\|\cdot\|_{W^{k,p}(I;X)}$ and $\|\cdot\|_{H^k(I;X)}$, respectively. 
 	
 	\item \textbf{Spaces of continuous functions.}
 	For any interval $I\subset\R$ and any Banach space $X$, $C(I;X)$ ($C_{w}(I;X)$) denotes the space of continuous functions (weakly continuous functions, respectively)  mapping from $I$ to $X$. Then, we denote by $C^\gamma(I;X)$, $\gamma\in(0,1]$, the space of $\gamma$-H\"{o}lder (respectively, Lipschitz) continuous functions with values in $X$, whereas by $C_c^k(I;X)$ we denote the spaces of $k$-continuously differentiable functions with compact support mapping $I$ into $X$.  
 	
 	\item \textbf{Radon measures, set of finite perimeters, and functions of bounded variation.}
 Assume $\Omega$ as before and consider any interval $I\subset\R$.
 We denote by $\mathcal{M}(\Omega \times I)$ ($\mathcal{M}_+(\Omega \times I )$) the set of (signed) Radon measures (positive Radon measures, respectively) and by $\langle \cdot, \cdot  \rangle$ the corresponding duality pairing. We introduce $BV(\Omega)$ as the set of functions of bounded variation on $\Omega$.  
 For a set $A \subset \R^d$, $\chi_A$ is the associated characteristic function, taking the value $1$ on $ A$ and $0$ in the complement.
 We say $A$ is a set of finite perimeter in $\Omega$ if $\chi_A\in BV(\Omega;\{0,1\}),$ and we denote its distributional derivative in $\Omega$ by $\nabla \chi_A$ and its total variation measure by $|\nabla \chi_A|$. We use $n_A:=\frac{\nabla \chi_A}{|\nabla \chi_A|}$ for the measure-theoretic inner normal of $A$. Note that we denoted by $\frac{\nabla \chi_A}{|\nabla \chi_A|}$ the Radon-Nikodym derivative of $\nabla\chi_A$ with respect to its total variation. Using $\partial^\ast A$ as the reduced boundary of the set $A$ in $\Omega$, we have $\int_\Omega \cdot \, \d |\nabla \chi_A| =  \int_{\partial^\ast A} \cdot \, \d \mathcal{H}^{d-1}$, where $\mathcal{H}^{d-1}$ is the Hausdorff measure. For any one-parameter family of open sets, $A(t) \subset \Omega$ , $t \in I$, we define $\chi_A(\cdot, t):=\chi_{A(t)}(\cdot)$ for any $t \in I$.
 For more information on functions of bounded variation and sets of finite perimeter, we refer to the classical \cite{AFP}.
 	
 	\item \textbf{Spaces of divergence-free functions.}
 We define the closed linear subspaces
 	\begin{align*}
 		\LL^p_\sigma(\Omega)
 		&:=\overline{\{\mathbf{u}\in \mathbf{C}^\infty_0(\Omega) \,\big\vert\, \operatorname{div}\ \mathbf{u}=0\}}^{\mathbf{L}^p(\Omega)}
 		\subset \LL^p(\Omega), \quad p\in[2,\infty),\\
 		\mathbf V_\sigma
 		&:= \LL^2_\sigma(\Omega) \cap \mathbf H^1(\Omega),
 		\\
 		\mathbf W^{k,p}_\sigma(\Omega)
 		&:= \mathbf W^{k,p}(\Omega)\cap \mathbf V_\sigma,\quad p\in[2,\infty),\quad k\in \mathbb N.
 	\end{align*}
 We then set for simplicity $\mathbf H_\sigma:=\mathbf L^2_\sigma(\Omega)$.
 Moreover, we denote 
 $ W^{k,p}_{N} $ and $\W^{k,p}_{\sigma,N} $
 the spaces $ W^{k,p} $ and $\W^{k,p}_{\sigma} $ endowed with the Neumann boundary conditions  
 $\partial_{\mathbf n_{\partial\Omega}} u =0$
 and
 $ (\operatorname{Id}-\mathbf n_{\partial\Omega}\otimes \mathbf n_{\partial\Omega})(\nabla \mathbf u+\nabla \mathbf u^T)\mathbf n_{\partial \Omega}=\mathbf 0$ a.e. on $\partial\Omega$, respectively.
 
 \end{enumerate} 
 
  \section{Main results} \label{sec:main}
  In this section we first introduce the assumptions on our model \eqref{eq:intro1}-\eqref{eq:intro4}, then we provide the statements of our main results.

   \subsection{Assumptions on the coefficients} \label{subsec:hp}
  We first recall that $W(u)$ is a standard double-well potential of the form
  \begin{align}
  	W(s) = \frac14 (s^2-1)^2, \quad \text{ whence  } \quad W^\prime(s)= s^3-s.
  \end{align}
  
  We assume the specific heat $c_V$, and thus the thermal component $Q$ in the energy density \eqref{eq:internalenergy}, to be given by 
  \begin{equation}\label{eq:defQ}
  	c_V(s)= s^{\alpha}, \quad 
  	Q(s)= \int_0^s c_V(y) \, \d y = \frac{1}{1+\alpha}s ^{1+\alpha}, \quad \text{ for }  0 < \alpha \leq 1,
  \end{equation}
  whereas the heat conductivity to be of the form
  \begin{equation}\label{eq:defkappa}
  	\kappa(s)= \kappa_1 + \kappa_2 s^\beta, \quad \kappa_1, \kappa_2>0, \quad \text{ for } 	\beta >0 .
  \end{equation}
  We observe that the bound $\alpha \leq 1$  implies that the thermal component $\Lambda$ of the entropy is concave, whereas the power-like behavior for the heat conductivity is typical of several type of fluids. Additional technical limitations on the exponents $\alpha$ and $\beta$ will be introduced in the statements below.

  We assume  the viscosity of the mixture $\nu\in W^{1,\infty}(\mathbb R)$ (depending on the temperature) to be bounded from above and below, i.e., 
  \begin{equation}
  	0 < \nu_1 \leq \nu(s) \leq \nu_2 ,\quad \nu_1, \nu_2>0. \label{hp:nu}
  \end{equation}
  Notice that the dependence of the viscosity on the temperature $\theta$, as well as its Lipschitz regularity, is often assumed in the literature, even in the case of non-Newtonian fluids (see, e.g., \cite{GPPV}).
  
  We then consider two classes of latent heat functions $\ell: [0, \infty) \rightarrow [0, \infty) $, $\ell \in C^2[0, \infty)$, namely 
  \begin{itemize}
  	\item linear latent heat functions:
  	\begin{align} \label{hp:linearell}
  		\ell(s)= \lambda s, \quad \lambda >0,
  \end{align}
  \item sublinear and bounded latent heat functions satisfying: 
  \begin{enumerate}[label=(\alph*)]
  \item $\ell\in C^1([0,\infty))$;
  	\item $\ell(0)=0$, $0 \leq \ell(s) \leq L$ for all $s \in [0, \infty)$ with $L>0$; 
  \item $\ell'$ is monotone nonincreasing on $[0, \infty)$; 
  \item $\lambda{s} \geq \ell(s)$ for some $\lambda>0$ for all $s \in [0, \infty)$; 
  \item $ \ell(s) \geq \delta_0 \ell^\prime(0) s$  for all $s \in [0, 1)$, 
  for 
   some constant $0< \delta_0 <1$ ; 
  	\item  $({1+s^{\gamma}})^{-1} \leq C_\ell {\ell^\prime (s)} $ for all $s \in [1, \infty)$,  for some constant $C_\ell>0$ and exponent $0<\gamma\leq 2$ .
  \end{enumerate}
An example of  sublinear latent heat function is given by $\ell(s)=\arctan s$, which satisfies the properties above. 
We observe that condition (e) {entails that $\ell'(0)>0$ as well as that $\ell'$, which is monotone nonincreasing by condition (b), can converge to $0$ at infinity with a sufficiently slow rate}. This last condition is introduced mainly because of mathematical reasons (see Section 5.1.2).
\end{itemize}

  As a consequence of the assumptions above, the thermal part $\Lambda $ of the (modified) entropy density is given by  
  \begin{equation}  \label{eq:defLambda}
  \Lambda(s)= \int_0^{s} \frac{Q^\prime(y)}{\ell(y)} \, \d y = \int_0^{s}  \frac{y^{\alpha}}{\ell(y)} \, \d y ,
  \end{equation}
  which, in case of linear latent heat $\ell(s)= \lambda s$, reduces to 
  $$
  \Lambda(s) = \frac{1}{\alpha\lambda }s^\alpha .$$

  \subsection{Well-posedness of local strong solutions}
  \label{WPlocal}
  Our first main result consists in showing that there exists a unique local-in-time strong solution satisfying \eqref{eq:strongmom}-\eqref{eq:strongheat}, given sufficiently regular initial data.

  To this aim, we need to introduce a suitable functional setting in the spirit of \cite{AGP1,AWe}, in order to apply a contraction mapping principle. 
 In particular, we first define the spaces 
 {\begin{align*}
 	&	Z_{\tT}^3:=\{f\in L^r(0,\tT;W^{2,r}(\Omega))\cap W^{1,r}(0,\tT;L^r(\Omega)):\ \partial_{\mathbf n_{\partial\Omega}} f =0\quad\text{a.e. on }\partial\Omega\},\quad\text{ for some }   r\in(\tfrac{7}2,+\infty),
 	\\&Z^1_{\tT}:=\{\mathbf f\in L^p(0,\tT; \W^{2,p}_\sigma(\Omega))\cap W^{1,p}(0,\tT;\L^p_\sigma(\Omega)):\  
    (\operatorname{Id}-\mathbf n_{\partial\Omega}\otimes \mathbf n_{\partial\Omega})(\nabla \mathbf f+\nabla \mathbf f^T)\mathbf n_{\partial \Omega}=\mathbf 0 \quad \text{a.e. on }\partial\Omega\} ,\\&\qquad\quad \text{for some }p\in (r,p^*),\\&
 	Z_{\tT}^2:=\{g\in L^q(0,\tT;W^{2,q}(\Omega))\cap W^{1,q}(0,\tT;L^q(\Omega)):\ \ \partial_{\mathbf n_{\partial\Omega}} g =0\quad\text{a.e. on }\partial\Omega\},\\&\qquad\quad\text{ for some } q\in(2p,+\infty),
 \end{align*}}
 where $p^*=\frac{3}2\frac{r}{5-r}>r$ for $r\in(\frac72,5)$, $p^*=+\infty$ for $r\geq 5$. Then we introduce, for a given triple $$(\v_0,\vphi_{\eps,0},\theta_0)\in ( {\W^{2,p}_{\sigma,N}(\Omega)},\L^p_\sigma(\Omega))_{1-\frac1p,p}\times (W^{2,q}_N(\Omega),L^q(\Omega))_{1-\frac1q,q}\times (W^{2,r}_N(\Omega),L^r(\Omega))_{1-\frac1r,r},$$ the metric spaces 
  $$
  X^1_{{T}}:=\{\v\in \Zut:\ \v(\cdot, 0)=\v_0\},\quad X^2_{{T}}:=\{\vphi\in \Zdt:\ \vphi(\cdot, 0)=\vphi_{\eps,0}\},\quad  X^3_{{T}}:=\{\theta\in \Ztt:\ \theta(\cdot, 0) =\theta_0\},
  $$
  $$
  \XT:= \XT^1\times \XT^2\times \XT^3,
  $$
  endowed with the corresponding $Z_T^i$ norms, $i=1,2,3$, defined in Section \ref{sec::strong}. We thus have the following

  \begin{teor}[Existence and uniqueness of a local-in-time strong solution]

  	\label{thm:strong}
  	Fix $\varepsilon>0$. Let the initial data be such that \begin{equation}
  	(\v_0,\vphi_{\eps,0},\theta_0)\in (  \W^{2,p}_{\sigma,N}(\Omega),\L^p_\sigma(\Omega))_{1-\frac1p,p}\times (W^{2,q}_N(\Omega),L^q(\Omega))_{1-\frac1q,q}\times (W^{2,r}_N(\Omega),L^r(\Omega))_{1-\frac1r,r}, \label{eq:reginitialstrong}
  \end{equation}
  	 $\norm{\vphi_{\eps,0}}_{L^\infty(\Omega)}\leq 1$ and $\inf_{x\in \Omega}\theta_0\geq c>0$.
  	
  	Assume $\alpha\in (0,1]$ and $\beta\geq1$. 
  	Under the assumptions stated in Section \ref{subsec:hp}, there exists $\tT_\varepsilon>0$, possibly depending on $\eps$, such that non-isothermal Navier-Stokes/Allen-Cahn system
  	\eqref{eq:strongmom}-\eqref{eq:strongheat} with boundary conditions \eqref{eq:BCv}-\eqref{eq:BCtheta} admits a unique solution $(\v_\eps,\vphi_\eps,\theta_\eps)$ such that 
  	$$
  	(\v_\eps,\vphi_\eps,\theta_\eps) \in \XT^1\times \XT^2\times \XT^3,
  	$$
  	with $$
  	\theta_\eps(x,t)\geq \frac c2>0,\quad \forall (x,t)\in\overline{\Omega}\times [0,T_\eps].
  	$$
  \end{teor}
  \begin{osse}
Note that, since we have, by standard embedding and interpolation results, $\vphi_\eps\in C(\overline{\Omega}\times[0,T])$, if we assume, as in the physical case, $\vphi_{\eps,0}\in[-1,1]$, then, by continuity, for any $\delta>0$ there exists $0<T_\delta\leq T$, such that 
$$
\vert \vphi_\eps(x,t)\vert\leq 1+\delta,\quad \forall (x,t)\in \overline{\Omega}\times[0,T_\delta].
$$
  \end{osse}

\subsection{Existence of global-in-time weak solutions}
 Our second main result consists in establishing global-in-time existence of solutions to the initial-boundary value problem associated to the weak formulation of the non-isothermal Navier-Stokes/Allen-Cahn system \eqref{eq:strongmom}-\eqref{eq:strongheat}. This will be proved by introducing a suitable approximation scheme and showing weak sequential stability of families of approximating solutions. 

First, we introduce a notion of weak solution for the model.
\begin{defi}[Weak solution to the non-isothermal Navier-Stokes/Allen-Cahn system] \label{def:weaksol}
	Fix $\eps>0$. Let $T >0$ be a given final time. 
	Consider an initial phase field $\vphi_{\eps,0} \in H^1(\Omega)$, an initial velocity field $\ve_{0} \in \Hs$, and an initial temperature $\theta_{0} \in L^{1+\alpha}(\Omega)$ such that $\theta_0 \geq 0$. 
	A weak solution to the non-isothermal Navier-Stokes/Allen-Cahn model \eqref{eq:strongmom}-\eqref{eq:strongheat} is a quadruplet of measurable functions $(\ve_\eps,p_\eps,  \varphi_\eps, \theta_\eps): \Omega \times [0, T) \rightarrow \R^d \times \R \times \R  \times \R $ 
	belonging to the following regularity class: 
	\begin{align}
		&\label{a0}\ve_\eps \in W^{1,\frac{4}{3}}(0,T;\Vs')\cap  L^\infty(0,T; \Hs)	\cap L^2(0,T; \Vs),\\& p_\eps\in L^\frac53(0,T;L^\frac53(\Omega)), \\
		&	\varphi_\eps \in W^{1, \frac53}(0,T; L^\frac53(\Omega))	 \cap  L^\infty(0,T; H^1(\Omega))	\cap L^2(0,T; H^2(\Omega)), \\
		&	\theta_\eps \in  L^\infty(0,T; L^{1+\alpha}(\Omega))	\cap  L^\beta(0,T; L^{3\beta}(\Omega)) \cap  L^2(0,T; H^1(\Omega)),\\&
		{\int_\Omega (\Lambda(\theta_\eps(\cdot, t)) + \varphi_\eps(\cdot, t)) 	\, \dx\geq \int_\Omega  (\Lambda(\theta_0) + \varphi_{\eps,0})\, \dx ,\quad \forall t\in(0,T),}\label{extra}\\&
		\label{thetagen}\theta_\eps\geq 0\textit{ a.e. in } \Omega \times (0,T),
	\end{align}
and satisfying
	the weak momentum balance
	\begin{align} \label{eq:weakbalmom}
		&\int_{0}^T \int_{\Omega} (\ve_\eps \cdot \partial_t \xe + (\ve_\eps  \otimes \ve_\eps  ): \nabla \xe  
		) \, \dx \dt + \int_\Omega \ve_0 \cdot \xe(\cdot, 0) \, \dx\notag \\
		&=	\int_{0}^T \int_{\Omega} \mathbb{S}(\theta_\eps)\ve_\eps: \nabla \xe \, \dx \dt+\int_{0}^T \int_{\Omega} p_\eps\opdiv\xe \, \dx \dt- 	\int_{0}^T \int_{\Omega} \eps (\nabla \varphi_\eps  \otimes \nabla \varphi_\eps ): \nabla \xe \, \dx \dt ,
	\end{align}
	{where $\mathbb{S}(\theta_\eps)\ve_\eps= \nu(\theta_\eps)(\nabla \ve_\eps + \nabla \ve^\mathsf{T}_\eps)$,} for all $\xe \in C^\infty_c([0,T); C^1(\overline{\Omega}; \mathbb{R}^3 ))$ such that $\xe \cdot \no=0$ on $\partial \Omega \times [0,T)$, 
    the Allen-Cahn equation
	\begin{align}\label{eq:weakAC}
		\partial_t \varphi_\eps + \ve_\eps \cdot \nabla \varphi_\eps  =  \Delta \varphi_\eps  -  \frac{1}{\eps^2} \partial_\varphi W(\varphi_\eps ) + \frac{1}{\eps} \ell(\theta_\eps)  \quad \text{a.e. in } \Omega \times (0,T),
	\end{align}
the total energy inequality
\begin{align}\label{ineq:enconsl}
	E_{\operatorname{tot}, \eps}[\vphi_\eps(\cdot, T'), \ve_{\eps}(\cdot, T'), \theta_{\eps}(\cdot, T')]\leq  E_{\operatorname{tot}, \eps}[\vphi_{\eps,0}, \ve_{0}, \theta_{0}] , \quad \text{for a.a. } T^\prime \in [0, T),
\end{align}		
where
\begin{align*}&E_{\operatorname{tot}, \varepsilon}[\vphi_\eps(\cdot, T'), \ve_{\eps}(\cdot, T'), \theta_{\eps}(\cdot, T')]:= \int_{\Omega} e_\varepsilon(\cdot, T') \, \dx + \int_{\Omega} \frac12 |\ve_\eps(\cdot, T')|^2\, \dx ,
\end{align*}	
and the weak entropy balance \begin{align} \label{eq:weakent}
  		&\int_{0}^T \int_\Omega(\Lambda(\theta_\eps) + \varphi_\eps) \partial_t \zeta \, \dx \dt
  		+ \int_{0}^T \int_\Omega (\Lambda(\theta_\eps)  + \varphi_\eps) \ve_\eps \cdot \nabla \zeta \, \dx \dt
  		+ \int_{0}^T \int_\Omega h(\theta_\eps ) \Delta \zeta \, \dx \dt
  		\notag \\ 
  		&\leq - \int_{0}^T \int_\Omega \bigg( \frac{\nu(\theta_\eps)}{\ell(\theta_\eps)} |\nabla \ve_\eps + \nabla \ve_\eps^\mathsf{T}|^2 + \frac{\ell^\prime (\theta_\eps)}{\ell^2(\theta_\eps)}  {\kappa(\theta_\eps)}  |\nabla \theta_\eps|^2
  		+ \frac{\varepsilon}{\ell(\theta_\eps)} \Big| \frac{D\varphi_\eps}{Dt} \Big|^2 \bigg) \zeta\, \dx \dt \notag  \\
  		& \quad - \int_{\Omega}(\Lambda(\theta_0) +  \varphi_{\eps,0}) \zeta( \cdot, 0) \, \dx ,
  	\end{align}
  	where $h(\theta_\eps)= \int_{1}^{\theta_\eps} \frac{\kappa(s)}{\ell(s)} \, \d s $ 
  	for all $\zeta \in C^\infty_c([0,T); C^2(\overline{\Omega}))$ such that $\zeta \geq 0$ and $\nabla \zeta  \cdot \no=0$ on $\partial \Omega \times [0,T)$.
\end{defi}

 \begin{osse}
	The regularity class stated in Definition \ref{def:weaksol} refers to the three-dimensional case $d=3$.
	We observe that better regularity for weak solutions can be established in the two-dimensional case $d=2$.
\end{osse}

In this setting, our main result reads as follows
   \begin{teor}[Existence of global-in-time weak solutions]
   	\label{theo:weaksol}
   Let $\tfrac12 < \alpha \leq 1$ and $\beta \geq 2$. 
  Assume that $\vphi_{\eps,0} \in H^1(\Omega)$, $\ve_{0} \in \Hs$, $\theta_{0} \in L^{1+\alpha}(\Omega)$ such that $\theta_0 \geq 0$, and
  $$
  |\vphi_{\eps,0}| \leq 1  
  \quad \text{ a.e. in } \Omega ,
  $$	
  and assume that there exists $\tau_0>0$ such that
 \begin{align}
 	\frac{\int_\Omega \Lambda(\theta_0)\, \dx}{\vert\Omega\vert}\geq 2+\tau_0.
 	\label{assumption}
 \end{align}
 Under the assumptions stated in Section \ref{subsec:hp}, for the class of bounded sublinear latent heat functions $\ell$, the non-isothermal Navier-Stokes/Allen-Cahn model admits at least one weak solution, {defined on $(0,+\infty)$}, in the sense of Definition \ref{def:weaksol} for any $T>0$, satisfying the following properties:
   \begin{enumerate}
  		\item[1)] (Upper and lower bound) The $L^\infty$-bound of $\vphi$ is conserved by the flow in the sense that
  		for any $\tau>0$ there exists {$\eps_\tau=\frac{W'(1+\tau)}{L}$}, such that we have, for any $T>0$,
  		\begin{align}\label{eqn:boundsACl}
  			&	-1 \leq \vphi_{\eps} \leq 1+ \tau  \text{ a.e. in } \Omega \times(0,T),\quad \forall\ \eps\in(0,\eps_\tau).
  		\end{align}
  		\item[2)](Strict positivity of $\theta_\eps$) 
 There exists $\eps_{\tau_0}>0$ such that $\theta_\eps>0$ almost everywhere in $\Omega\times(0,T)$ for any $\eps\in(0,\eps_{\tau_0})$ and any $T>0$.
 \item[3)] (Localized total energy inequality)
	\begin{align}\label{eq:weaktoten}
		&\int_{0}^T \int_{\Omega}  \Big(e_\eps + \frac12 |\ve_\eps|^2\Big) \partial_t \eta \, \dx \dt
		+	\int_{0}^T \int_{\Omega}  \Big(e_\eps \ve_\eps + \frac12 |\ve_\eps|^2\ve_\eps \Big)\cdot \nabla \eta  \, \dx \dt
		\notag\\&
		+ \int_{0}^T \int_{\Omega} \hat\kappa(\theta_\eps  ) \Delta \eta \, \dx \dt + 	\int_{0}^T \int_{\Omega} p_\eps \ve_\eps  \cdot \nabla \eta \, \dx \dt
		- \int_{0}^T \int_{\Omega} (\mathbb{S}(\theta_\eps)\ve_\eps ):  (\ve_\eps \otimes \nabla \eta ) \, \dx \dt\notag\\& 
		+  \int_{\Omega}  \Big(e_0 + \frac12 |\ve_0|^2\Big) \eta(\cdot, 0)\, \dx \geq 0,
	\end{align}
	where $\hat \kappa (\theta_\eps)= \int_0^{\theta_\eps} \kappa(s) \d s = \kappa_1\theta_\eps  + \kappa_2 \frac{1}{\beta +1} \theta_\eps^{\beta +1}$,
 for all $\eta \in C^\infty_c([0,T); C^2(\overline{\Omega}))$ such that $\eta \geq 0$ and $\nabla \eta  \cdot \no=0$ on $\partial \Omega \times [0,T)$.
 \end{enumerate}

  \end{teor}

 {\begin{osse}
 		Notice that, since by standard results, $\ve_\eps\in C_{w}([0,T],\Hs)$, and $\varphi_\eps\in C([0,T];H^1(\Omega))$, then we have $\ve_\eps(0)=\ve_0$, $\varphi_\eps(0)=\varphi_{\eps,0}$, where the limit is intended in $\Hs$ weak and $H^1(\Omega)$, respectively. Due to a lack of control on the time derivative of $\theta_\eps$, the initial datum $\theta_0$ is only enforced weakly in the sense of 
        \eqref{eq:weakent}.
 		\label{initialdata}
 \end{osse}} 
 
 \begin{coro}\label{lemma:furtherproperties}
  Any weak solution to the non-isothermal Navier-Stokes/Allen-Cahn model \eqref{eq:strongmom}-\eqref{eq:strongheat} given by Theorem \ref{theo:weaksol} satisfies the following property:
  	\begin{enumerate}
\item[3)] (Entropy production inequality)  
For any $T>0$, defining the entropy as
\begin{align*}
    S[\varphi_\varepsilon(\cdot, T'), \theta_\varepsilon(\cdot, T')] := \int_{\Omega} \Lambda(\theta_\eps(\cdot, T')) +  \varphi_\varepsilon(\cdot, T') \, \dx ,
\end{align*}		
for almost all $T^\prime \in [0, T)$, we have
\begin{align}\notag
	&S[\varphi_\varepsilon(\cdot, T'), \theta_\varepsilon(\cdot, T')] -  S[\varphi_{\varepsilon,0}, \theta_{0}]  \\
	& {\geq} \int_{0}^{T'}  \int_{\Omega }	\frac{\ell^\prime (\theta_\varepsilon)}{\ell^2(\theta_\varepsilon)} \kappa(\theta_\varepsilon) |\nabla \theta_{\varepsilon}|^2 \, \dx \mathrm{d}t
	+  \int_{0}^{T^\prime}  \int_{\Omega } \frac{\nu(\theta_\eps )}{\ell(\theta_\varepsilon)} |\nabla \ve_\eps + \nabla \ve_\eps^\mathsf{T}|^2\, \dx \mathrm{d}t \notag \\
	& \quad + \int_{0}^{T^\prime}  \int_{\Omega }	 \frac{\varepsilon }{\ell(\theta_\varepsilon)} \Big| \frac{D\varphi_\eps}{Dt} \Big|^2 \, \dx \mathrm{d}t. \label{eq:entropyldissip}
\end{align}	
  	\end{enumerate}
  \end{coro}

    \begin{osse} 
  	\label{essa}
  	{Observe that the assumption \eqref{assumption} means that the average of the initial absolute temperature in the domain has to be strictly above zero, which is physically sound.}
  	{This assumption \eqref{assumption} is needed due to the lack of mass conservation in the Allen-Cahn equation. Indeed, if we knew that $\int_\Omega \vphi_\eps=\int_\Omega \vphi_{\eps,0}$, then this condition would become simply $\int_\Omega \Lambda(\theta_0)\, \dx>0$ (see \eqref{below} below).}
  	{Furthermore, we observe that 
    the entropy production inequality \eqref{eq:entropyldissip} follow by the proof of its weak counterpart 
    \eqref{eq:weakent} (cf. Section \ref{sec:weaksol}). More precisely, one just needs to test the strong formulation 
    \eqref{eq:entropybal} with a space independent test function coinciding with $1$ on the time interval $(0,T')$ and then integrate over $(0,T')$.}
  \end{osse}

    We can prove an additional property of weak solutions, which is stated into the following result.  

 \begin{lemm}[Uniform strict separation of $\theta_\eps$ from zero]\label{lemma:DeGiorgi}
 Any weak solution to the non-isothermal Navier-Stokes/Allen-Cahn model \eqref{eq:strongmom}-\eqref{eq:strongheat} given by Theorem \ref{theo:weaksol} satisfies the following property:
 \begin{equation}
  	 	\label{thetamin}\text{ if }\alpha\in(\tfrac12,1)\text{ and }\exists c_0>0: \inf_{x\in \Omega}\theta_0\geq c_0,\text{ then }\ \theta_\eps \geq c_\eps(T) >0 \textit{ a.e. in } \Omega \times (0,T),  \, \forall T >0.
  	 \end{equation}
\end{lemm}

 Observe that both the upper and lower bound for the phase-field variable $\varphi_{ \eps}$ from Lemma \ref{lemma:furtherproperties}, as well as the strict positivity of  $\theta_\eps$, will be needed in order to pass to the sharp interface limit ($\varepsilon \rightarrow 0$). for this reason, we will consider again only the class of sublinear and bounded latent heat functions complying with the assumptions listed in Section \ref{subsec:hp}.

  \subsection{$BV$ solutions for non-isothermal two-phase Navier-Stokes/mean curvature flow}
  Our third main result concerns the analysis of the sharp interface limit ($\varepsilon \rightarrow 0$) of the non-isothermal Navier-Stokes/Allen-Cahn system \eqref{eq:strongmom}-\eqref{eq:strongheat} in the weak solutions setting. 
  
  In our analysis of the sharp interface limit of \eqref{eq:strongmom}-\eqref{eq:strongheat}, we will consider only the class of sublinear and bounded latent heat functions $\ell$ (see Section \ref{subsec:hp} for the assumptions). {Indeed, a key ingredient to prove the result is the uniform boundedness of $\varphi_\eps$, which can be ensured only in case of bounded latent heat function $\ell$ (cf. Lemma \ref{lemma:furtherproperties}).

We first provide a notion of weak solutions to the expected sharp interface limit model, i.e., a non-isothermal two-phase Navier-Stokes/mean curvature flow.
  \begin{defi}[$BV$ solutions of non-isothermal two-phase Navier-Stokes/mean curvature flow]\label{def:BVsol}
  	Let $T>0$ be a given final time and $\sigma >0$ be the surface tension.  Let $A(t) \subset \Omega$  be a one-parameter family of open sets with finite perimeter, $t \in [0, T],$ such that {$\chi_A \in L^1(0,T; BV(\Omega))$}. Let $\ve: \Omega \times (0,T) \rightarrow \R^d $  and $\theta: \Omega \times (0,T) \rightarrow (0, \infty)$ be two measurable functions such that 
  	{$\ve \in L^\infty (0,T; \Hs) \cap L^2(0,T; \Vs)\cap W^{1,1}(0,T;\mathbf H^{-m}(\Omega))$, $m>3$, $\theta\in L^\infty (0,T; L^{1+\alpha}(\Omega)) \cap L^2(0,T; H^1(\Omega))$, and $\theta >0$ almost everywhere in $\Omega \times (0,T)$}.
  	We say that the pair $(\chi_A, \ve, \theta)$ is a $B$V solution to the non-isothermal  Navier-Stokes/two-phase mean curvature flow system if
  	\begin{enumerate}
  		\item (Balance of momentum) For any test vector field $\xe \in C^\infty_c([0,T); C^1(\overline{\Omega}; \mathbb{R}^d ))$ with $\operatorname{div}\xe= 0$ and $\xe \cdot  \no = 0$ on $\partial \Omega \times (0, T) $, it holds
  		\begin{align}\label{eqn:momentumsharp}
  			&\int_{0}^T \int_{\Omega} (\ve \cdot \partial_t \xe + (\ve \otimes \ve ): \nabla \xe  
  			) \, \dx \dt +  \int_\Omega \ve_0 \cdot \xe(\cdot, 0) \, \dx \notag \\
  			&=	\int_{0}^T \int_{\Omega} \mathbb{S}(\theta) \ve : \nabla \xe \, \dx \dt + \sigma 	\int_0^T  \int_{\partial^\ast A(t)} (\operatorname{Id} -\n_{A(t)} \otimes \n_{A(t)})  : \nabla \xe \,  \mathrm{d}\mathcal{H}^{d-1} \mathrm{d}t ,
  		\end{align}
  		where $\n_{A(t)} := \frac{\nabla \chi_{A(t)}}{|\nabla \chi_{A(t)}|}$.
  		\item (Existence of normal velocity) There exists a $(\sigma \mathcal{H}^{d-1}\llcorner \partial^*A(t) )\mathrm{d}t$-measurable function V
  		such that
  		\begin{equation}\label{eqn:VL2}
  			\int_0^T \int_{\partial^*A(t)}  | V|^2 \, \mathrm{d}\mathcal{H}^{d-1}\mathrm{d}t < \infty 
  		\end{equation}
  		and $ V (\cdot,  t)$ is the normal speed of $\partial^*A(t)$ with respect to $\n_{A(t)} := \frac{\nabla \chi_{A(t)}}{|\nabla \chi_{A(t)}|}$ in the sense that for almost every $T' \in (0,T)$ and all $\eta \in C^\infty_c(\overline \Omega \times [0,T))$
  		\begin{equation}\label{eqn:transport}
  			\int_{A(T')} \sigma \eta(\cdot, T') \, \mathrm{d}x - 	\int_{A(0)} \sigma \eta(\cdot, 0) \, \mathrm{d}x  = \int_0^{T'} \int_{A(t)} \sigma (\partial_t \eta + \ve \cdot \nabla \eta )\, \mathrm{d}x\mathrm{d}t + \int_0^{T'}
  			\int_{\partial^*A(t)}  \sigma V  \eta  \, \mathrm{d}\mathcal{H}^{d-1}\mathrm{d}t.
  		\end{equation}
  		\item (Motion law) For almost every $T' \in  (0, T)$ and any test vector field $\Psi \in C^1_c (\overline\Omega \times [0,T); \R^d)$ with $\Psi \cdot  \no = 0$ on $\partial \Omega \times (0, T) $, it holds
  		\begin{align}\label{eqn:motionlaw}
  			&- 	\int_{0}^{T'} \int_{\partial^* A(t)} \sigma (\operatorname{Id} -  \n_{A(t)} \otimes  \n_{A(t)} ): \nabla \Psi \, \mathrm{d}\mathcal{H}^{d-1}\mathrm{d}t
  			+	
  			\int_{0}^{T'} \int_{\Omega}  \chi_A  \textnormal{div} (\Psi \ell(\theta ))\, \dx \dt,
  			\notag \\
  			&= 	\int_{0}^{T'} \int_{\partial^* A(t)} \Psi\cdot  n_{ A(t)}  \sigma V \, \mathrm{d}\mathcal{H}^{d-1}\mathrm{d}t. 
  		\end{align} 
  		\item (Total energy inequality) For almost every $T' \in (0,T)$, we have
  		\begin{equation}\label{eqn:energyconservation}
  			E_{\operatorname{tot}}[\chi_A(T'), \ve(T'), \theta (T')] \leq  E_{\operatorname{tot}}[\chi_0, \ve_0, \theta_0],
  		\end{equation}
  		where $$E_{\operatorname{tot}}[\chi_A(T'), \ve(T'), \theta(T') ] :=  \sigma\int_{\Omega} 1 \, \d |\nabla \chi_A(\cdot, T')| + \int_{\Omega} \frac12 |\ve(\cdot,T')|^2 \, \dx + \int_{\Omega} Q(\theta(\cdot,T')) \, \dx .$$
  		\item (Entropy production inequality) 
  		There exists a measure $\mu_\ell \in \mathcal{M}_+({\Omega} \times (0,T))$ such that
  		\begin{align}
  			L  \langle \mu_\ell , 1 \rangle  \geq 	\langle \ell(\theta) \mu_\ell , 1 \rangle
  			\geq   \int_0^T \int_{\partial^*A(t)} \sigma | V|^2 \, \mathrm{d}\mathcal{H}^{d-1}\mathrm{d}t ,
  		\end{align}
  		and,
  		for any test function $\zeta \in C^\infty_c({\Omega}\times (0,T))$ such that $\zeta \geq 0 $, it holds
  		\begin{align}\label{eqn:entropydissipBV}
  			&	
  			-\int_{0}^{T}  \int_{\Omega }	(\Lambda (\theta) + \chi_A)(t) (\partial_t \zeta+ \ve \cdot \nabla \zeta) \,  \dx \mathrm{d}t
  			-	\int_{0}^{T}  \int_{\Omega }	  h(\theta) \Delta\zeta \,  \dx \mathrm{d}t
  			\notag \\
  			&	\geq 
  			\int_{0}^{T}  \int_{\Omega }\frac{\ell^\prime (\theta)}{\ell^2(\theta)} \kappa(\theta) |\nabla \theta|^2 \zeta \, \dx \mathrm{d}t  
  			+\int_{0}^{T}  \int_{\Omega }	\frac{\nu(\theta)}{\ell(\theta)} |\nabla \ve+ \nabla \ve^\mathsf{T}|^2
  			\zeta \, \dx \mathrm{d}t 
  			+ 	\langle \mu_\ell, \zeta \rangle   ,
  		\end{align}
  		where $h(\theta):= \int_{1}^{\theta} \frac{\kappa(s)}{\ell(s)} \, \d s $.
  	\end{enumerate}
  \end{defi}
  \begin{osse}
  	Notice that $\ve\in C_{weak}([0,T],\Hs)$, so that it holds $\ve(0)=\ve_0$, where the limit is intended weakly in $\mathbf L^2(\Omega)$. As in Remark \ref{initialdata}, due to a lack of control on $\partial_t\theta$, the initial datum $\theta_0$ is only enforced weakly in the sense of \eqref{eqn:energyconservation}. 
  \end{osse}
  
  As a third main contribution of this work, we establish the existence of distributional solutions to the non-isothermal two-phase Navier-Stokes/mean curvature flow system in the sense of the previous definition. More precisely, we show that (weak) solutions to the non-isothermal Navier-Stokes/Allen-Cahn system \eqref{eq:strongmom}-\eqref{eq:strongheat} (with sublinear and bounded latent heat) converge subsequentially and conditionally towards such distributional solutions. Before stating our result, we define 
  \begin{itemize}
  	\item[-] the phase indicator function and the surface tension:
  	$$ \psi_{\eps}:= \psi(\varphi_\eps)  = \int_{-1}^{\varphi_{ \eps}} \sqrt{W(s)} \, \mathrm{d} s, \quad  \sigma:= \psi(1)- \psi(-1) =\int_{-1}^{1}\sqrt{2W(s)}\, \d s ;
  	$$
  	\item[-] the phase-field interface energy and the sharp interface energy:
  	$$
  		E_{{\operatorname{int}}, \eps}[\varphi_{\varepsilon}] := \int_{\Omega} \frac{\eps}{2} |\nabla \varphi_{ \eps}|^2 + \frac{1}{\eps} W(\varphi_{ \eps}) \, \dx , \quad 
  		E_{{\operatorname{int}}}[\chi_A] := \sigma \int_{\Omega} 1 \, \d |\nabla \chi_A|. 
  	$$
  \end{itemize}

  \begin{teor}[Conditional convergence of non-isothermal Navier-Stokes/Allen-Cahn system]\label{theo:convergenceBV}
  	Consider a sequence $(\varphi_{\varepsilon,0})_{\varepsilon >0}$ of initial phase field,  an initial velocity field $\ve_{0}$, and an initial temperature $\theta_{0}$, such that 
  	\begin{align*}
  		\sup_{\eps>0} E_{\operatorname{tot}, \eps}[\varphi_{\varepsilon,0},\ve_0,  \theta_{0 }] < \infty , 
  		\quad 	 \|\ve_{ 0 }\|_{\LL^2(\Omega)} < \infty , 
  		\quad 	\|Q (\theta_{ 0 })\|_{L^1(\Omega)} < \infty , \quad 
  		\sup_{\eps>0} \|\varphi_{\varepsilon,0}\|_{L^\infty(\Omega)}\leq 1 ,
  	\end{align*}
  	and such that there exists a set of finite perimeter $A(0)$ and a scalar field $\theta_0$ satisfying
  	\begin{align}
  	\psi_{\eps,0}:=	\psi (\varphi_{\eps,0}) &\rightarrow \sigma \chi_{A(0)}\quad \text{ in } L^1(\Omega) \text{ as } \eps \rightarrow 0,\\
  		E_{{\operatorname{int}}, \eps}[\varphi_{\varepsilon,0}]&\rightarrow E_{{\operatorname{int}}}[\chi_{A(0)}] \quad \text{ as } \eps \rightarrow 0. \label{eqn:convE0}
  	\end{align}
  Let $(\ve_\eps,\varphi_\varepsilon, \theta_\varepsilon)_{\varepsilon>0}$ denote the associated sequence of weak solutions to the non-isothermal Navier-Stokes/Allen-Cahn system in the sense of Definition \ref{def:weaksol} satisfying all the assumptions of Lemma \ref{lemma:furtherproperties}, with $\tfrac12 < \alpha \leq 1$, $2 \leq \beta $. 
  Let $T>0$ be such that the total energy inequality \eqref{ineq:enconsl} and the entropy production inequality \eqref{eq:entropyldissip} hold. 
  	Consider the class of bounded sublinear latent heat functions $\ell$ (see  assumptions in Section \ref{subsec:hp}). 
  	
  	If 
  	\begin{align}	
  		\label{eqn:hpQ} 
  		\limsup_{\eps\to 0} \int_{0}^{T}  \int_\Omega Q(\theta_{\varepsilon}) \,\mathrm{d}x  \mathrm{d}t &\leq   \int_{0}^{T}  \int_\Omega Q(\theta ) \, \mathrm{d}x \mathrm{d}t ,
  	\end{align}
  	then there exists a (nonrelabeled) subsequence $\eps \rightarrow 0$, a one-parameter family of sets of finite perimeter $A(t)$, $t \in [0, T]$, a vector field $\ve: \Omega \times (0,T) \rightarrow \R^d$, and a function $\theta: \Omega \times (0,T) \rightarrow (0, \infty)$ such that
  	\begin{equation}\label{eqn:convergencepsi} 
  		\psi_\eps := \psi(\varphi_\eps ) = \int_{-1}^{\varphi_{\varepsilon}} \sqrt{2W(s)}\, \mathrm{d}s\rightarrow \sigma \chi_A \quad \text{ strongly in } L^1(\Omega \times (0,T)) \text{ as } \eps \rightarrow 0,
  	\end{equation}
  	where {$\chi_A \in H^1(0,T; (H^1(\Omega))') \cap L^\infty_{w*}(0,T; BV(\Omega)) $}, and 
  		\begin{align}
  			&	\ve_{\varepsilon} \rightarrow \ve \quad\text{ strongly in } \LL^2(\Omega \times (0,T)) \text{ as } \eps \rightarrow 0,  \label{eqn:convergencevel} \\
  			&	Q(\theta_{\varepsilon}) \rightarrow Q(\theta) \quad \text{ strongly in } L^1(\Omega \times (0,T)) \text{ as } \eps \rightarrow 0. \label{eqn:convergenceQ}
  	\end{align}
   In addition, assuming the following energy convergence hypotheses
  	\begin{align}\label{eqn:hpenergyconv} 
  		\lim_{\eps\to 0} \int_{0}^{T} E_{{\operatorname{int}}, \eps} [\varphi_{\varepsilon}(\cdot, t)]\, \mathrm{d}t &=  \int_{0}^{T} E_{\operatorname{int}}[\chi_A(\cdot, t)]\, \mathrm{d}t , 
  	\end{align}
  	then $A(t)$, $t \in [0, T]$,  is a $BV$ solution for the non-isothermal two-phase Navier-Stokes/mean curvature flow system in the sense of Definition \ref{def:BVsol} above.
  \end{teor}


 \section{Local-in-time existence of strong solutions}\label{sec::strong}
{In this section we prove the well-posedness of local strong solution as stated in Theorem \ref{thm:strong}. We consider the case of a three-dimensional bounded domain ($d=3$), since the case $d=2$ can be easily deduced by adapting this proof.}

{To this aim, recalling the definitions of $Z_T^i$, $i=1,\ldots,3$, from Section \ref{WPlocal}}, we endow them with the norms
\begin{align}
	& \Vert \f \Vert_{Z^1_{\tT}}:=\Vert \f \Vert_{L^p(0,\tT; \W^{2,p}_\sigma(\Omega))}+\Vert \f\Vert_{ W^{1,p}(0,\tT;\L^p_\sigma(\Omega))}+\Vert \f(0) \Vert_{(\W^{2,p}_\sigma(\Omega),\L^p_\sigma(\Omega))_{1-\frac1p,p}},\\&
	\Vert f\Vert_{\Zdt}:=\Vert f \Vert_{L^q(0,\tT;W^{2,q}_N(\Omega))}+\Vert f \Vert_{W^{1,q}(0,\tT;L^q(\Omega))}+\Vert f(0) \Vert_{(W^{2,q}_N(\Omega),L^q(\Omega))_{1-\frac1q,q}},\\&
	\Vert f\Vert_{\Ztt}:=\Vert f \Vert_{L^r(0,\tT;W^{2,r}_N(\Omega))}+\Vert f \Vert_{W^{1,r}(0,\tT;L^r(\Omega))}+\Vert f(0) \Vert_{(W^{2,r}_N(\Omega),L^r(\Omega))_{1-\frac1r,r}},
\end{align}
 where $r\in(\tfrac{7}2,+\infty),\ p\in (r,p^*)$, with $p^*=\frac{3}2\frac{r}{5-r}>r$ for $r\in(\frac72,5)$, $p^*=+\infty$ for $r\geq 5$, and $q\in (2p,+\infty)$.
Note that these norms are meaningful, since, by, e.g., \cite[Lemma 2]{Saal1}, we can immediately deduce that for any $T_0>0$ there exists $C(T_0)>0$ such that
\begin{align}
	\nonumber	&\Vert \f\Vert_{BUC([0,\tT];\W_\sigma^{2-\frac2p,p}(\Omega))}\leq C(T_0)\Vert \f\Vert_{\Zut},\quad \Vert f\Vert_{BUC([0,\tT];W^{2-\frac2q,q}(\Omega))}\leq C(T_0)\Vert f\Vert_{\Zdt},\\&
	\Vert f\Vert_{BUC([0,\tT];W^{2-\frac{2}{r},r}(\Omega))}\leq C(T_0)\Vert f\Vert_{\Ztt},
	\quad \forall \tT\leq T_0.
	\label{uniform}
\end{align}
Note that, since $q>2p> 6$, it holds
\begin{align}
	(W^{2,q}_N(\Omega),L^q(\Omega))_{1-\frac1q,q}\hookrightarrow W^{\frac 85,5}(\Omega)\hookrightarrow W^{1,\zeta}(\Omega),\quad \forall \zeta\geq2.
	\label{embedding1}
\end{align}
This also entails that
\begin{align}
	Z^2_T\hookrightarrow C(\overline{\Omega}\times[0,T]).
	\label{boundphi}
\end{align}
Note that we also have
\begin{align}
	\label{holder}
	\W_\sigma^{2-\frac2p,p}(\Omega)\hookrightarrow \W_\sigma^{1,3}(\Omega),
\end{align}
since $p>3$.
Furthermore by Sobolev-Gagliardo-Nirenberg's inequality, we deduce, since $3<r<p$,
\begin{align}
	\norm{\f}_{L^{2p}(0,T;\W_\sigma^{1,2p}(\Omega))}\leq C\norm{\f}_{L^\infty(0,T;\W_\sigma^{2-\frac 2p,p}(\Omega))}^\frac12\norm{\f}_{L^p(0,T;\W_\sigma^{2,p}(\Omega))}^\frac12\leq C\norm{\f}_{\Zut},\quad \forall \f\in \Zut,
	\label{essential1}
\end{align}
so that, since $p>r$, 
\begin{align}
	\norm{\f}_{L^{2r}(0,T;\W_\sigma^{1,2r}(\Omega))}\leq CT^{\frac12(\frac1r-\frac 1p)}	\norm{\f}_{L^{2p}(0,T;\W_\sigma^{1,2p}(\Omega))}\leq CT^{\frac12(\frac1r-\frac 1p)}\norm{\f}_{\Zut},\quad \forall \f\in \Zut,
	\label{essential}
\end{align}
and again as for \eqref{uniform} the constant does not depend on the specific $T<T_0$, for $T_0>0$ fixed. Since we have assumed $r>3$, we also have
\begin{align}
	Z_T^3\hookrightarrow L^{r}(0,T;W^{1,\infty}(\Omega)).
	\label{WW}
\end{align}

\noindent In conclusion, since $r>3$, note that we also have the following interpolation and embedding results (see, for instance, \cite{Triebel1}):
\begin{align*}
	(W^{2-\frac2r,r}(\Omega),L^r(\Omega))_{\vartheta,r}\hookrightarrow W^{1,r}(\Omega)\hookrightarrow C^{1-\frac 3r}(\overline{\Omega}),
\end{align*}
Now, since $\Ztt\hookrightarrow W^{1,r}(0,\tT;L^r(\Omega))\hookrightarrow C^{1-\frac 1r}([0,\tT];L^r(\Omega))$, by \cite[Lemma 1]{AWe} we also infer, by the interpolation result above and \eqref{uniform},
\begin{align}
	\nonumber&\Ztt\hookrightarrow C^{\vartheta(1-\frac1r)}([0,\tT];(W^{2-\frac2r,r}(\Omega),L^r(\Omega))_{\vartheta,r})\\&
	\hookrightarrow C^{\vartheta(1-\frac1r)}([0,\tT];W^{1,r}(\Omega))\nonumber
	\\&\hookrightarrow C^{\vartheta(1-\frac1r)}([0,\tT];C^{1-\frac 3r}(\overline\Omega)),	\label{holder1}
\end{align}
so that, in particular, any function in $Z_T^3$ is space-time H\"{o}lder continuous. 

Additionally, we note that, since $r<p\leq p^*$,
$$
W^{2-\frac2r,r}(\Omega)\hookrightarrow W^{1,2p}(\Omega),
$$
entailing
\begin{align}
Z^3_T\hookrightarrow BUC([0,T];W^{1,2p}(\Omega)).
\label{Z3}
\end{align}

We then introduce the following spaces: 
$$
X^1_{{T}}:=\{\v\in \Zut:\ \v_{|t=0}=\v_0\},\quad X^2_{{T}}:=\{\vphi\in \Zdt:\ \vphi_{|t=0}=\vphi_0\},\ X^3_{{T}}:=\{\theta\in \Ztt:\ \theta_{|t=0}=\theta_0\},
$$
$$
\XT:= \XT^1\times \XT^2\times \XT^3,
$$
endowed with the corresponding $Z_T^i$ metrics,
as well as 
$$
Y^1_{{T}}:=L^p(0,{T};\L^p_\sigma(\Omega)), \quad Y_{{T}}^2=L^q(0,T;L^q(\Omega)),\quad  Y_{{T}}^3=L^r(0,T;L^r(\Omega))
$$
$$
\YT:= Y_{{T}}^1\times Y_{{T}}^2 \times Y_{{T}}^3.
$$
Observe that, by \eqref{holder1} together with the assumption $\inf_{x\in\overline{\Omega}}\theta_0\geq c>0$, for any $R>0$ there exists $T_0>0$ such that, also, 
\begin{align}
	\inf_{(x,t)\in \overline{\Omega}\times [0,T]}	\theta(x,t)\geq \frac c2,\quad \forall \theta \in \XT^3,\quad \norm{\theta}_{\XT^3}\leq R,\quad \forall T\leq T_0.
	\label{sd}
\end{align}
From now on we consider $T$ sufficiently small such that \eqref{sd} holds.
 
 \subsection{Proof of Theorem \ref{thm:strong}}
Since $\eps >0$ is fixed in this section, for the sake of simplicity, we omit the index $\eps$ and write $(\ve, \varphi, \theta)$ instead of $(\ve_\eps, \varphi_{\varepsilon}, \theta_\eps)$. 
 
 First, we consider the Leray projector $\mathbf P_\sigma:\mathbf L^2(\Omega)\to \Hs$, and obtain 
 \begin{align}
 	\label{eq:strong3}\partial_t \ve + \P_\sigma ((\ve \cdot \nabla) \ve) -\mathbf P_\sigma\text{div}(\mathbb S(\theta)\ve)= -\P_\sigma\opdiv(\varepsilon \nabla \vphi \otimes \nabla \vphi ),
 \end{align}
 where $\mathbb S(\theta)\ve=\nu(\theta)(\nabla\ve+\nabla\ve^T).$
 Analogously, we also introduce the homogeneous Neumann Laplace operator $A$.
 We are now ready to start the proof. To this aim we rewrite system \eqref{eq:intro1}-\eqref{eq:intro4} as follows 
 \begin{align*}
 	\mathcal{L}(\ve,\vphi,\theta)=\mathcal{F}(\ve,\vphi,\theta),
 \end{align*}
 where we defined the  linear operator $\mathcal{L}:\XT\to\YT$ as 
 \begin{align}
 	\label{L}
 	\mathcal{L}(\ve,\vphi,\theta):=\begin{bmatrix}
 		\partial_t \ve  -\mathbf P_\sigma\text{div}(\mathbb S(\theta_0)\ve)\\
 		\partial_t\vphi+A\vphi \\
 		\partial_t\theta-\frac1{\theta_0^\alpha}\opdiv(\kappa(\theta_0)\nabla\theta)
 	\end{bmatrix},
 \end{align}
 Note that the operators $-\mathbf P_\sigma\text{div}(\mathbb S(\theta_0)\ve),\ A,\ -\frac1{\theta_0^\alpha}\opdiv(\kappa(\theta_0)\nabla\theta)$, {over suitable domains specified below},  have $L^q-L^p$ maximal regularity on the finite interval $[0,T]$ (see, e.g., \cite[Theorem 3.1]{ShibataShimada}, {\cite{DHP}}).
 
 We also defined the nonlinear operator $\mathcal{F}:\XT\to\YT$ as
 \begin{align} 
 	\label{F}
 	\mathcal{F}(\ve,\vphi,\theta):=\begin{bmatrix}
 		\mathbf P_\sigma\text{div}(\mathbb S(\theta)\ve)-\mathbf P_\sigma\text{div}(\mathbb S(\theta_0)\ve)-\P_\sigma((\ve\cdot\nabla)\ve)-\P_\sigma\opdiv(\varepsilon\nabla\vphi\otimes\nabla\vphi),\\
 		-\ve\cdot\nabla\vphi-\frac1{\varepsilon^2}W'(\vphi)+\frac{\ell(\theta)}\varepsilon
 		\\
 		\frac{1}{\theta_0^\alpha}(\theta_0^\alpha-\theta^\alpha)\partial_t\theta-\frac{1}{\theta_0^\alpha}\ve\cdot\nabla Q(\theta)+\frac{1}{\theta_0^\alpha}\opdiv((\kappa(\theta)-\kappa(\theta_0))\nabla\theta)-\frac{1}{\theta_0^\alpha}\ell(\theta)\frac{D\vphi}{Dt}\\+\frac{1}{\theta_0^\alpha}\vert \nu(\theta)(\nabla\ve+\nabla\ve^T)\vert^2+\frac{1}{\theta_0^\alpha}\varepsilon\norma{\frac{D\vphi}{Dt}}^2
 	\end{bmatrix}.
 \end{align}
 We now state the following lemmas, which will be proven later on:   
 \begin{lemm}
 	Let the assumptions of Theorem \ref{thm:strong} hold true. Then
 	there is a constant $C(T, R)>0$ such that
 	$$\Vert\mathcal{F}(\ve_1, \vphi_1,\theta_1) -\mathcal{F}(\ve_2, \vphi_2,\theta_2)\Vert_{\YT}
 	\leq C(T,R)\Vert(\ve_1 - \ve_2, \vphi_1 -\vphi_2,\theta_1-\theta_2)\Vert_{\XT},$$
 	for all $(\ve_i, \vphi_i,\theta_i ) \in \XT$ with $\Vert(\ve_i, \vphi_i,\theta_i )\Vert_{\XT}
 	\leq R$, $R>0$, and $i = 1, 2$. Furthermore it holds
 	$C(T, R) \to0$ as ${T} \to0$.
 	\label{prop1}
 \end{lemm}
 \begin{lemm}
 	\label{prop2}
 	Let ${T} > 0$ and $\mathcal{L}$ as defined above. Then, $\mathcal{L}: \XT \to \YT$
 	is an isomoprhism. Moreover, for every $T_0>0$ there is a constant $C(T_0) > 0$ such that
 	$\Vert\mathcal{L}^{-1}\Vert_{\mathcal{L}(\YT ,\XT )} \leq C(T_0)$ for any $T\in(0,T_0]$.
 \end{lemm}
 With  Propositions \ref{prop1}-\ref{prop2} at hand, which correspond to \cite[Propostion 1-Theorem 4]{AWe}, Theorem \ref{thm:strong} is proved by following \textit{verbatim} the same fixed point argument (contraction mapping theorem) as in the proof of \cite[Theorem 2]{AWe}. The fixed point argument is applied to the problem satisfied by a solution $(\ve,\vphi,\theta)\in \XT$, for some $\tT>0$, i.e.,
 $$
 (\ve,\vphi,\theta)=\mathcal{L}^{-1}\mathcal{F}(\ve,\vphi,\theta) \quad\text{in }\XT.
 $$
 We do not repeat the details here, and we refer to \cite[Theorem 2]{AWe} (see also, for instance, \cite{AGP1}) for this rather classical argument.
 
 \subsection{Intermediate results and proofs}
 
 \subsubsection{Proof of Lemma \ref{prop1} }   	
 Let us consider $(\ve_i, \vphi_i,\theta_i ) \in \XT$ with $\Vert(\ve_i, \vphi_i,\theta_i )\Vert_{\XT}
 \leq R$, $R>0$, and $i = 1, 2$. We need to compute the quantity $\Vert\mathcal F(\ve_1,\vphi_1,\theta_1)-\mathcal F(\ve_2,\vphi_2,\theta_2)\Vert_{\YT}$. We now analyze term by term this quantity. Let us start with the first row. 
 First, we have, recalling that $\P_\sigma$ is also a continuous operator from $\L^p(\Omega)$ to $\L^p_\sigma(\Omega)$
 \begin{align*}
 	&\norm{	-\P_\sigma((\ve_1\cdot\nabla)\ve_1)+	\P_\sigma((\ve_2\cdot\nabla)\ve_2)}_{\Ya}\\&
 	\leq C \norm{((\v_1-\v_2)\cdot \nabla)\v_1}_{\Yan}+C \norm{((\v_2\cdot \nabla)(\v_1-\v_2)}_{\Yan}\\&
 	\leq C\norm{\v_1-\v_2}_{L^p(0,T;\L^\infty(\Omega))}\norm{\nabla\v_1}_{L^\infty(0,T;\L^p(\Omega))}+C\norm{\v_2}_{L^p(0,T;\L^\infty(\Omega))}\norm{\nabla(\v_1-\v_2)}_{L^\infty(0,T;\L^p(\Omega))}\\&
 	\leq C\norm{\v_1-\v_2}_{L^\infty(0,T;\H^1(\Omega))}^\frac12\norm{\v_1-\v_2}_{L^{\frac p2}(0,T;\H^2(\Omega))}^\frac12\norm{\nabla\v_1}_{L^\infty(0,T;\L^p(\Omega))}\\&
    \quad +C\norm{\v_2}_{L^\infty(0,T;\H^1(\Omega))}^\frac12\norm{\v_2}_{L^{\frac p2}(0,T;\H^2(\Omega))}^\frac12\norm{\nabla(\v_1-\v_2)}_{L^\infty(0,T;\L^p(\Omega))}
 	\\&\leq C\norm{\v_1-\v_2}_{L^\infty(0,T;\H^1(\Omega))}^\frac12\norm{\v_1-\v_2}_{L^{ p}(0,T;\H^2(\Omega))}^\frac12T^{\frac{1}{2p}}\norm{\nabla\v_1}_{L^\infty(0,T;\L^p(\Omega))}\\&
    \quad +C\norm{\v_2}_{L^\infty(0,T;\H^1(\Omega))}^\frac12\norm{\v_2}_{L^{ p}(0,T;\H^2(\Omega))}^\frac12T^{\frac{1}{2p}}\norm{\nabla(\v_1-\v_2)}_{L^\infty(0,T;\L^p(\Omega))}\leq C(R)T^\frac{1}{2p}\norm{\v_1-\v_2}_{X_T^1},
 \end{align*}
 where we exploited the fact that $\W^{2,p}(\Omega)\hookrightarrow \H^2(\Omega)$ and we applied Agmon's, H\"{o}lder, and Poincaré's inequalities.
 Secondly, we study the Korteweg force. We have
 \begin{align*}
 	&\norm{-\P_\sigma\opdiv(\varepsilon\nabla\vphi_1\otimes\nabla\vphi_1)+\P_\sigma\opdiv(\varepsilon\nabla\vphi_2\otimes\nabla\vphi_2)}_{\Ya}\\&
 	\leq C\norm{\opdiv(\varepsilon\nabla\vphi_1\otimes\nabla(\vphi_1-\vphi_2))}_{\Yan}+C\norm{\opdiv(\varepsilon\nabla(\vphi_1-\vphi_2)\otimes\nabla\vphi_2)}_{\Yan}\\&
 	\leq 
 	C\norm{\vphi_1-\vphi_2}_{L^\infty(0,T;W^{1,2p}(\Omega))}(\norm{\vphi_1}_{L^p(0,T;W^{2,2p}(\Omega))}+\norm{\vphi_2}_{L^p(0,T;W^{2,2p}(\Omega))})\\&\quad +\norm{\vphi_1-\vphi_2}_{L^p(0,T;W^{2,2p}(\Omega))}(\norm{\vphi_1}_{L^\infty(0,T;W^{1,2p}(\Omega))}+\norm{\vphi_2}_{L^\infty(0,T;W^{1,2p}(\Omega))})
 	\\&
 	\leq 	C(R)T^{\frac1{p}-\frac1q}\norm{\vphi_1-\vphi_2}_{X_T^2},
 \end{align*}
 where we have exploited \eqref{embedding1}, as well as the estimate
 $$
 \norm{f}_{L^{p}(0,T;W^{2,2p}(\Omega))}\leq \norm{f}_{L^{q}(0,T;W^{2,q}(\Omega))}T^{\frac1{p}-\frac 1q},
 $$ 
 since we chose $q>2p$.
In conclusion, we estimate the residual between the evaluation of the stress tensors, namely, recalling $\nu\in W^{1,\infty}(\mathbb R)$ and \eqref{essential1}, together with \eqref{holder1}-\eqref{Z3}, 
 \begin{align*}
 &	\norm{\mathbf P_\sigma\text{div}(\mathbb S(\theta_1)\ve_1)-\mathbf P_\sigma\text{div}(\mathbb S(\theta_2)\ve_2)-\mathbf P_\sigma\text{div}(\mathbb S(\theta_0)\ve_1)+\mathbf P_\sigma\text{div}(\mathbb S(\theta_0)\ve_2)}_{\Ya}\\&
 	\leq C\norm{(\nu'(\theta_1)-\nu'(\theta_2))(\nabla\ve_1+\nabla\ve_1^T)\nabla\theta_2}_{\Yan}+C\norm{\nu'(\theta_1)(\nabla \ve_1+\nabla\ve_1^T)\nabla(\theta_1-\theta_2)}_{\Yan}\\& \quad 
 	+C\norm{(\nu(\theta_1)-\nu(\theta_2))\Delta\ve_1}_{\Yan}+C\norm{\nu'(\theta_2)(\nabla(\ve_1-\ve_2)-\nabla(\ve_1-\ve_2)^T)\nabla\theta_2}_{\Yan}\\& \quad 
    +C\norm{(\nu(\theta_2)-\nu(\theta_0))\Delta(\ve_1-\ve_2)}_{\Yan}
 +C\norm{\nu'(\theta_0)(\nabla(\ve_1-\ve_2)-\nabla(\ve_1-\ve_2)^T)\nabla\theta_0}_{\Yan}\\&
 \leq C(R)\norm{\theta_1-\theta_2}_{L^{2p}(0,T;L^\infty(\Omega))}\norm{\theta_2}_{L^\infty(0,T;W^{1,2p}(\Omega))}\norm{\ve_1}_{L^{2p}(0,T;\mathbf W^{1,2p}(\Omega))}\\& \quad 
 +C(R)\norm{\theta_1-\theta_2}_{L^\infty(0,T;W^{1,2p}(\Omega))}\norm{\ve_1}_{L^{p}(0,T;\mathbf W^{1,2p}(\Omega))}\\& \quad 
 +C(R)\norm{\theta_1-\theta_2}_{L^\infty(0,T;L^\infty(\Omega))}\norm{\ve_1}_{L^p(0,T;\mathbf W^{2,p}(\Omega))}\\&
\quad + C(R)\norm{\theta_2}_{L^{2p}(0,T;L^\infty(\Omega))}\norm{\theta_2}_{L^\infty(0,T;W^{1,2p}(\Omega))}\norm{\ve_1-\ve_2}_{L^{2p}(0,T;\mathbf W^{1,2p}(\Omega))}\\&\quad +
C(R)\norm{\theta_2-\theta_0}_{L^\infty(0,T;L^\infty(\Omega))}\norm{\ve_1-\ve_2}_{L^p(0,T;\mathbf W^{2,p}(\Omega))}\\&
\leq C(R)(T^\frac1{2p}+T^{\vartheta(1-\tfrac1r)})\norm{(\v_1-\v_2,\theta_1-\theta_2)}_{X_T^1\times X_T^3}.
 \end{align*}
 We can now pass to the second line of the operator $\mathcal F$. Namely, we have
 \begin{align*}
 	&	\norm{-\v_1\cdot\nabla\vphi_1+\v_2\cdot\nabla\vphi_2}_{\Yb}\\&
 	\leq   \norm{\v_1\cdot\nabla(\vphi_1-\vphi_2)}_{\Yb}+\norm{(\v_1-\v_2)\cdot\nabla\vphi_2}_{\Yb}\\&
 	\leq \norm{\v_1}_{L^{\infty}(0,T;\L^{2q}(\Omega))}\norm{\vphi_1-\vphi_2}_{L^q(0,T;W^{1,2q}(\Omega))}
    +\norm{\v_1-\v_2}_{L^{\infty}(0,T;\L^{2q}(\Omega))}\norm{\vphi_2}_{L^q(0,T;W^{1,2q}(\Omega))}\\&
 	\leq C(R)T^\frac1q\norm{(\v_1-\v_2,\vphi_1-\vphi_2)}_{X_T^1\times X_T^2},
 \end{align*}
 where we exploited \eqref{holder} as well as $X_T^1\hookrightarrow L^\infty(0,T;\L^{r}(\Omega))$ for any $r\geq2$, due to the embedding $\W^{1,3}(\Omega)\hookrightarrow \L^r(\Omega)$ for any $r\geq 2$.
 
 Proceeding with the estimates, using $W\in C^2(\R)$, we have, recalling that, owing to \eqref{boundphi}, $\norm{\vphi_i}_{L^\infty(\Omega\times(0,T))}\leq C(R)$,
 \begin{align*}
 	&	\frac1{\varepsilon^2}\norm{W'(\vphi_1)-W'(\vphi_2)}_{\Yb}\\&\leq 	\frac1{\varepsilon^2} \max_{s\in[-C(R),C(R)]}\norma{W''(s)}\norm{\vphi_1-\vphi_2}_{\Yb}\\&\leq C(R)T^\frac1q\norm{\vphi_1-\vphi_2}_{X_T^2}.
 \end{align*}
 Then, we have, concerning the sublinear term in $\theta$, since $\norm{\theta_i}_{L^\infty(0,T;L^\infty(\Omega))}\leq C(R)$, $i=1,2$,
 \begin{align*}
 	\norm{\frac1\varepsilon(\ell(\theta_1)-\ell(\theta_2))}_{\Yb}\leq C(R)T^\frac{1}q\norm{\theta_1-\theta_2}_{X_T^3},
 \end{align*}
 where we used \eqref{holder1}.
 
 \noindent We are now left with the third line of $\mathcal F$, which is the most delicate. First, we have, recalling $\inf_{x\in \overline{\Omega}}\theta_0(x)\geq c$,
 \begin{align}
 	\label{oo}	\nonumber&\norm{\frac{1}{\theta_0^\alpha}\left((\theta_0^\alpha-\theta_1^\alpha)\partial_t\theta_1-(\theta_0^\alpha-\theta_2^\alpha)\partial_t\theta_2\right)}_{\Yc}\\&\nonumber\leq 	\frac 1{c^\alpha}\norm{(\theta_0^\alpha-\theta_1^\alpha)(\partial_t\theta_1-\partial_t\theta_2)}_{\Yc}+\frac1{c^\alpha}\norm{(\theta_1^\alpha-\theta_2^\alpha)\partial_t\theta_2}_{\Yc}\\&
 	\leq 	\frac 1{c^\alpha}\norm{\theta_0^\alpha-\theta_1^\alpha}_{L^\infty(0,T;L^\infty(\Omega))}\norm{\partial_t\theta_1-\partial_t\theta_2}_{\Yc}
    \\& \nonumber \quad 
    +	\frac 1{c^\alpha}\norm{\theta_1^\alpha-\theta_2^\alpha}_{L^\infty(0,T;L^\infty(\Omega))}\norm{\partial_t\theta_1}_{\Yc}.
 \end{align}
 Now, by \eqref{holder1}, since $\alpha\in(0,1]$, it holds by \eqref{holder} 
 \begin{align}
 	\nonumber	&\norm{\theta_1^\alpha(t)-\theta_1^\alpha(s)}_{L^\infty(\Omega)}\leq \norm{\theta_1(t)-\theta_1(s)}^\alpha_{L^\infty(\Omega)}\\&\nonumber\leq \norm{\theta_1}^\alpha_{C^{\vartheta(1-\frac1r)}([0,T];L^\infty(\Omega))}\vert t-s\vert^{\vartheta \alpha(1-\frac1 r)}\\&\leq C(R)\vert t-s\vert^{\vartheta\alpha(1-\frac1 r)},\quad \forall 0\leq s\leq t\leq T,\label{time}
 \end{align}
 so that we deduce, choosing $s=0$ and recalling $\theta_1(0)=\theta_0$,
 \begin{align}
 	\norm{\theta_1^\alpha-\theta_0^\alpha}_{L^\infty(0,T;L^\infty(\Omega))}\leq C(R)T^{{\vartheta\alpha(1-\frac1 r)}}.
 	\label{th1}
 \end{align}
 Moreover, note that, thanks to \eqref{sd} (recall that $\alpha\in(0,1]$)
 \begin{align*}
 	\vert \theta_1^\alpha(t)-\theta_2^\alpha(t)\vert \leq \alpha\int_{0}^1(s\theta_1(t)+(1-s)\theta_2(t))^{\alpha-1}\vert \theta_1(t)-\theta_2(t)\vert ds\leq C\vert \theta_1(t)-\theta_2(t)\vert.
 \end{align*}
 Then, since $\theta_1(0)=\theta_2(0)=\theta_0$, 
 \begin{align*}
 	&	\norm{\theta_1^\alpha(t)-\theta_2^\alpha(t)	}_{L^\infty(\Omega)}\leq C	\norm{\theta_1(t)-\theta_2(t)	}_{L^\infty(\Omega)}\\&= \norm{(\theta_1(t)-\theta_2(t))-(\theta_1(0)-\theta_2(0))}_{L^\infty(\Omega)}\\&\leq C\norm{\theta_1-{\theta_2}}_{C^{\vartheta(1-\frac1r)}([0,T];L^\infty(\Omega))}\vert t\vert^{\vartheta(1-\frac1 r)}\\&\leq C\norm{\theta_1-{\theta_2}}_{C^{\vartheta(1-\frac1r)}([0,T];L^\infty(\Omega))} T^{\vartheta(1-\frac1 r)}, 	\quad \forall 0\leq t\leq T,
 \end{align*}
 entailing 
 \begin{align}
 	\norm{\theta_1^\alpha-\theta_2^\alpha	}_{L^\infty(0,T;L^\infty(\Omega))}\leq \norm{\theta_1-\theta_2}_{\XT^3}T^{\vartheta(1-\frac1 r)}.
 	\label{th2}
 \end{align}
 Plugging \eqref{th1} and \eqref{th2} into \eqref{oo}, we infer 
 \begin{align*}
 	&	\norm{\frac{1}{\theta_0^\alpha}\left((\theta_0^\alpha-\theta_1^\alpha)\partial_t\theta_1-(\theta_0^\alpha-\theta_2^\alpha)\partial_t\theta_2\right)}_{\Yc}\\&\leq C(R)T^{\vartheta\alpha(1-\frac1 r)}\norm{\partial_t\theta_1-\partial_t\theta_2}_{\Yc}+CT^{\vartheta(1-\frac1 r)}\norm{\theta_1-\theta_2}_{\XT^3}\norm{\partial_t\theta_2}_{\Yc}\\&
 	\leq C(R)T^{\vartheta\alpha(1-\frac1 r)}\norm{\theta_1-\theta_2}_{X_T^3}+C(R)T^{\vartheta(1-\frac1 r)}\norm{\theta_1-\theta_2}_{\XT^3}.
 \end{align*}
 Then, recalling the definition of $Q(\cdot)$,
 \begin{align*}
 	&	\norm{\frac{1}{\theta_0^\alpha}(-\v_1\theta_1^\alpha\cdot\nabla\theta_1+\v_2\theta_2^\alpha\cdot\nabla\theta_2)}_{\Yc}\\&
 	\leq   C\norm{\v_1\cdot(\theta_1^\alpha-\theta_2^\alpha)\nabla\theta_1}_{\Yc}+C\norm{\v_1\cdot\theta_2^\alpha\nabla(\theta_1-\theta_2)}_{\Yc}\\&\quad +C\norm{(\v_1-\v_2)\cdot\theta_2^\alpha\nabla\theta_2}_{\Yc}\\&
 	\leq \norm{\theta_1^\alpha-\theta_2^\alpha}_{L^\infty(\Omega\times(0,T))}\norm{\v_1}_{L^{2r}(0,T;\L^{\infty}(\Omega))}\norm{\theta_1}_{L^{2r}(0,T;W^{1,r}(\Omega))}\\& \quad +\norm{\theta_2}^\alpha_{L^\infty(\Omega\times(0,T))}\norm{\v_1}_{L^{2r}(0,T;\L^{\infty}(\Omega))}\norm{\theta_1-\theta_2}_{L^{2r}(0,T;W^{1,r}(\Omega))}\\&
\quad+\norm{\theta_2}^\alpha_{L^\infty(\Omega\times(0,T))}\norm{\v_1-\v_2}_{L^{2r}(0,T;\L^{\infty}(\Omega))}\norm{\theta_2}_{L^{2r}(0,T;W^{1,r}(\Omega))}\\&
 	\leq T^{1-\frac1r}\norm{\theta_1-\theta_2}_{X_T^3}\norm{\v_1}_{L^\infty(0,T;\H^1(\Omega))}^\frac12\norm{\v_1}_{L^r(0,T;\H^2(\Omega))}^\frac12T^\frac1{2r}\norm{\theta_1}_{L^\infty(0,T;W^{1,r}(\Omega))}\\& \quad +C(R)\norm{\v_1}_{L^\infty(0,T;\H^1(\Omega))}^\frac12\norm{\v_1}_{L^r(0,T;\H^2(\Omega))}^\frac12T^\frac1{2r}\norm{\theta_1-\theta_2}_{L^\infty(0,T;W^{1,r}(\Omega))}\\&
 	\quad + C(R)\norm{\v_1-\v_2}_{L^\infty(0,T;\H^1(\Omega))}^\frac12\norm{\v_1-\v_2}_{L^r(0,T;\H^2(\Omega))}^\frac12T^\frac1{2r}\norm{\theta_2}_{L^\infty(0,T;W^{1,r}(\Omega))}
 	\\&\leq C(R)T^\frac1{2r}\norm{(\v_1-\v_2,\theta_1-\theta_2)}_{X_T^1\times X_T^3},
 \end{align*}
 where we used Agmon's and H\"{o}lder's inequalities, the fact that $r<p$ and thus $\v_i\in L^r(0,T;\H^{2}(\Omega))$, $i=1,2$, as well as \eqref{uniform} thanks to the embedding $W^{2-\frac2r,r}(\Omega)\hookrightarrow W^{1,r}(\Omega)$, and \eqref{th2}.
 
 Proceeding with the estimates, we have
 \begin{align}
 	\nonumber&\norm{\frac{1}{\theta_0^\alpha}\left(\opdiv((\kappa(\theta_1)-\kappa(\theta_0))\nabla\theta_1)\right)-\frac{1}{\theta_0^\alpha}\left(\opdiv((\kappa(\theta_2)-\kappa(\theta_0))\nabla\theta_2)\right)}_{\Yc}\\&\nonumber
 	\leq C\norm{\opdiv((\kappa(\theta_1)-\kappa(\theta_0))(\nabla\theta_1-\nabla\theta_2))}_{\Yc}+C\norm{\opdiv((\kappa(\theta_1)-\kappa(\theta_2))\nabla\theta_2)}_{\Yc}\\&\nonumber
 	\leq 
 	C\norm{(\nabla\kappa(\theta_1)-\nabla\kappa(\theta_0))\cdot \nabla(\theta_1-\theta_2)}_{L^{r}(0,T;L^{r}(\Omega))}\\&\nonumber \quad +C\norm{\kappa(\theta_1)-\kappa(\theta_0)}_{L^{\infty}(0,T;L^{\infty}(\Omega))}\norm{\theta_1-\theta_2}_{L^{r}(0,T;W^{2,r}(\Omega))}\\&\quad +
 	\nonumber	C\norm{(\nabla\kappa(\theta_1)-\nabla\kappa(\theta_2))	\cdot \nabla\theta_2)}_{L^{r}(0,T;L^{r}(\Omega))}\\&\quad +C\norm{\kappa(\theta_1)-\kappa(\theta_2)}_{L^{\infty}(0,T;L^{\infty}(\Omega))}\norm{\theta_2}_{L^{r}(0,T;W^{2,r}(\Omega))}=\sum_{i=1}^4J_i\label{long}
 \end{align}
 Recalling the definition of $\kappa$, with $\beta\geq1$,
 \begin{align*}
 	&J_1=C\norm{(\nabla\kappa(\theta_1)-\nabla\kappa(\theta_0))\cdot \nabla(\theta_1-\theta_2)}_{L^{r}(0,T;L^{r}(\Omega))}\\&\leq \beta\norm{(\theta_1^{\beta-1}-\theta_0^{\beta-1})\nabla\theta_1\cdot \nabla(\theta_1-\theta_2)}_{\Yc}+\beta\norm{\theta_0^{\beta-1}\nabla(\theta_1-\theta_0)\cdot \nabla(\theta_1-\theta_2)}_{\Yc }\\&
 	\leq \beta \norm{\theta_1^{\beta-1}-\theta_2^{\beta-1}}_{L^\infty(\Omega\times(0,T))}\norm{\theta_1}_{L^\infty(0,T;W^{1,r}(\Omega))}\norm{\theta_1-\theta_2}_{L^r(0,T;W^{1,\infty}(\Omega))}\\& \quad
 	+\beta\norm{\theta_0^{\beta-1}}_{L^\infty(\Omega\times(0,T))}\norm{\theta_1-\theta_0}_{L^\infty(0,T;W^{1,r}(\Omega))}\norm{\theta_1-\theta_2}_{L^r(0,T;W^{1,\infty}(\Omega))}.
 \end{align*}
 Now, analogously to \eqref{th2}, we infer, exploiting \eqref{holder1},
 \begin{align}
 	&	\label{hol0}\norm{\theta_1^{\beta-1}-\theta_0^{\beta-1}}_{L^\infty(0,T;L^\infty(\Omega))}\leq C(R)T^{{\vartheta(1-\frac1 r)}},\\&	\norm{\theta_1^{\beta-1}-\theta_2^{\beta-1}}_{L^\infty(0,T;L^\infty(\Omega))}\leq C(R)T^{{\vartheta(1-\frac1 r)}}\norm{\theta_1-\theta_2}_{\XT^3},
 	\label{th1b}
 \end{align}
 but also 
 \begin{align}
 	\label{hol}
 	&\norm{\theta_1-\theta_0}_{L^\infty(0,T;W^{1,r}(\Omega))}\leq C\norm{\theta_1}_{C^{\vartheta(1-\frac1r)}([0,T];W^{1,r}(\Omega))}T^{\vartheta(1-\frac1r)},\\&
 	\norm{\theta_1-\theta_2}_{L^\infty(0,T;W^{1,r}(\Omega))}\leq C\norm{\theta_1-\theta_2}_{C^{\vartheta(1-\frac1r)}([0,T];W^{1,r}(\Omega))}T^{\vartheta(1-\frac1r)}\label{hol1b}
 \end{align}
 and thus, recalling \eqref{WW}, we get 
 \begin{align*}
 	J_1 &
 	\leq \beta \norm{\theta_1^{\beta-1}-\theta_2^{\beta-1}}_{L^\infty(\Omega\times(0,T))}\norm{\theta_1}_{L^\infty(0,T;W^{1,r}(\Omega))}\norm{\theta_1-\theta_2}_{L^r(0,T;W^{1,\infty}(\Omega))}\\& 
 	\quad +\beta\norm{\theta_0^{\beta-1}}_{L^\infty(\Omega\times(0,T))}\norm{\theta_1-\theta_0}_{L^\infty(0,T;W^{1,r}(\Omega))}\norm{\theta_1-\theta_2}_{L^r(0,T;W^{1,\infty}(\Omega))}\\&
 	\leq C(R)T^{\vartheta(1-\frac1r)}\norm{\theta_1-\theta_2}_{\XT^3}.
 \end{align*}
 In a completely analogous way, thanks to \eqref{th1b} and \eqref{hol1b}, we also infer 
 \begin{align*}
 	J_3=C\norm{(\nabla\kappa(\theta_1)-\nabla\kappa(\theta_2))	\cdot \nabla\theta_2)}_{L^{r}(0,T;L^{r}(\Omega))}\leq C(R)T^{\vartheta(1-\frac1r)}\norm{\theta_1-\theta_2}_{\XT^3}.
 \end{align*}
 Concerning $J_2$ and $J_4$, 
 we observe that, analogously to the computation in \eqref{th2}, it is easy to have
 \begin{align}
 	&\norm{\kappa(\theta_1)-\kappa(\theta_0)	}_{L^\infty(0,T;L^\infty(\Omega))}\leq C\norm{\theta_1}_{\XT^3}T^{\vartheta(1-\frac1 r)}\leq C(R)T^{\vartheta(1-\frac1 r)},
 	\label{th2v}\\&
 	\norm{\kappa(\theta_1)-\kappa(\theta_2)	}_{L^\infty(0,T;L^\infty(\Omega))}\leq CT^{\vartheta(1-\frac1 r)}\norm{\theta_1-\theta_2}_{\XT^3}.
 	\label{th2vc}
 \end{align}
 Therefore, since $\theta_i\in L^r(0,T;W^{2,r}(\Omega))$, $i=1,2$, 
 \begin{align*}
 	J_2+J_4\leq C(R)T^{\vartheta(1-\frac1 r)} \norm{\theta_1-\theta_2}_{\XT^3},
 \end{align*}
 and this concludes the estimate of the terms in \eqref{long}.
 We now recall that
 $$
 \frac{D\vphi}{Dt}=\partial_t\vphi+\v\cdot\nabla\vphi.
 $$
 Then we have (here and in the sequel we neglect the factor $\tfrac1{\theta_0^\alpha}$ in front, since this is always bounded in $L^\infty(\Omega\times(0,T))$), recalling that $\ell$ is sublinear in its argument,
 \begin{align*}
 	&\norm{-\ell(\theta_1)\frac{D\vphi_1}{Dt}+\ell(\theta_2)\frac{D\vphi_2}{Dt}}_{\Yc}\\&=	\norm{-\ell(\theta_1)(\partial_t\vphi_1+\v_1\cdot\nabla\vphi_1)+\ell(\theta_2)(\partial_t\vphi_2+\v_2\cdot\nabla\vphi_2)}_{\Yc}
 	\\&
 	\leq 
 	C(R)\norm{(\theta_1-\theta_2)\partial_t\vphi_1}_{\Yc}+\norm{\ell(\theta_2)\partial_t(\vphi_1-\vphi_2)}_{\Yc}\\& \quad
 	+C(R)\norm{(\theta_1-\theta_2)\v_1\cdot\nabla\vphi_1}_{\Yc}+\norm{\ell(\theta_2)(\v_1-\v_2)\cdot\nabla\vphi_1}_{\Yc}\\&
    \quad +\norm{\ell(\theta_2)\v_2\cdot\nabla(\vphi_1-\vphi_2)}_{\Yc}\\&
 	\leq 
 	C(R)\norm{\theta_1-\theta_2}_{L^{2r}(0,T;L^{2r}(\Omega))}\norm{\partial_t\vphi_1}_{L^{2r}(0,T;L^{2r}(\Omega))}\\&\quad +C\norm{\ell(\theta_2)}_{L^{2r}(0,T;L^{2r}(\Omega))}\norm{\partial_t(\vphi_1-\vphi_2)}_{L^{2r}(0,T;L^{2r}(\Omega))}
 	\\& \quad
 	+C(R) \norm{\theta_1-\theta_2}_{L^\infty(0,T;L^\infty(\Omega))}\norm{\v_1}_{L^r(0,T;\L^\infty(\Omega))}\norm{\vphi_1}_{L^\infty(0,T;W^{1,r}(\Omega))}\\&
 	\quad+C\norm{\ell(\theta_2)}_{L^\infty(0,T;L^\infty(\Omega))}\norm{\v_1-\v_2}_{L^r(0,T;\L^\infty(\Omega))}\norm{\vphi_1}_{L^\infty(0,T;W^{1,r}(\Omega))}\\&
 	\quad+C\norm{\ell(\theta_2)}_{L^\infty(0,T;L^\infty(\Omega))}\norm{\v_2}_{L^r(0,T;\L^\infty(\Omega))}\norm{\vphi_1-\vphi_2}_{L^\infty(0,T;W^{1,r}(\Omega))}\\&
 	\leq C(R)T^\frac1{2r}\norm{\theta_1-\theta_2}_{L^\infty(0,T;L^{\infty}(\Omega))}\norm{\partial_t\vphi_1}_{L^{2r}(0,T;L^{2r}(\Omega))}\\&\quad+CT^\frac1{2r}\norm{\ell(\theta_2)}_{L^\infty(0,T;L^\infty(\Omega))}\norm{\partial_t(\vphi_1-\vphi_2)}_{L^{2r}(0,T;L^{2r}(\Omega))}\\&
 	\quad+C(R)\norm{\theta_1-\theta_2}_{L^\infty(0,T;L^{	\infty}(\Omega))}\norm{\v_1}_{L^\infty(0,T;\H^1(\Omega))}^\frac12\norm{\v_1}_{L^\frac r2(0,T;\H^2(\Omega))}^\frac12\norm{\vphi_1}_{L^\infty(0,T;W^{1,r}(\Omega))}\\&
 	\quad+C\norm{\ell(\theta_2)}_{L^\infty(0,T;L^\infty(\Omega))}\norm{\v_1-\v_2}_{L^\infty(0,T;\H^1(\Omega))}^\frac12\norm{\v_1-\v_2}_{L^\frac r2(0,T;\H^2(\Omega))}^\frac12\norm{\vphi_1}_{L^\infty(0,T;W^{1,r}(\Omega))}\\&
 	\quad+C\norm{\ell(\theta_2)}_{L^\infty(0,T;L^\infty(\Omega))}\norm{\v_2}_{L^\infty(0,T;\H^1(\Omega))}^\frac12\norm{\v_2}_{L^\frac r2(0,T;\H^2(\Omega))}^\frac12\norm{\vphi_1-\vphi_2}_{L^\infty(0,T;W^{1,r}(\Omega))}\\&
 	\leq C(R)T^\frac1{2r}\norm{(\vphi_1-\vphi_2,\theta_1-\theta_2)}_{X_T^2\times X_T^3}+C(R)T^{\frac12-\frac r{4p}}\norm{(\v_1-\v_2,\vphi_1-\vphi_2,\theta_1-\theta_2)}_{X_T},
 \end{align*}
 where again we exploited Agmon's and H\"{o}lder's inequalities, \eqref{uniform}-\eqref{embedding1}, together with the estimate
 $$
 \norm{\f}_{L^\frac r2(0,T;\H^2(\Omega))}\leq T^{\frac2r-\frac 1{p}}\norm{\f}_{L^p(0,T;\W^{2,p}(\Omega))}.
 $$
 Note that $q>2r$, so that $\partial_t\vphi_i\in L^{2r}(0,T;L^{2r}(\Omega))$, $i=1,2$.
 
 \noindent
 Proceeding with the computations, we have
 \begin{align*}
 	&	\norm{\nu(\theta_1)^2\vert \nabla\v_1+\nabla\v_1^T	\vert^2-\nu(\theta_2)^2\vert \nabla\v_2+\nabla\v_2^T	\vert^2}_{\Yc}\\&
 	\leq \norm{\nu(\theta_1)(\nabla\v_1+\nabla\v_1^T)+\nu(\theta_2)(\nabla\v_2+\v_2)^T):(\nu(\theta_2)(\nabla(\v_1-\v_2)+\nabla(\v_1-\v_2)^T))}_{\Yc}\\&
 	\quad+\norm{\nu(\theta_1)(\nabla\v_1+\nabla\v_1^T)+\nu(\theta_2)(\nabla\v_2+\v_2)^T):((\nu(\theta_1)-\nu(\theta_2))(\nabla\ve_1+\nabla\ve_1^T))}_{\Yc}\\&
 	\leq C(\norm{\v_1}_{L^{2r}(0,T;\W^{1,2r}_\sigma(\Omega))}+\norm{\v_2}_{L^{2r}(0,T;\W_\sigma^{1,2r}(\Omega))})\norm{\v_1-\v_2}_{L^{2r}(0,T;\W_\sigma^{1,2r}(\Omega))}\\&
 	\quad+C(\norm{\v_1}_{L^{2r}(0,T;\W^{1,2r}_\sigma(\Omega))}+\norm{\v_2}_{L^{2r}(0,T;\W_\sigma^{1,2r}(\Omega))})\norm{\theta_1-\theta_2}_{L^\infty(0,T;L^\infty(\Omega))}\norm{\ve_1}_{L^{2r}(0,T;\mathbf W^{2,r}(\Omega))}\\&
 	\leq C(R)(T^{\frac1r-\frac1p}+T^{\vartheta(1-\tfrac1r)})\norm{(\v_1-\v_2,\theta_1-\theta_2)}_{X_T^1\times X_T^3},
 \end{align*}
 where we have used H\"{o}lder's inequality, \eqref{essential}, and \eqref{holder1}.
 
 In conclusion, we note that, for $i=1,2$,
 \begin{align*}
 	\norm{\frac{D\vphi_i}{Dt}}_{L^{2r}(0,T;L^{2r}(\Omega))}\leq \norm{\partial_t\vphi_i}_{L^{2r}(0,T;L^{2r}(\Omega))}+\norm{\v_i}_{L^\infty(0,T;\L_\sigma^{4r}(\Omega))}\norm{\vphi_i}_{L^{2r}(0,T;W^{1,4r}(\Omega))}\leq C(R),
 \end{align*}
 recalling that $q>2p>6$, and the embeddings \eqref{embedding1} and \eqref{holder}.
 Therefore, we can write 
 \begin{align*}
 	&\norm{\varepsilon\left(\norma{\frac{D\vphi_1}{Dt}}^2-\norma{\frac{D\vphi_2}{Dt}}^2\right)}_{L^{r}(0,T;L^{r}(\Omega))}\\&=\norm{\left(\frac{D\vphi_1}{Dt}+\frac{D\vphi_2}{Dt}\right)(\partial_t(\vphi_1-\vphi_2)+\v_1\cdot\nabla(\vphi_1-\vphi_2)+(\v_1-\v_2)\cdot\nabla\vphi_2)}_{\Yc}\\&
 	\\&	\leq
 	\varepsilon \left(	\norm{\frac{D\vphi_i}{Dt}}_{L^{2r}(0,T;L^{2r}(\Omega))}+	\norm{\frac{D\vphi_i}{Dt}}_{L^{2r}(0,T;L^{2r}(\Omega))}\right)\\&\quad \times(\norm{\partial_t(\vphi_1-\vphi_2)}_{L^{2r}(0,T;L^{2r}(\Omega))}+\norm{\v_1}_{L^\infty(0,T;\L_\sigma^{4r}(\Omega))}\norm{\vphi_1-\vphi_2}_{L^{2r}(0,T;W^{1,4r}(\Omega))}\\&\quad +\norm{\v_1-\v_2}_{L^\infty(0,T;\L_\sigma^{4r}(\Omega))}\norm{\vphi_2}_{L^{2r}(0,T;W^{1,{4r}}(\Omega))})\\&
 	\leq C(R)(T^{\frac1{2r}-\frac{1}{q}}\norm{\partial_t(\vphi_1-\vphi_2)}_{L^q(0,T;L^q(\Omega))}+C(R)T^\frac{1}{2r}\norm{\vphi_1-\vphi_2}_{L^\infty(0,T;W^{1,4r}(\Omega))}\\&\quad +C(R)T^\frac1{2r}\norm{\v_1-\v_2}_{L^\infty(0,T;\L_\sigma^{4r}(\Omega))})\\&
 	\leq C(R)T^{\frac1{2r}-\frac{1}q}\norm{\vphi_1-\vphi_2}_{X_T^2}+C(R)T^\frac1{2r}\norm{(\v_1-\v_2,\vphi_1-\vphi_2)}_{X_T^1\times X_T^2},
 \end{align*}
 where we have used \eqref{embedding1}, \eqref{holder} and the embedding $\W_\sigma^{1,3}(\Omega)\hookrightarrow \L^\iota_\sigma(\Omega)$, for any $\iota\geq 2$.
 
 To sum up, putting together all the estimates, we have shown that there is a constant $C(T, R)>0$ such that
 $$\Vert\mathcal{F}(\v_1, \vphi_1,\theta_1) -\mathcal{F}(\v_2, \vphi_2,\theta_2)\Vert_{\YT}
 \leq C(T, R)\Vert(\v_1 - \v_2, \vphi_1 -\vphi_2,\theta_1-\theta_2)\Vert_{\XT},$$
 for all $(\v_i, \vphi_i,\theta_i ) \in \XT$ with $\Vert(\v_i, \vphi_i,\theta_i )\Vert_{\XT}
 \leq R$, $R>0$, and $i = 1, 2$. As it is clear from the estimates, the positive constant $C(T,R)$ vanishes as $T\to0$.
 This concludes the proof of Lemma \ref{prop1}.

 \subsubsection{Proof of Lemma \ref{prop2}}
 Let us first observe that the three rows of $\mathcal{L}$ are independent from each other, so that we can study them separately. In particular, for the first line, the result directly comes from the fact that the Stokes operator with slip boundary conditions, $\mathbf A_S:=-\P_\sigma\text{div}(\mathbb S(\theta_0)\cdot)$, since $\nu(\theta_0)\in W^{1,r}(\Omega)$, $r>3$, has the $L^q$-$L^p$ maximal regularity property (see, for instance, \cite[Theorem 3.1]{ShibataShimada}): given an initial datum in a suitable trace space, i.e., $\v_0\in (\mathfrak D(\mathbf A_S),\L^p_\sigma(\Omega))_{1-\frac1p,p}$, and a function $\f\in L^p(0,\tT;\L^p_\sigma(\Omega))$ (for $p$ given in the definition of $\Zut$), there exists a unique $\v\in \Zut$ such that
 \begin{align*}
 	\begin{cases}
 		\partial_t\v-\P_\sigma\text{div}(\mathbb S(\theta_0)\ve)=\f,&\quad\text{ in }\Omega\times(0,\tT),\\
 		\opdiv\v=0, &\quad\text{ in }\Omega\times(0,\tT),\\
 		\v(\cdot, 0)=\v_0,&\quad\text{ in }\Omega,\\
 		 \mathbf v \cdot \no =0,&\quad\text{ on }\partial\Omega\times(0,T),\\
 			(\operatorname{Id}-\no\otimes\no)((\mathbb S(\theta_0)\ve)\no)=\mathbf 0,&\quad\text{ on }\partial\Omega\times(0,T).
 	\end{cases}
 \end{align*}   
 Analogously, the homogeneous Neumann Laplace operator A has the $L^q$-$L^p$ maximal regularity property (see, for instance, {\color{black}\cite{DHP}}), so that, given an initial datum $\vphi_0\in (W^{2,q}_N(\Omega),L^q(\Omega))_{1-\frac1q,q}$, and a function $g\in L^q(0,\tT;L^q(\Omega))$ (for $q\in(2p,6]$ given in the definition of $\Zdt$), there exists a unique $\vphi\in \Zdt$ such that
 \begin{align*}
 	\begin{cases}
 		\partial_t\vphi+A\vphi=g,&\quad\text{ in }\Omega\times(0,\tT),\\
 		\vphi(\cdot, 0)=\vphi_0,&\quad\text{ in }\Omega,\\
 		\partial_{\no}\vphi=0,&\quad\text{ on }\partial\Omega\times(0,T).
 	\end{cases}
 \end{align*}   	   	
 Moreover, given an initial datum  $\theta_0\in (W^{2,r}_N(\Omega),L^r(\Omega))_{1-\frac1r,r}$, and a function $h\in L^r(0,\tT;L^r(\Omega))$, there exists a unique $\theta\in \Ztt$ such that
 \begin{align*}
 	\begin{cases}
 		\partial_t\theta-\frac{1}{	\theta_0^\alpha}\opdiv(\kappa(\theta_0)\nabla\theta)=h,&\quad\text{ in }\Omega\times(0,\tT),\\
 		\theta(\cdot, 0)=\theta_0,&\quad\text{ in }\Omega,\\
 		\partial_{\no}\theta=0,&\quad\text{ on }\partial\Omega\times(0,T).
 	\end{cases}
 \end{align*}   	
 This is possible since $\theta_0\in W^{1,s}(\Omega)$, for some $s>3$, and thus the principal part of the operator has H\"{o}lder continuous coefficients. Indeed, the operator in non-divergence form can be written as
 $$
 -\frac{1}{	\theta_0^\alpha}\opdiv(\kappa(\theta_0)\nabla\theta)=	-\frac{1}{	\theta_0^\alpha}\nabla\kappa(\theta_0)\cdot\nabla\theta-\frac{1}{	\theta_0^\alpha}\kappa(\theta_0)\Delta\theta,
 $$ 
 and its principal part $-\frac{1}{	\theta_0^\alpha}\kappa(\theta_0)\Delta\theta$ has indeed H\"{o}lder continuous coefficients (bounded below since by assumption $\inf_{x\in \overline{\Omega}}\theta_0(x)>c>0$). Therefore the operator still has the $L^p$ maximal regularity property (see, e.g., {\cite{DHP}}).
 
 As a consequence, the operator $\mathcal{L}$ is  invertible as a linear operator from $\YT$ to $\XT$, i.e., by the bounded inverse theorem, there exists $C(\tT)>0$, possibly depending on $\tT$, such that 
 $$
 \norm{\mathcal{L}^{-1}}_{\mathcal{L}(\YT,\XT)}\leq C(\tT),\quad \forall \tT>0.
 $$
 To conclude the proof of Lemma \ref{prop2} it is enough to show that the constant above does not change with $\tT$. This can be obtained by a standard extension argument (see for instance the proof of \cite[Lemma 7]{AWe}), leading to show that, given any $T_0>0$,
 $$
 \norm{\mathcal{L}^{-1}}_{\mathcal{L}(\YT,\XT)} \leq \norm{\mathcal{L}^{-1}}_{\mathcal{L}(Y_{T_0},X_{T_0})}\leq C(T_0),\quad \forall 0<\tT\leq T_0,
 $$
 concluding the proof of the proposition.
 

 \section{Global-in-time existence of weak solutions} 
 \label{sec:weaksol}

This section is devoted to the proof of the global-in-time existence of weak solutions, as stated in Theorem \ref{theo:weaksol}. 
We start by briefly sketching our strategy.

We introduce a suitable approximation of the non-isothermal Navier-Stokes/Allen-Cahn system \eqref{eq:strongmom}-\eqref{eq:strongheat} depending on the parameters $k \in \N$, $\delta>0$, and $m\in \N$, and we find a solution to the weak formulation of the corresponding approximate problem by adapting a Galerkin scheme (cf. Section \ref{subsec:approx}).
We derive regularity estimates holding uniformly for the approximate solutions $(\ve_\delta,p_\delta,   \varphi_\delta, \theta_\delta)$, and we first comment on the convergence $k\to \infty$ (cf.~Section \ref{subsec:convk}). Note that, by a slight abuse of notation for the sake of readability, we do not rename the limit as $k\to \infty$.
Then, we show the uniform boundedness of $\varphi_\delta$, as well as the strict positivity of the temperature $\theta_\delta$, and we pass to the limit as $\delta\to0$ (cf.~Section \ref{subsec:delta}).
Secondly, we deduce further regularity estimates to pass to the limit $m\to\infty$, obtaining a quadruple $(\ve, p,\varphi, \theta)$ solving \eqref{eq:weakbalmom}-\eqref{eq:weakent} and \eqref{eq:weaktoten} (cf.~Section \ref{subsec:limitm}).

Finally, we apply a De Giorgi iteration scheme and show the strict separation from zero for a solution to \eqref{eq:strongheat} whenever $\alpha <1$ in \eqref{eq:defQ} (cf.~Section \ref{Degiorgi}). In particular, we prove Lemma \ref{lemma:DeGiorgi}, deducing that $\theta_\eps \geq c_\eps> 0$ almost everywhere in $\Omega\times (0,T)$, as long as also the initial datum $\theta_{0}$ is uniformly strictly separated from zero.

We will consider the case when $\ell$ is bounded, as the in assumption (b) of Section \ref{subsec:hp}, excluding the linear case \eqref{hp:linearell}.
Note that boundedness of $\ell$ is not only required at the level of the approximation scheme, but also for the final passage to limit of the entropy production inequality.

 Throughout the whole section, we will adopt the following choices for the sake of simplicity.
We consider a smooth bounded domain in $\R^d$ with $d=3$, as for the case $d=2$ the same (and even better) regularity properties hold true.
 Since $\eps >0$ is fixed in this section, we omit the indices $\eps$ and write $(\ve, p, \varphi, \theta)$ instead of $(\ve_\eps, p_\eps,  \varphi_{\varepsilon}, \theta_\eps)$. However, we keep writing the parameter $\eps $ in the Section \ref{subsec:delta} in order to keep track of which regularity estimates are uniform (or not) in $\eps $. This will be useful in the following section for our analysis of the sharp interface limit (i.e., letting $\eps \rightarrow 0$).
 Moreover, the same letter $C$ (and, sometimes, with an index when accounting for the dependence on a parameter) will denote positive constants whose actual value may vary from line to line and which do not depend on the approximation scheme adopted.


\subsection{The $(k,\delta,m)$-approximation scheme}
\label{subsec:approx}
The idea of the scheme is to first approximate the functions $\ell$ and $Q$, so that they do no longer vanish at $s=0$, introducing suitable approximations $\ell_\delta$ and $Q_\delta$ for $\delta>0$.  The second approximation layer consists of the regularization of the Allen-Cahn equation by adding the Cahn-Hilliard-type term $\frac1k\Delta\widetilde\mu$, $k\in\N$ where $\widetilde \mu:=-\eps \Delta \varphi+\frac1\varepsilon \partial_\varphi W(\varphi)$, together with a cutoff of the quadratic terms in the heat equation. Finally, we further regularize the heat equation by adding the term $-\frac1m \Delta Q_\delta(\theta)$, $m\in \N$, so that the equation becomes parabolic in the variable $Q_\delta(\theta)$. In addition, a symmetric smooth convolution kernel $J_m$ is introduced into the quadratic terms. At the same approximation level $m$, we also use $J_m$ to regularize the advective term in the Navier-Stokes equation. 

To find a global solution to the approximate problem, for fixed $k\in\N$, $\delta>0$, and $m\in \N$, we use a Galerkin scheme. Note that an energy identity is satisfied for fixed $k\in\N$, $\delta>0$, and $m\in\N$.
Then, we first pass to the limit as $k\to\infty$, as the Cahn-Hilliard-type approximation prevents the application of a maximum principle, which is exploited in our the proof. 
As a drawback, letting $k\to\infty$ yields (only) an energy inequality, obtained through lower semicontinuity properties. We then pass to the limit as $\delta\to0$, and finally as $m\to\infty$.

\subsubsection{Approximation of the system}
We introduce step-by-step all the suitable candidates to approximate the latent heat $\ell$, the heat $Q$, and heat conductivity $\kappa$ at the level $\delta>0$. Then, we present the approximation system, including the level $m \in \N$.
\paragraph{\it Construction of $\ell_\delta$}
Our aim is to give a suitable approximation for the latent heat function $\ell\geq0$ by a strictly positive function to rigorously derive the entropy balance \eqref{eq:entropybal}, as \textit{a priori} we do not know that $\theta>0$. 
Moreover, we introduce a function $\log_\delta$ approximating the logarithm and playing a crucial role in proving the strict positivity of the temperature (cf. Section \ref{subsubsec:strictpostheta}).

For $ \delta \in (0, \delta_0 \ell'(0))$, where $\delta_0\in(0,1)$ is arbitrarily chosen, consider the Lipschitz function:
\begin{align}
\log_\delta(s):=\begin{cases}
\log s,\quad\text{ for }s\geq \delta,\\
\frac{1}{\delta}s+\log\delta-1,\quad\text{ for }s< \delta,
\end{cases}\label{logdelta}
\end{align}
which is a nonsingular approximation of the logarithmic function. We compute its first derivative, which reads
\begin{align}
\log_\delta'(s):=\begin{cases}
\frac1s,\quad\text{ for }s > \delta,\\
\frac{1}{\delta},\quad\text{ for }s< \delta,
\end{cases} \label{log'delta}
\end{align}
and, clearly, is not differentiable at $s=\delta$.

Introducing the Lipschitz approximation of $\tilde l(s) = s$ as $$\widetilde l_\delta(s):= 
 \begin{cases} 
s,\quad\text{ for }s\geq  \delta ,\\
\frac{\delta}{2- \frac{s}{\delta}},\quad\text{ for }s< \delta, \end{cases}$$ we have
\begin{align}
\widetilde l_\delta'(s)
=\begin{cases} 
1,\quad\text{ for }s> \delta,\\
\frac{1}{(2-\frac{s}{\delta})^2},\quad\text{ for }s< \delta.
\end{cases}
\end{align}
Notice that we obtain 
\begin{align}
\frac{\widetilde l_\delta'(s)}{\widetilde l_\delta(s)^2}=\begin{cases} 
\frac {1}{s^2},\quad\text{ for }s> \delta,\\
\frac{1}{\delta^2},\quad\text{ for }s< \delta,
\end{cases}
\end{align}
and, in particular,  
\begin{align}
\frac{\widetilde l_\delta'(s)}{\widetilde l_\delta(s)^2}=\log_\delta'(s)^2.\label{logdelta1}
\end{align}
Moreover, observe that, as $\delta\to 0$,
\begin{align}
\widetilde l_\delta(s)\to \begin{cases} s,\quad\text{ for }s> 0,\\
0\quad\text{ for }s< 0,
\end{cases}\label{f1tilde}
\end{align}
 as well as 
\begin{align}
\widetilde l_\delta'(s)\to \begin{cases} 1,\quad\text{ for }s> 0,\\
0\quad\text{ for }s< 0.
\end{cases}\label{f2tilde}
\end{align}

We recall that, by assumptions on $\ell$ (cf. Section \ref{subsec:hp}), for some $c_l>0$, we can find $0<\lambda_0\leq \lambda_1$ such that
\begin{align}
\lambda_0s\leq \ell(s)\leq \lambda_1 s, \quad \forall s\in(0,c_l).
\label{control}
\end{align}
Also, we can extend the function $\ell$ by setting $\ell(s)=\ell'(0)s$ for $s<0$, so that $\ell$ is a Lipschitz function, as $l(0)=0$. 
However, we observe that this extension does not play any role in the approximation, as in the limit $\delta\to 0$ we will only consider the domain $(0,\infty)$. As an approximation for $\ell$, we define \begin{align} \label{elllambda}
\ell_\delta(s):=\widetilde{l}_\delta(\ell(s))
= \begin{cases} 
\ell(s),\quad\text{ for } \ell(s)\geq  \delta ,\\
\frac{\delta}{2- \frac{\ell(s)}{\delta}},\quad\text{ for } \ell(s)< \delta, \end{cases}
, 
\end{align}
giving
\begin{align*}
\ell_\delta'(s)=\widetilde{l}_\delta'(\ell(s))\ell'(s).
\end{align*}
From assumption (e), we can deduce that $\min_{s\in [0,c_l)}{\ell'(s)}= C_l>0$. Moreover, from \eqref{f1tilde}-\eqref{f2tilde} it follows that, as $\delta\to0$,
\begin{align}
\ell_\delta(s)\to \begin{cases} \ell(s),\quad\text{ for }s> 0,\\
0,\quad\text{ for }s< 0,
\end{cases}\label{f1}
\end{align}
as well as 
\begin{align}
\ell_\delta'(s)\to \begin{cases} \ell'(s),\quad\text{ for }s> 0,\\
0,\quad\text{ for }s< 0.
\end{cases}\label{f2}
\end{align}
In conclusion, from \eqref{logdelta1} we infer 
\begin{align}
\frac{ \ell_\delta'(s)}{ \ell_\delta(s)^2}=\frac{\widetilde l_\delta'(\ell(s))\ell'(s)}{\widetilde l_\delta(\ell(s))^2}=\log_\delta'(\ell(s))^2\ell'(s)\geq C_l \log_\delta'(\ell(s))^2\geq \underbrace{C_{l}\min\{1,\frac1{\lambda_1^2}\}}_{:=C_0>0}\log_\delta'(s)^2,\label{logdelta2}
\end{align}
where we used \eqref{log'delta} and \eqref{control}. Indeed, for $\ell(s)<\delta$,
\begin{align*}
\log_\delta'(\ell(s))^2=\frac 1\delta=\log_\delta'(s),
\end{align*}
whereas, for $\ell(s)>\delta$,
\begin{align*}
\log_\delta'(\ell(s))^2=\frac 1{\ell(s)^2}\geq \frac1{\lambda_1^2} \log_\delta'(s)^2.
\end{align*}
It is now clear that $\ell_\delta$ is a good candidate to approximate the latent heat function $l$ as $\ell_\delta(s)>0$ for any $s\in \R$. 

\paragraph{\it Definition of $Q_\delta$ and $\kappa_\delta$}
We introduce the monotone (thus invertible) function
\begin{align}\label{defQdelta}
Q_\delta(s):=\begin{cases}
    \frac{(s^2+\delta^2)^{\frac{1+\alpha}2}}{1+\alpha},\quad s\geq 0,\\\, \\
    \frac{2\delta^{{1+\alpha}}}{(1+\alpha)(2-s)},\quad s<0,
\end{cases}
\end{align}
 as a good candidate to approximate the heat function $Q$  given by \eqref{eq:defQ} so that $Q(s)>0$ for any $s\in \R$.
We also extend (without relabeling) the heat conductivity function $\kappa$ defined in \eqref{eq:defkappa} by means of
\begin{align}\label{defkappadelta}
\kappa_\delta(s):=\begin{cases}
\kappa(s)\quad \text{ if }s>\delta,\\
\kappa_1+\kappa_2\delta^2\norma{s}^{\beta-2},\quad \text{ if }s\leq \delta.
\end{cases}
\end{align}

\paragraph{\it Approximation problem}
We aim at solving, in a weak formulation, the following equations. 
\begin{itemize}
\item[-] The approximate Navier-Stokes equation:
	\begin{align} \label{eq:weakbalmomapprox}
		&\int_{0}^T \int_{\Omega} (\ve_\delta \cdot \partial_t \xe + (\ve_\delta  \otimes \P_\sigma ( J_m\ast\ve_\delta)  ): \nabla \xe  
		) \, \dx \dt \notag \\
		&=	\int_{0}^T \int_{\Omega} \mathbb{S}(\theta_\delta)\ve_\delta: \nabla \xe \, \dx \dt+\int_{0}^T \int_{\Omega} p_\delta\opdiv\xe \, \dx \dt\\&\nonumber\quad - 	\int_{0}^T \int_{\Omega} \eps (\nabla \varphi_\delta  \otimes \nabla \varphi_\delta ): \nabla \xe \, \dx \dt- \int_\Omega \ve_0 \cdot \xe(\cdot, 0) \, \dx ,
	\end{align}
    for all $\xe \in C^\infty_c([0,T); C^1(\overline{\Omega}; \mathbb{R}^3 ))$ such that $\xe \cdot \no=0$ on $\partial \Omega \times [0,T)$, where $\mathbf v_\delta$ satisfies the incompressibility condition $\opdiv \mathbf v_\delta =0$ a.e. in $(0,T)\times \Omega$, and $J_m$ is a suitable symmetric smooth mollifier, parametrized by $m\in\N$.
    
\item[-]The approximate Allen-Cahn system: 
	\begin{align}
	&	\partial_t \varphi_\delta +  \ve_\delta \cdot \nabla \varphi_\delta  = - \frac1\eps \mu_\delta +\frac1k \Delta \widetilde\mu_\delta,\notag  \\
        &\mu_\delta =\widetilde \mu_\delta- \ell_\delta(\theta_\delta),
        \label{eq:weakACapprox}
	\end{align}
 {a.e. in } $\Omega \times (0,T)$, for $k\in\N$,
    where 
    $$
    \widetilde\mu_\delta:=- \eps  \Delta \varphi_\delta  +  \frac{1}{\eps} \partial_\varphi W(\varphi_\delta ).
    $$
  \item[-]  The approximate heat equation:
    \begin{align} 
  		&\int_{0}^T \int_\Omega q_\delta\partial_t \zeta \, \dx \dt
  		+ \int_{0}^T \int_\Omega q_\delta    \ve_\delta\cdot \nabla \zeta \, \dx \dt
  		\notag\\&\notag+   \int_{0}^T \int_\Omega  \hat \kappa_\delta(\theta_\delta)\Delta \zeta \, \dx \dt-\frac1m\int_0^T\int_\Omega \nabla q_\delta\cdot\nabla \zeta \, \dx\dt
  		 \\  \label{eq:weakheatapprox}
  		&= - \int_{\Omega}q_\delta(0) \zeta( \cdot, 0) \, \dx + \int_{0}^T \int_\Omega \ell_\delta(\theta_\delta) \frac1 \varepsilon \mu_\delta \chi_{\{q_{\delta}\geq0\}}\zeta \, \dx \dt
        \\&
        \quad - \int_{0}^T \int_\Omega  \Big( \nu(\theta_\delta)|\nabla \ve_\delta + \nabla \ve_\delta^\mathsf{T}|^2\Big)  \chi_{\{q_{\delta}\geq0\}} K_k(\zeta)\, \dx \dt 
         \notag\\&
        \quad - \int_{0}^T \int_\Omega \Big( 
  		\frac1 \varepsilon | \mu_\delta |^2 +\frac 1k |\nabla \widetilde\mu_\delta |^2  \Big)\chi_{\{q_{\delta}\geq0\}} K_k(\zeta)\, \dx \dt \notag
        ,
  	\end{align} 
  	for all $\zeta \in C^\infty_c( [0,T); C^2(\overline{\Omega}))$ such that $\nabla \zeta \cdot \no =0$ along $\partial \Omega $ and for any $ T >0$, which now will be solved for $q_\delta$ with initial datum $q_\delta(0)=Q_\delta(\theta_0)$. 
    Here, $\theta_\delta:=Q_\delta^{-1}((q_\delta)^++\delta)\in (Q_\delta^{-1}(\delta),+\infty)$ and $\hat \kappa_\delta (\theta_\delta)= \int_0^{\theta_\delta} \kappa_\delta(s) \, \d s = \kappa_1\theta_\eps  + \kappa_2 \frac{1}{\beta +1} \theta_\eps^{\beta +1}$ (see \eqref{eq:defkappa}). Moveover, for $k\in \N$, the cutoff function $K_k$ is defined as
\begin{align*}
K_k(s):=\begin{cases}
s,\quad\text{ if }s\in[-k,k],\\
k,\quad \text{if }s>k,\\
-k,\quad \text{if }s<-k.
\end{cases}
\end{align*}
 Note that the solution $q_\delta$ also depends on $m,k\in \N$. However, we omit this dependence in the next two steps for the sake of readability.
\end{itemize}
\subsubsection{Galerkin scheme.} In order to find a solution to the approximating problem \eqref{eq:weakbalmomapprox}-\eqref{eq:weakheatapprox}, for any $\delta>0$, we first further approximate it by means of a Galerkin scheme. 

Let $\mathbf V_n$ be the space generated by the first $n$ eigenfuctions $\{\bv_j\}_{j=1,\ldots,n}$ for the Stokes problem 
$$
\begin{cases}
    -\frac12 \div{}(\nabla \bv_j+\nabla \bv_j^T) +\nabla \pi_j=\lambda_j\bv_j,&\quad\text{ in }\Omega\\
    \operatorname{div} \bv_j=0,&\quad \text{ in }\Omega\\
    \bv_j\cdot\mathbf n_{\partial\Omega}=0,&\quad \text{ on }\partial\Omega
    \\(\operatorname{Id}-\mathbf n_{\partial\Omega}\otimes \mathbf n_{\partial\Omega})\frac12(\nabla \bv_j+\nabla \bv_j^T)=0,&\quad \text{ on }\partial\Omega,
\end{cases}
$$
which are orthonormal in $\mathbf L^2(\Omega)$ and orthogonal in $\mathbf H^1(\Omega)$. Let also $V_n$ be the space generated by first $n$ eigenfunctions $\{w_j\}_{j=1,\ldots,n}$ of the homogeneous Neumann Laplacian operator, which are orthonormal in $L^2(\Omega)$ and orthogonal in $H^1(\Omega)$. We also introduce the projectors $P_{\bV_n}$ and $P_{V_n}$ as the be the 
orthogonal projectors onto the spaces.  

We give the following \textit{ansatz} for the solution functions:
\begin{align*}
\bv_\delta^n:=\sum_{j=1}^n\alpha_j(t)\bv_j,\quad \vphi_\delta^n:=\sum_{j=1}^n\beta_j(t)w_j,\quad q_\delta^n:=\sum_{j=1}^n\gamma_j(t)\bv_j,
\end{align*}
where $\alpha_j, \beta_j, \gamma_j: [0,\infty) \rightarrow \R$. 
Moreover, we define \begin{align} \label{defthetadelta}
\theta_\delta^n=Q_\delta^{-1}((q_\delta^n)^++\delta)\in(Q_\delta^{-1}(\delta),+\infty), 
\end{align} 
which does not necessarily belong to $V_n$. We remark that $- \varepsilon \nabla \cdot  (\nabla \vphidn \otimes \nabla \vphidn ) = \widetilde\mu_\delta^n\nabla \vphidn - \nabla \Big( \frac{\varepsilon  }2 |\nabla \vphidn|^2+ \frac1 \varepsilon W(\vphidn ) \Big)$.
For almost every $t >0$, the triple $(\bv_\delta^n, \vphi_\delta^n, q_\delta^n)$  is assumed to solve:
\begin{align} \label{eq:weakbalmom1}
		&\int_{\Omega} (\partial_t\ve_\delta^n \cdot \xe + (\ve_\delta^n  \otimes \P_\sigma(J_m\ast\vdn)  ): \nabla \xe  
		) \, \dx
		+	\int_{\Omega} \mathbb{S}(\theta_\delta^n)\ve_\delta^n: \nabla \xe \, \dx \\ 
        & \notag = \int_{\Omega} \widetilde\mu_\delta^n\nabla \vphidn \cdot\xe \, \dx  , \quad    \forall \xe\in \bV_n,     \\& \notag \bv_\delta^n(0)=P_{\bV_n}\bv_0,\\&\,\notag\\&\label{varphi1}
        \eps\int_\Omega\partial_t \varphi_\delta^n\zeta \, \dx +  \eps\int_\Omega\ve_\delta^n\cdot \nabla \varphi_\delta^n\zeta \, \dx +\int_\Omega \mu_\delta^n\zeta \, \dx+\frac\eps k\int_\Omega \nabla\widetilde\mu_\delta^n\cdot\nabla \zeta \, \dx=0 ,\quad \forall \zeta\in V_n,\\&\notag
        \widetilde \mu_\delta^n=-\eps\Delta\vphi_\delta^n + \frac{1}{\eps} P_{V_n}\partial_\varphi W(\varphi_\delta ^n),
        \\&\notag
        \mu_\delta^n=\widetilde\mu_\delta^n - P_{V_n}(\ell_\delta(\theta_\delta^n)\chi_{\{\qdn\geq0\}}), 
        \\& \notag\vphi_\delta^n(0)=P_{V_n}\varphi_{\eps,0},\\&\,\nonumber\\&
        \label{eq:weakent1}
  		 \int_\Omega \partial_tq_\delta^n  \zeta \, \dx 
  		- \ \int_\Omega q_\delta^n   \ve_\delta^n\cdot \nabla \zeta \, \dx 
  		- \int_\Omega \hat \kappa_\delta(\theta_\delta^n)\Delta \zeta \, \dx +\frac1m\int_\Omega \nabla q_\delta^n\cdot\nabla \zeta \, \dx
  		 \\ 
  		&= -\frac1\eps \int_\Omega \ell_\delta(\theta_\delta^n)\mu_\delta^n\chi_{\{q_{\delta}^n\geq0\}}\zeta \, \dx  \non
        \\& \quad +\int_\Omega  \Big( {\nu(\thetadn)}|\nabla \ve_\delta^n + (\nabla \ve_\delta^n)^\mathsf{T}|^2\Big)  \chi_{\{q_{\delta}^n\geq0\}} K_k(\zeta)\, \dx  \notag 
          \\& \quad +\int_\Omega \Big(  
  		 \frac1\varepsilon \large| \mu_\delta^n \large|^2 +\frac 1k |\nabla \widetilde\mu_\delta^n |^2  \Big)\chi_{\{q_{\delta}^n\geq0\}} K_k(\zeta)\, \dx  \notag ,
         \quad  \forall \zeta\in V_n,\\&\notag q_\delta^n(0)=P_{V_n}Q_\delta(\vartheta_0).
	\end{align}
The local existence of a unique solution $(\bv_\delta^n, \vphi_\delta^n, q_\delta^n)$ can be established by means of the Cauchy-Lipschitz theorem for ODEs. 

As a next step, we derive the energy estimates, leading to the global existence of the solution. By testing equation \eqref{eq:weakbalmom1} with $\xe=\bv_\delta^n$, equation \eqref{varphi1} with $\zeta=\frac1\eps\widetilde\mu_\delta^n$, and equation \eqref{eq:weakent1} with $\zeta=\frac1{2k}\qdn$, we obtain the energy identity
\begin{align}   \nonumber
&\ddt \left(\frac12\norm{\vdn}^2_{\mathbf L^2(\Omega)}+\frac\eps2\norm{\nabla\vphidn}^2_{L^2(\Omega)}+\frac1\eps\int_\Omega W(\vphidn)\dx+\frac1{4k}\norm{\qdn}^2_{L^2(\Omega)}\right)
\\& \nonumber
+\int_\Omega \nu(\thetadn)\norma{\nabla \vdn+(\nabla\vdn)^T}^2\dx
+\frac1\eps\norm{\widetilde\mu_\delta^n}^2_{L^2(\Omega)}+\frac1k\norm{\nabla\widetilde\mu_\delta}^2
\\& \nonumber
+\frac1{2k}\int_\Omega \kappa_\delta(\thetadn)\nabla \thetadn\cdot\nabla \qdn\dx+\frac 1 {2km}\norm{\nabla \qdn}^2_{L^2(\Omega)}
\\& 
\label{enid}  
= \frac1{2k}\int_\Omega \Big( {\nu(\thetadn)}|\nabla \ve_\delta^n + (\nabla \ve_\delta^n)^\mathsf{T}|^2  + 
  		 \frac1\varepsilon \large| \mu_\delta^n \large|^2+\frac 1k |\nabla \widetilde\mu_\delta^n |^2 \Big) K_k((\qdn)^+)\, \dx
         \\& \notag
         \quad -  \frac1{2k\varepsilon} \int_\Omega \ell_\delta(\theta_\delta^n)\mu_\delta^n\chi_{\{q_{\delta}^n\geq0\}}q_\delta^n \, \dx +\frac1 {\varepsilon} \int_\Omega  \widetilde \mu_\delta^n\ell_\delta(\thetadn)\chi_{\{\qdn\geq0\}} \, \dx.   
\end{align}
By the definition of $\kappa_\delta$ and $\thetadn$, since $Q_\delta^{-1}$ is monotone nondecreasing, it holds 
\begin{align}
\frac1{2k}\non&\int_\Omega \kappa_\delta(\thetadn)\nabla \thetadn\cdot\nabla \qdn\dx=\frac1{2k}\int_\Omega \kappa_\delta(Q_\delta^{-1}((\qdn)^++\delta)(Q_\delta^{-1})'((\qdn)^++\delta)\norma{\nabla (\qdn)^+}^2\dx\\&=\frac1{2k}\int_\Omega\norma{\nabla F_\delta((\qdn)^++\delta)}^2\dx\geq 0,\label{Fdelta}
\end{align}
where
\begin{align} \label{defFdelta}
F_\delta(s):=\int_0^s\sqrt{\kappa_\delta(Q_\delta^{-1}(\tau))(Q_\delta^{-1})'(\tau)}\d\tau.
\end{align}
Recalling that $\norma{K_k(s)}\leq k$, $\forall s \in \R $, we have by standard inequalities
\begin{align*}
&\frac1{2k}\int_\Omega  \Big( {\nu(\thetadn)}|\nabla \ve_\delta^n + (\nabla \ve_\delta^n)^\mathsf{T}|^2  + 
  		 \frac1\varepsilon \large| \mu_\delta^n \large|^2  +\frac 1k |\nabla \widetilde\mu_\delta^n |^2 \Big) K_k((\qdn)^+)\, \dx
        \\&\leq \frac1{2}\int_\Omega \Big( {\nu(\thetadn)}|\nabla \ve_\delta^n + (\nabla \ve_\delta^n)^\mathsf{T}|^2 + 
  		\frac1\varepsilon \large| \mu_\delta^n \large|^2 +\frac 1k |\nabla \widetilde\mu_\delta^n |^2 \Big) \, \dx.
\end{align*}
 Notice that this estimate justifies the necessity of introducing the cutoff $K_k$ at the Galerkin approximation step. After this step, the cutoff is removed by letting $k\to\infty$.

Now, since $\ell_\delta(s)\leq C$ for any $s\in\R$, by Cauchy-Schwarz and Young's inequalities we infer 
\begin{align*}
\norma{  -\frac 1{2k\eps} \int_\Omega \ell_\delta(\theta_\delta^n)\mu_\delta^n\chi_{\{q_{\delta}^n\geq0\}}q_\delta^n \, \dx +\frac1\varepsilon  \int_\Omega \widetilde \mu_\delta^n\ell_\delta(\thetadn)\chi_{\{\qdn\geq0\}} \, \dx}\leq \frac C\varepsilon \left(1+\norm{\qdn}_{L^2(\Omega)}^2\right)+\frac{1}{2\varepsilon}\norm{\widetilde\mu_\delta^n}_{L^2(\Omega)}^2
\end{align*}
As a consequence, we obtain from \eqref{enid}, after an application of Gronwall's lemma, the following bounds:
\begin{align}
&\non\norm{\vdn}_{L^\infty(0,T;\mathbf L^2(\Omega))\cap L^2(0,T;\mathbf H^1(\Omega))}+\norm{\vphidn}_{L^\infty(0,T;H^1(\Omega))}\\&+\norm{\mu_\delta^n}_{L^2(0,T;H^1(\Omega))}+\norm{\qdn}_{L^\infty(0,T;L^2(\Omega))\cap L^2(0,T;H^1(\Omega))}\leq C_{\delta,m,k}(T),\label{uniformn}
\end{align}
for any $T>0$, which entails, by a standard ODE blow-up argument, that the maximal time of existence of the solution is $+\infty$.  
Moreover, from \eqref{Fdelta} we deduce
\begin{align}
\norm{F_\delta((q_\delta^n)^++\delta)}_{L^2(0,T;L^6(\Omega))}\leq C\norm{F_\delta((q_\delta^n)^++\delta)}_{L^2(0,T;H^1(\Omega))}\leq C_{\delta,m,k}.\label{Fdeltacontrol}
\end{align} 
Note that, from \eqref{defQdelta} it follows that, for $s\geq \delta>Q_\delta(0)$ (for $\delta>0$ sufficiently small), we have
$
Q_\delta^{-1}(s)=\left((1+\alpha)^{\frac{2}{1+\alpha}}s^{\frac 2{1+\alpha}}-\delta^2\right)^\frac12
$, thus $(Q_\delta^{-1})'(s) = (Q_\delta^{-1}(s))^{-1} (1+\alpha)^{\frac{1-\alpha}{1+\alpha}}s^{\frac{1-\alpha}{1+\alpha}}$.
Recalling that $\beta\geq2$, \eqref{defkappadelta} and \eqref{defFdelta}, for $s\geq \delta >Q_\delta(0)$, we compute
\begin{align*}
&\sqrt{{\kappa_2}({1+\alpha})^{\frac{1-\alpha}{1+\alpha}}}(1+\alpha)^{\frac{1}{1+\alpha}(\frac\beta2-\frac12)}s^{\frac{\beta+ 2\alpha+1}{2(1+\alpha)}}\\
&=\int_0^s\sqrt{{\kappa_2}
({1+\alpha})^{\frac{1-\alpha}{1+\alpha}}}
\tfrac{\beta + 2 \alpha +1}{2(1+\alpha)}
((1+\alpha)^{\frac{2}{1+\alpha}}\tau^{\frac 2{1+\alpha}})^{\frac\beta4-\frac14}\d\tau\\&
\leq C\int_0^s\sqrt{{\kappa_2}({1+\alpha})^{\frac{1-\alpha}{1+\alpha}}}((1+\alpha)^{\frac{2}{1+\alpha}}\tau^{\frac 2{1+\alpha}}-\delta^2)^{\frac\beta4-\frac14}\d\tau+Cs\delta^{\frac\beta2-\frac12}
\\&\leq 
C\int_0^s\sqrt{{\kappa_2}({1+\alpha})^{\frac{1-\alpha}{1+\beta}} (Q_\delta^{-1}(\tau))^{\beta -1} }\d\tau+Cs\delta^{\frac\beta2-\frac12}
\\&\leq CF_\delta(s)+Cs.
\end{align*}
Recalling \eqref{defthetadelta}, we observe that, for $\delta>0$ sufficiently small, $Q_\delta^{-1} (\delta)>0$, thus $\theta_\delta^n >0.$
Hence, using \eqref{uniformn} and \eqref{Fdeltacontrol}, we can infer that 
\begin{align}
\norm{{(\theta_\delta^n)}^{\frac{\beta}{2}+\frac{2\alpha+1}{2}}}_{L^2(0,T;L^6(\Omega))}\leq C\Big(1+\norm{((q_\delta^n)^++\delta)^{\frac{\beta}{2(1+\alpha)}+\frac{2\alpha+1}{2(\alpha+1)}}}_{L^2(0,T;L^6(\Omega))}\Big)\leq C_{\delta,m,k}(T),
\end{align}
entailing, 
\begin{align}
\norm{\theta_\delta^n}_{L^{\beta^+}(0,T;L^{3\beta^+}(\Omega))}\leq C_{\delta,m,k}(T),\label{boundbeta}
\end{align}
where $\beta^+:= {\beta}+{2\alpha+1} > \beta$.

Concerning time derivatives, by comparison in \eqref{eq:weakbalmom1}-\eqref{eq:weakent1}, using \eqref{boundbeta} and recalling that $\ell_\delta$ is bounded, as well as $\norma{K_k(s)}\leq k,$ $\forall s \in \R $, we obtain
\begin{align}
\norm{\partial_t\bv_\delta^n}_{L^2(0,T; \bV_\sigma')}+\norm{\partial_t\vphidn}_{L^2(0,T;(H^1(\Omega))')}+\norm{\partial_t\qdn}_{L^2(0,T;(H^2(\Omega))')}\leq C(T),\label{dt1}
\end{align}
for any $T>0$, where the constant here depends on $k\in\N$, $\delta>0$, and $m\in\N$. As a consequence, by standard compactness arguments, we deduce the following convergences as $n\to\infty$:
\begin{align*}
&\vdn\rightarrow \bv_\delta   &&\text{weakly in }L^2_{loc}([0,\infty);\bV_\sigma) , \\&
\vdn \rightarrow  \bv_\delta  &&\text{weakly* in }L^\infty_{loc}([0,\infty);\mathbf H_\sigma),
\\&
\partial_t\vdn\rightarrow\partial_t\bv_\delta &&\text{weakly in }L^2_{loc}([0,\infty); \bV_\sigma '),
\\&
\vdn\to \bv_\delta  &&\text{strongly in }L^2_{loc}([0,\infty);\mathbf H_\sigma)\text{ and a.e. in } \Omega\times(0,\infty),
\\&
\vphidn \rightarrow \varphi_\delta 
 &&\text{weakly* in }L^\infty_{loc}([0,\infty);H^1(\Omega)),
\\&
\partial_t\vphidn \rightarrow \partial_t\vphi_\delta  &&\text{weakly in }L^2_{loc}([0,\infty);(H^1(\Omega))'),
\\&
\vphidn\to \vphi_\delta  && \text{strongly in }L^2_{loc}([0,\infty);L^2(\Omega)) \text{ and a.e. in } \Omega\times(0,\infty),\\&
\mu_\delta^n \rightarrow \mu_\delta  && \text{weakly in }L^2_{loc}([0,\infty);H^1(\Omega)),\\&
\qdn \rightarrow  q_\delta  && \text{weakly in }L^2_{loc}([0,\infty);H^1(\Omega)),\\&
\qdn \rightarrow  q_\delta  && \text{weakly* in }L^\infty_{loc}([0,\infty);L^2(\Omega)),\\&
\partial_t\qdn  \rightarrow  \partial_tq_\delta  &&\text{weakly in }L^2_{loc}([0,\infty);(H^2(\Omega))'),\\&
\qdn \to q_\delta &&\text{strongly in }L^2_{loc}([0,\infty);L^2(\Omega))\text{ and a.e. in } \Omega\times(0,\infty).
\end{align*}
Since the function $Q_\delta^{-1}$ is continuous, by means of the Vitali convergence theorem, from \eqref{boundbeta} we can infer
\begin{align*}
\thetadn\to \theta_\delta:=Q_\delta^{-1}((q_\delta)^++\delta) \quad\text{ strongly in }L^\beta_{loc}([0,\infty);L^{2\beta}(\Omega))\text{ and a.e. in } \Omega\times(0,T),
\end{align*}
as $n\to\infty$, for any $T>0$, where $\theta_\delta \geq 0$ a.e. in  $\Omega\times(0,T) $.

Observe that these convergences, thanks to the regularization parameters, are enough to pass to the limit \eqref{eq:weakbalmom1}, \eqref{varphi1}, and \eqref{eq:weakent1}, obtaining that the triple $(\bv_\delta,\vphi_\delta,\theta_\delta)$ solves \eqref{eq:weakbalmomapprox} (for divergence free test functions), \eqref{eq:weakACapprox}, and \eqref{eq:weakheatapprox}. 
Furthermore, by lower semicontinuity of the norms, we can pass to the limit in the energy identity \eqref{enid} (after an integration in time).
Finally, we introduce the pressure function $p_\delta$ by solving the ``very weak" formulation: 
	\begin{align}\label{eq:weakp}
		\int_{\Omega} p_\delta \Delta \eta \, \dx = \int_{\Omega} (\nu(\theta_\delta)D\bv_\delta- \eps \nabla \varphi_\delta \otimes \nabla \varphi _\delta - \ve_\delta \otimes \P_\sigma(\ve_\delta)) : \nabla \nabla \eta \, \dx   , 
	\end{align}
	for any $\eta \in \mathbf H^2(\Omega)$ with $\nabla \eta \cdot \no = {0}$ along ${\partial \Omega}$, for almost any $t>0$. The existence of a solution such that $\int_\Omega p_\delta \, \dx =0$ can be obtained from the regularity of the right-hand side. Therefore, with this formulation at hand, we can see that \eqref{eq:weakbalmomapprox} holds for any $\xe \in C^\infty_c([0,T); C^1(\overline{\Omega}; \mathbb{R}^3 ))$ such that $\xe \cdot \no=0$ on $\partial \Omega \times [0,T)$, without requiring $\xe$ to be divergence free. 

    \subsubsection{Approximate total energy balance and approximate entropy equation}
    Note that $\bv_\delta\in AC([0,T];\H_\sigma)$, for any $T>0$, hence we can test \eqref{eq:weakbalmomapprox} with $\xe=\bv_\delta\eta$, where $\eta\in C^\infty_c([0,T); C^2(\overline{\Omega}; \mathbb{R} ))$ with $\nabla \eta  \cdot \no=0$ on $\partial \Omega \times [0,T)$, as well as we can multiply \eqref{eq:weakACapprox} by $\frac1\eps\eta\mu_\delta$, where $\mu_\delta=-\eps\Delta\vphi_\delta+\frac1\eps\partial_\vphi W(\vphi_\delta)-\frac1\eps\ell_\delta(\theta_\delta)$ and integrate space-time, and choose $\zeta=\eta$ in \eqref{eq:weakheatapprox}.  
Notice that, since we introduced the approximating term $\frac1k \Delta\widetilde \mu_\delta$ into the scheme, we have
$\widetilde\mu_\delta\in L^2(0,T;H^1(\Omega))$, which in turn guarantees that $\eps \nabla \varphi_\delta \otimes \nabla \varphi_\delta $ is regular enough to be tested with the gradient of $\bv_\delta \eta$.
Summing up the resulting equations, we obtain 
    \begin{align}   \label{eq:weaktoten2b}
		&\int_{0}^T \int_{\Omega}  \Big(e_\delta + \frac12 |\ve_\delta|^2\Big) \partial_t \eta \, \dx \dt
		+	\int_{0}^T \int_{\Omega}  \Big(e_\delta \ve_\delta + \frac12 |\ve_\delta|^2\mathbf P_\sigma (J_m\ast\ve_\delta) \Big)\cdot \nabla \eta  \, \dx \dt 
        \\& \notag 
		+ \int_{0}^T \int_{\Omega} (\hat\kappa_\delta(\theta_\delta  )+\frac1m q_\delta) \Delta \eta \, \dx \dt	
		+ 	\int_{0}^T \int_{\Omega} p_\delta \ve_\delta  \cdot \nabla \eta \, \dx \dt
         \\& \notag 
		- \int_{0}^T \int_{\Omega} (\mathbb{S}(\theta_\delta)\ve_\delta ):  (\ve_\delta \otimes \nabla \eta ) \, \dx \dt
        + \frac1k\int_{0}^T \int_{\Omega} |\mu_\delta|^2 \Delta \eta
		\, \dx \dt
        - \frac1k\int_{0}^T \int_{\Omega} \ell(\theta_\delta) \nabla \tilde \mu_\delta \cdot \nabla \eta
		\, \dx \dt
         \\& \notag 
         -\int_0^T\int_\Omega \Big({\nu(\theta_\delta)}|\nabla \ve_\delta + (\nabla \ve_\delta)^\mathsf{T}|^2 + 
  		 \frac1\varepsilon | \mu_\delta |^2+\frac1k | \nabla\widetilde\mu_\delta |^2  \Big)\eta \, \dx \dt\\&
         +\int_0^T\int_\Omega \Big( {\nu(\theta_\delta)}|\nabla \ve_\delta + (\nabla \ve_\delta)^\mathsf{T}|^2 
  	+	 \frac1\varepsilon | \mu_\delta |^2 +\frac1k | \nabla\widetilde \mu_\delta |^2 \Big)\chi_{\{q_{\delta}\geq0\}} K_k(\eta)\, \dx \dt
        \notag\\& 
		+ \int_{\Omega}  \Big(e_0 + \frac12 |\ve_0|^2\Big) \eta(\cdot, 0)\, \dx =0,\notag
	\end{align}
	where $e_\delta=\frac\eps2\norma{\nabla \vphi_\delta}^2+\frac1\eps W(\vphi_\delta)+q_\delta$ and $\hat \kappa_\delta (\theta_\delta)= \int_0^{\theta_\delta} \kappa_\delta(s) \d s$,
    for all $\eta \in C^\infty_c([0,T); C^2(\overline{\Omega}))$ such that 
    $\nabla \eta  \cdot \no=0$ on $\partial \Omega \times [0,T)$.
    
Observe that we can test \eqref{eq:weakheatapprox} by $\zeta=-(q_\delta)^-$, hence, after some immediate integration by parts and recalling that $\mathrm{supp}\ \theta_\delta=\{q_\delta\geq0\}$, we can obtain
\begin{align*}
\frac 12 \norm{(q_\delta)^-( \cdot, t) }^2_{L^2(\Omega)}+\frac1m \int_0^t \norm{\nabla (q_\delta)^-}_{L^2(\Omega)}^2 \, \dt \leq \frac 12 \norm{(q_\delta)^-( \cdot, 0) }^2_{L^2(\Omega)},
\end{align*}
for almost all $t \in (0, \infty)$.
Since $q_\delta(0)=Q_\delta(\theta_0)\geq0$, we can  
\begin{align}
q_\delta\geq 0 \quad \text{ almost everywhere in }\Omega\times(0,\infty),
\label{negat}
\end{align}
so that $\theta_\delta=Q_\delta^{-1}(q_\delta+\delta)$, thus $q_\delta=Q_\delta(\theta_\delta)+\delta$, and $\chi_{\{q_\delta\geq0\}}\equiv 1$ almost everywhere in $\Omega\times(0,\infty)$. This means that we can substitute $ q_\delta$ (together with its derivatives) in \eqref{eq:weakheatapprox} with $ Q_\delta(\theta_\delta)$, noting that the initial condition for $\theta_\delta$ becomes $\theta_\delta(0)=Q_\delta^{-1}(Q_\delta(\theta_0)+\delta)\to \theta_0$ almost everywhere in $\Omega$, as $\delta\to0$, as expected.
Moreover, observe that one can even show that  $\theta_\delta \in L^2_{loc}(0, \infty; H^1(\Omega))$, for instance, by testing \eqref{eq:weakheatapprox} with $\zeta=  \theta_\delta$, although this estimate would depend on $\delta$.

Furthermore, by testing \eqref{eq:weakheatapprox} with $\zeta=\frac1 {\ell_\delta(\theta_\delta)}\eta$ (recall that by construction $\ell_\delta>0$ on $\R$), we obtain the entropy balance:
  	\begin{align} \label{eq:weakent2} 
  		&\int_{0}^T \int_\Omega(\Lambda_\delta(\theta_\delta) + \varphi_\delta) (\partial_t \eta + \ve_\delta \cdot \nabla \eta) \, \dx \dt
     + \int_{0}^T \int_\Omega (h_\delta(\theta_\delta )+h_{m,\delta}(\theta_\delta)) \Delta \eta \, \dx \dt
  		 \\ \notag
  		&= - \int_{0}^T \int_\Omega \bigg( \frac{K_k(\eta)}{\ell_\delta(\theta_\delta)} {\nu(\theta_\delta)}|\nabla \ve_\delta + \nabla \ve_\delta^\mathsf{T}|^2 + \frac{\ell^\prime_\delta (\theta_\delta)}{\ell^2_\delta(\theta_\delta)}  
        \Big(\kappa_\delta(\theta_\delta)+\frac1m Q'_\delta(\theta_\delta)\Big)  |\nabla \theta_\delta|^2
  		 \bigg) \eta\, \dx \dt   \\ \notag 
        & \quad -  \int_{0}^T \int_\Omega  
  		 \frac{1}{\varepsilon\ell_\delta(\theta_\delta)} |\mu_\delta |^2  K_k(\eta)\, \dx \dt -  \int_{0}^T \int_\Omega  
  		 \frac{1}{k\ell_\delta(\theta_\delta)} |\nabla \widetilde{\mu}_\delta |^2  K_k(\eta)\, \dx \dt  
          \\ \notag
        & \quad  - \int_{\Omega}(\Lambda_\delta(\theta_\delta(0)) +  \varphi_{\eps,0}) \eta( \cdot, 0) \, \dx , \notag
  	\end{align}
  	where 
    \begin{align}
    \label{Lambdadelta}
    \Lambda_\delta(s):=\int_0^s \frac{Q_\delta' (\tau)}{\ell_\delta(\tau)}\d\tau \quad \text{and} \quad h_\delta(\theta_\delta):= \int_{1}^{\theta_\delta} \frac{\kappa_\delta(s)}{\ell_\delta(s)} \, \d s,
    \end{align}  
    for all $\eta \in C^\infty_c(\overline{\Omega} \times [0,T))$ such that $\eta \geq 0$ and {$\nabla \eta \cdot \no =0$ along $\partial \Omega$}, for any $ T >0$. Here, we also set $$h_{m,\delta}(s):=\frac1m\int_1^s\frac{Q_\delta(\tau)}{\ell_\delta(\tau)}\d\tau\leq C(\norma{s}+\int_1^s{Q_\delta(\tau)}\d\tau)$$ by the construction of $Q_\delta,\ell_\delta$ and the assumptions on $\ell$ (cf. Section \ref{subsec:hp}). Notice that all the functions $Q_\delta, \ell_\delta$ are Lipschitz continuous, hence the identities above hold due to the differentiation formula for the composition of Lipschitz functions with $H^1(\Omega)$-functions.

\subsection{{Limit as $k\to\infty$ and $\delta\to 0$.} } \label{subsec:delta}
This section concerns the limits $k\to\infty$ and $\delta\to 0$.
We first derive some regularity estimates directly following from the total energy balance \eqref{eq:weaktoten2b}, the entropy balance \eqref{eq:weakent2}, and the heat equation \eqref{eq:weakent1}. 
Then, since the most significant limit passage is that as $\delta\to0$, we only comment on the limit 
$k\to\infty$ (which is, in fact, very similar) and focus instead on the case $\delta\to0$.
After letting $k\to\infty$, we can apply the maximum principle to show the uniform boundedness of $\varphi_\delta$ as $\delta\to0$, as well as that the strict positivity of $\theta_\delta$ is preserved in the limit as $\delta\to0$. Finally, thanks to the convergence properties, we can show the weak sequential stability of the problem as $\delta\to 0$.
 
 \subsubsection{Energy estimates} \label{SubsecEnEntEst}
Consider the total energy balance \eqref{eq:weaktoten2b} and test it by a sequence of functions approximating the identity, integrating by parts in time and then passing to the limit  (i.e., $\eta\equiv1$, so that $K_k(\eta)=\eta$), one can deduce the energy inequality
\begin{align}
 &\int_\Omega  e_\delta(T)+\frac12\norma{\ve_\delta(T)}^2\dx\leq \int_\Omega e_\delta(0)+\frac12\norma{\ve_\delta(0)}^2\dx, \label{energineqA}
\end{align}
for any $T>0$.
From this inequality, we can infer the following estimates (all the constants appearing here are independent of all the approximating parameters $k,\delta$, and $m$):
\begin{align}\label{E1}
	\|\ve_\delta\|_{L^\infty(0,T; \LL^2(\Omega))}&\leq C,\\
	\|Q(\theta_\delta)\|_{L^\infty(0,T; L^1(\Omega))} &\leq C , \label{eq:boundQinfty} \\
	\|W(\varphi_\delta)\|_{L^\infty(0,T; L^1(\Omega))} &\leq  C \eps, \\
	\sqrt{\eps }\|\nabla \varphi_\delta \|_{L^\infty(0,T; \mathbf L^2(\Omega))} &\leq C ,
\end{align}
whence 
\begin{align}
\|\theta_\delta\|_{L^\infty(0,T; L^{1+\alpha}(\Omega))} &\leq C, \\
\|\varphi_\delta \|_{L^\infty(0,T; L^4(\Omega))} &\leq C,\\
\sqrt{\eps } \|\varphi_\delta \|_{L^\infty(0,T; H^1(\Omega))} &\leq C.
\end{align}
Using the classical Sobolev embedding theorem, we deduce
\begin{align}
	\sqrt{\eps } \|\varphi_\delta \|_{L^\infty(0,T; L^6(\Omega))} &\leq C.
\end{align}
and therefore
\begin{equation}
\sqrt{\eps }	\|\partial_\varphi W(\varphi_\delta)\|_{L^\infty(0,T; L^2(\Omega))} \leq C. \label{eq:estWprime}
\end{equation}
Additionally, since $K_k(\eta)\leq \eta$ whenever $\eta\geq0$, we can obtain from \eqref{eq:weaktoten2b}
\begin{align}   \label{eq:weaktoten2bf}
		&\int_{0}^T \int_{\Omega}  \Big(e_\delta + \frac12 |\ve_\delta|^2\Big) \partial_t \eta \, \dx \dt
		+	\int_{0}^T \int_{\Omega}  \Big(e_\delta \ve_\delta + \frac12 |\ve_\delta|^2\mathbf P_\sigma (J_m\ast\ve_\delta) \Big)\cdot \nabla \eta  \, \dx \dt 
        \\& \notag 
		+ \int_{0}^T \int_{\Omega} (\hat\kappa_\delta(\theta_\delta  )+\frac1m q_\delta) \Delta \eta \, \dx \dt	
		+ 	\int_{0}^T \int_{\Omega} p_\delta \ve_\delta  \cdot \nabla \eta \, \dx \dt
         \\& \notag 
		- \int_{0}^T \int_{\Omega} (\mathbb{S}(\theta_\delta)\ve_\delta ):  (\ve_\delta \otimes \nabla \eta ) \, \dx \dt
         + \frac1k\int_{0}^T \int_{\Omega} |\mu_\delta|^2 \Delta \eta
		\, \dx \dt
        - \frac1k\int_{0}^T \int_{\Omega} \ell(\theta_\delta) \nabla \tilde \mu_\delta \cdot \nabla \eta
		\, \dx \dt
         \\& +\int_{\Omega} \Big(e_0 + \frac12 |\ve_0|^2\Big) \eta(\cdot, 0)\, \dx \geq 0,\notag
	\end{align}
 for all $\eta \in C^\infty_c([0,T); C^2(\overline{\Omega}))$ such that $\eta \geq 0$ and
    $\nabla \eta  \cdot \no=0$ on $\partial \Omega \times [0,T)$.
 
 \subsubsection{Entropy estimates depending on the latent heat $\ell$} \label{subsec:estl}
 We recall that we consider the class of sublinear and bounded latent heat functions complying with the assumptions stated in Section \ref{subsec:hp}. 

First, using the definitions of $\ell_\delta$, $Q_\delta$ and $\Lambda_\delta$ (cf. \eqref{elllambda}, \eqref{defQdelta} and \eqref{Lambdadelta}), as well as the assumptions (d) and (e) on $\ell$ listed in Section \ref{subsec:hp}, we deduce  
\begin{itemize}
\item for $0\leq \theta_\delta\leq  \ell^{-1} (\delta)$,
\begin{align} 
	\non	\Lambda_\delta(\theta_\delta) &= \int_{0}^{\theta_\delta} \frac{Q^\prime_\delta(s)}{\ell_\delta(s)} \,  \mathrm{d}s =\int_0^{\theta_\delta}
    \frac1{\delta^2}s(s^2+\delta^2)^{\frac{\alpha-1}2}(2\delta-\ell(s)) \, \d s
		\\&
        \leq \frac2{\delta}\int_0^{\theta_\delta} s (s^2+\delta^2)^{\frac{\alpha-1}2} \,  \mathrm{d}s = 
        \frac2{\delta(1+\alpha)}[({\theta_\delta}^2+\delta^2)^{\frac{1+\alpha}2}-\delta^{1+\alpha}]\non\\&
         \leq \frac2{\delta(1+\alpha)}[((\delta_0 \ell^\prime(0) )^{-1}\ell({\theta_\delta})^2+\delta^2)^{\frac{1+\alpha}2}-\delta^{1+\alpha}]
        \notag  \\ & \leq \frac{2\delta^\alpha}{(1+\alpha)}[((\delta_0 \ell^\prime(0) )^{-1}+1)^{\frac{1+\alpha}2}-1]\leq C,
		 \label{eqn:boundSQ01}
	\end{align} 
    \item for $\ell^{-1} (\delta)<\theta_\delta \leq 1$, 
	\begin{align} 
		\Lambda_\delta(\theta_\delta) &= \Lambda_\delta(\ell^{-1}(\delta))+\int_{\ell^{-1}(\delta)}^{\theta_\delta}\frac{Q_\delta^\prime(s)}{\ell(s)} \,  \mathrm{d}s 
		\leq 
        C+\frac1{\delta_0 \ell'(0)}\int_{\ell^{-1}(\delta)}^{\theta_\delta} \frac{Q^\prime_\delta(s)}{s} \,  \mathrm{d}s \non\\&\leq \frac1{\delta_0\ell'(0)}\int_{\ell^{-1}(\delta)}^{1} (s^2+\delta^2)^{\frac{1+\alpha}2} \,  \mathrm{d}s 
     \leq C,
		 \label{eqn:boundSQ0}
	\end{align} 
	\item for $\theta_\delta >  1$,   
	\begin{align} 
		\Lambda_\delta(\theta_\delta) &= \Lambda_\delta(1) + \int_{1}^{\theta_\delta} \frac{Q^\prime_\delta(s)}{\ell(s)} \,  \mathrm{d}s \notag \\
		&	\leq C + \frac1{\ell(1)} \int_{1}^{\theta_\delta} {Q^\prime_\delta(s)}\,  \mathrm{d}s = C + \frac1{\ell(1)} (Q_\delta(\theta_\delta)-Q_\delta(1)).\label{eqn:boundSQ}
	\end{align} 
\end{itemize}

   Second, we consider the entropy balance \eqref{eq:weakent2} and test it by a sequence of functions approximating the identity. Integrating by parts in time and then passing to the limit  (i.e., $\eta\equiv1$, recalling $K_k(1)=1$), we infer
  \begin{align}\label{eq:stimakappa}
    \int_{0}^T \int_\Omega   \frac{\ell^\prime_\delta (\theta_\delta)}{\ell^2_\delta(\theta_\delta)}  
        \kappa_\delta(\theta_\delta)  |\nabla \theta_\delta|^2
  		\, \dx \dt \leq C.
   \end{align}
   
   Finally, we prove some additional regularity bounds for $\theta$ following from \eqref{eq:stimakappa}, namely
  \begin{align}
  \label{eq:boundtheta1}	\|\nabla \theta_\delta\|_{L^2(0,T; \LL^2(\Omega))} &\leq C, \\
  \label{eq:boundtheta2}	||\nabla \log_\delta(\theta_\delta)\|_{L^2(0,T; \LL^2(\Omega))} &\leq C ,\\
\label{eq:boundtheta3}  	\|\nabla \theta_\delta^{\beta/2}\|_{L^2(0,T; \LL^2(\Omega))} &\leq C ,
  \end{align}
  where $\log_\delta$ is defined in \eqref{logdelta}. As a consequence, we obtain
  \begin{equation*}
  \|\theta_\delta^{\beta/2}\|_{L^2(0,T; L^6(\Omega))} \leq C \implies \|\theta_\delta\|_{L^\beta(0,T; L^{3\beta}(\Omega))} \leq C. \label{eq:boundthetabeta0}
  \end{equation*}
Observe that  there exists $C_\beta\in(0,1)$ such that
$$
\frac{1+ s^\beta}{1+s^\gamma }\geq C_\beta s^{\beta-\gamma},\quad \forall s \geq 1,
$$
whenever $\gamma \leq \beta$.
For instance, choosing $C_\beta=\tfrac12$, we have
$
1+s^\beta\geq s^\beta\geq  \frac12 s^{\beta-\gamma}+\frac12 s^\beta, \, \forall s \geq1.
$
Hence, using the assumptions (a) and (e) on $\ell $ listed in Section \ref{subsec:hp}, as well as the fact that $\ell_\delta(s)=\ell(s)$ for $s\geq1$ (see \eqref{elllambda}), we can infer from \eqref{eq:weakent} that
\begin{align*} 
	\int_{0}^{T}  \int_{\theta_\delta \geq1}	\frac{\ell^\prime_\delta (\theta_\delta)}{\ell_ \delta^2(\theta_\delta)} \kappa_\delta(\theta_\delta) |\nabla \theta_\delta|^2 \, \dx \mathrm{d}t &\geq \frac{C}{L^2} \int_{0}^{T}  \int_{\theta_\delta\geq 1}	\frac{1+ \theta_\delta^\beta}{1+\theta_\delta^\gamma }  |\nabla \theta_\delta|^2 \, \dx \mathrm{d}t\\
	& \geq  \frac{C}{L^2} \int_{0}^{T}  \int_{\theta_\delta\geq 1}	\theta_\delta^{\beta-\gamma} |\nabla \theta_\delta|^2 \, \dx \mathrm{d}t
	.
\end{align*}
On the other hand, using the assumption (b) and (c) on $\ell $ from Section \ref{subsec:hp}, recalling \eqref{logdelta2} and the definition of $\kappa_\delta$ we get
\begin{align*} 
	\int_{0}^{T}  \int_{\delta<\theta_\delta<1 }	\frac{\ell^\prime _\delta(\theta_\delta)}{\ell_\delta^2(\theta_\delta)} \kappa_\delta(\theta_\delta) |\nabla \theta_\delta|^2 \, \dx \mathrm{d}t&\geq C\int_{0}^{T}  \int_{\delta<\theta_\delta<1 }	\frac{\kappa_1+\kappa_2 \theta_\delta^\beta}{\theta_\delta^{2}}  |\nabla \theta_\delta|^2 \, \dx \mathrm{d}t\\&
    \geq C\min\{\kappa_1,\kappa_2\}\int_{0}^{T}  \int_{\delta<\theta_\delta<1 }\frac{1+ \theta_\delta^\beta}{\theta_\delta^{2}}  |\nabla \theta_\delta|^2 \, \dx \mathrm{d}t
    . 
\end{align*}
Additionally, for the remaining case $\theta_\delta<\delta$, recalling again \eqref{logdelta2} and the definition of $\kappa_\delta$, we have 
\begin{align*} 
	\int_{0}^{T}  \int_{\theta_\delta<\delta }	\frac{\ell^\prime _\delta(\theta_\delta)}{\ell_\delta^2(\theta_\delta)} \kappa_\delta(\theta_\delta) |\nabla \theta_\delta|^2 \, \dx \mathrm{d}t&\geq \frac{C}{\delta^2} \int_{0}^{T}  \int_{\theta_\delta<\delta }	({\kappa_1+ \delta^2\kappa_2\norma{\theta_\delta}^{\beta-2}}) |\nabla \theta_\delta|^2 \, \dx \mathrm{d}t\\&
    \geq {C} \int_{0}^{T}  \int_{\theta_\delta<\delta }	\norma{\nabla \log_\delta(\theta_\delta)}^2+ |\nabla {\theta_\delta}^{\beta/2}|^2 \, \dx \mathrm{d}t. 
\end{align*}
Therefore, for $\beta\geq \gamma$ we can deduce
\begin{align*} 
	\int_{0}^{T}  \int_{\Omega } |\nabla \theta_\delta|^2 \, \dx \mathrm{d}t &\leq \int_{0}^{T}  \int_{\theta_\delta\geq 1}	\theta_\delta^{\beta-\gamma} |\nabla \theta_\delta|^2 \, \dx \mathrm{d}t +\int_{0}^{T}  \int_{\theta_\delta< 1}\min\Big\{\frac{1}{\theta_\delta^{2}},\frac1{\delta^2}\Big\}  |\nabla \theta_\delta|^2 \, \dx \mathrm{d}t\leq C,
\end{align*}
and 
\begin{align} 
\non	&\int_{0}^{T}  \int_{\Omega }	 |\nabla \log_\delta\theta_\delta|^2 \, \dx \mathrm{d}t \\&\leq C \int_{0}^{T}  \int_{\theta_\delta \geq 1}	 |\nabla \theta_\delta |^2 \, \dx \mathrm{d}t +\int_{0}^{T}  \int_{\delta<\theta_\delta<  1}	\frac{1}{\theta_\delta^2}  |\nabla \theta_\delta|^2 \, \dx \mathrm{d}t+ \int_{0}^{T}  \int_{\theta_\delta<  \delta} |\nabla \log_\delta(\theta_\delta)|^2 \, \dx \mathrm{d}t\leq C.\label{logdelta3}
\end{align}
Furthermore, for $\beta \geq 2 \geq \gamma$ we have
\begin{align*} 
	\int_{0}^{T}  \int_{\Omega } |\nabla \theta_\delta^{\beta/2}|^2 \, \dx \mathrm{d}t 
	&	= \int_{0}^{T}  \int_{\Omega } \theta_\delta^{\beta-2}|\nabla \theta_\delta|^2 \, \dx \mathrm{d}t \\
	&\leq \int_{0}^{T}  \int_{\theta_\delta\geq 1}	\theta_\delta^{\beta-\gamma} |\nabla \theta_\delta|^2 \, \dx \mathrm{d}t +\int_{0}^{T}  \int_{\theta_\delta< 1}	\theta_\delta^{\beta-2} |\nabla \theta_\delta|^2 \, \dx \mathrm{d}t\leq C.
\end{align*}

 \subsubsection{Heat estimates}  \label{SubsecHeatEst} 
Consider \eqref{eq:weakheatapprox} and test it by a sequence of functions approximating the identity, integrate by parts in time and then pass to the limit  (i.e., $\eta\equiv1$, note that $K_k(\eta)=\eta$ in this case). 
Using the Young inequality, we can obtain
\begin{align}
&\int_0^T\int_\Omega \Big( {\nu(\theta_\delta)}|\nabla \ve_\delta + (\nabla \ve_\delta)^\mathsf{T}|^2  
  	+	 \frac1\varepsilon | \mu_\delta |^2 +\frac1k | \nabla\widetilde \mu_\delta |^2 \Big) \, \dx \dt\notag\\&\leq \frac1 {2\varepsilon}\int_{0}^{T} \int_{\Omega} \left| \mu_\delta \right|^2\dx\dt+ \norm{q_\delta}_{L^\infty(0,T;L^1(\Omega))}+C\notag\\&
\leq \frac1 {2\varepsilon}\int_{0}^{T} \int_{\Omega} \left| \mu_\delta \right|^2\dx\dt+C, \label{est:heat1}
\end{align}
 due to the control on $e_\delta$ given by the energy inequality \eqref{energineqA}. Here, the constant $C>0$ is independent of $k,\delta,$ and $m$.

 As a consequence of this inequality, we have
 \begin{align}
 \norm{\ve_\delta}_{L^2(0,T;\mathbf H^1(\Omega))}+\frac1{\sqrt\eps}\norm{\mu_\delta}_{L^2(0,T;L^2(\Omega))}\leq C,
 \label{evelox}
 \end{align}
 together with 
 \begin{align}
 \frac1{\sqrt k }\norm{\nabla\widetilde\mu_\delta}_{L^2(0,T;\mathbf L^2(\Omega))}&\leq C,\label{restimate}
 \end{align}
 for $C>0$ independent of $k,\delta,m$.
Finally, by comparison in \eqref{eq:weakACapprox}, using \eqref{est:heat1}, we obtain
\begin{equation}
\sqrt{\eps } \|	\Delta \varphi_\delta \|_{L^2(0,T; L^2(\Omega))} \leq C,  \quad \text{whence} \quad  \sqrt{\eps }  \|\varphi_\delta \|_{L^2(0,T; H^2(\Omega))} \leq C,  \label{eq:estvarphiH2}
\end{equation}
due to the standard elliptic regularity theory.

\subsubsection{Further estimates}\label{SubsecFurtherPB}   
One can show that
\begin{align}
	L^\infty(0,T; L^2(\Omega))	\cap L^2(0,T; H^1(\Omega))  \hookrightarrow
	L^\frac{10}3(0,T; H^\frac{3}{5}(\Omega))	 
\hookrightarrow
	L^\frac{10}3(0,T; L^{\frac{10}{3}}(\Omega))	. \label{eq:emb}
\end{align}
Hence, from the estimates above we deduce
\begin{align}
	& \sqrt{\eps } \|\nabla \varphi_\delta \|_{L^{\frac{10}3}(0,T; \LL^{\frac{10}3}(\Omega))} + 
	\|\ve_\delta \|_{L^{\frac{10}3}(0,T; \LL^{\frac{10}3}(\Omega))}+	\|\theta_\delta \|_{L^{\frac{10}3}(0,T; L^{\frac{10}3}(\Omega))} \leq C. \label{eq:estfur1}
\end{align}   
Additionally, we have
\begin{align}
	& \eps \|\nabla \varphi_\delta  \otimes  \nabla \varphi_\delta  \|_{L^{\frac53}(0,T; \LL^{\frac53}(\Omega))} + \|\ve_\delta  \otimes  \ve_\delta  \|_{L^{\frac53}(0,T; \LL^{\frac53}(\Omega))} + \sqrt{\eps } \|\ve_\delta  \cdot \nabla \varphi_\delta  \|_{L^{\frac53}(0,T; L^{\frac53}(\Omega))} \leq C, \label{eq:estfur2}\\
	&\eps \|(\nabla \varphi_\delta  \otimes \nabla \varphi_\delta ) \ve_\delta  \|_{L^{\frac{10}9}(0,T; \LL^{\frac{10}9}(\Omega))} + \eps \| |\ve_\delta |^2 \ve_\delta  \|_{L^{\frac{10}9}(0,T; \LL^{\frac{10}9}(\Omega))} \leq C,  \label{eq:estfur3} \\
	& \|Q_\delta(\theta_\delta ) \|_{L^{\frac{10}{3(1+\alpha)}}(0,T; L^{\frac{10}{3(1+\alpha)}}(\Omega))} \leq C, \quad \|Q_\delta(\theta_\delta )  \ve_\delta  \|_{L^{\frac{10}{3(2+\alpha)}}(0,T; \LL^{\frac{10}{3(2+\alpha)}}(\Omega))} \leq C , \label{eq:estfur4}\\&
    \|Q_\delta(\theta_\delta ) \|_{L^{2}(0,T; H^1(\Omega))} \leq C(m), \label{Qdeltam}
\end{align}
where the last constant might depend on $m \in \N $., but does not depend on $k\in\N$ nor $\delta>0$.

 By means of \eqref{E1} and \eqref{eq:estfur3}, we then infer from \eqref{eq:weakp} the estimate for the pressure $p_\delta$, satisfying 
 \begin{align}\label{pressure1}
\norm{p_\delta}_{L^\frac53(0,T;L^\frac53(\Omega))}\leq C.
 \end{align}
 As an estimate depending on $m\in \N$ (but uniform in $\delta>0$), one can infer by comparison in \eqref{eq:weakbalmomapprox}-\eqref{eq:weakheatapprox}
\begin{align}
\norm{\partial_t\bv_\delta}_{L^2(0,T; \bV_\sigma')}+\norm{\partial_t\varphi_\delta}_{L^2(0,T;(H^1(\Omega))')}+\norm{\partial_tQ_\delta(\theta_\delta)}_{L^2(0,T;H^3(\Omega)')}\leq C(T,m),\label{dt12}
\end{align}
for any $T>0$, uniformly in $\delta>0$.

Using \eqref{E1}, \eqref{eq:estfur1},\eqref{eq:estfur2}, \eqref{eq:estfur4}, and \eqref{pressure1}, we also obtain
 \begin{align}
 			\| (\mathbb{S}(\theta_\delta) \ve_\delta ) \ve_\delta\|_{\LL^\frac54(\Omega \times (0,T))} &\leq C,\label{eq:convSv1}\\
	\| e_\delta \ve_\delta\|_{\LL^\frac{10}9(\Omega \times (0,T))} &\leq C, \label{eq:convev1}\\
  		\| p_\delta \ve_\delta\|_{\LL^\frac{10}9( \Omega \times (0,T))} &\leq C. \label{eq:convpv1}
  \end{align}
Observe that 
  \begin{align}
  	 L^\infty(0,T; L^{1+\alpha}(\Omega))	\cap L^\beta(0,T; L^{3\beta}(\Omega))  \hookrightarrow
  	L^q 
  ( \Omega \times (0,T)) 
  	\label{eq:embtheta},
  \end{align}
  where 
  $$q >  \beta  + 1 \iff \alpha > 1/2.$$
  Hence, we deduce
  \begin{align} 
  	 \| \theta_\delta \|_{L^{(1+\beta)^+}( \Omega \times (0,T))} \leq C, \label{eq:boundthetabetaA}\\
  \| \hat{\kappa}(\theta_\delta) \|_{L^{1^+}( \Omega \times (0,T))} \leq C, \label{eq:estkappahatA}
  \end{align}
  where the notations $(1+\beta)^+$ and $1^+$ indicate an exponent strictly greater than $(1+\beta)$  and $1$, respectively.
  
  From the estimates above and by comparison in \eqref{eq:weaktoten2b} it also follows that
  \begin{align}\label{finquimenouno}
  	\left\| \partial_t \left(\frac12 |\ve_\delta|^2 + e_\delta\right) \right\|_{L^1(0,T;{ ({W}_\sigma^{1,10}(\Omega))'})} \leq C. 
  \end{align}
  Furthermore, using \eqref{eq:estWprime}, \eqref{eq:boundtheta1}, \eqref{E1}, and \eqref{eq:estvarphiH2}, together with \eqref{eq:estfur1}, we can deduce that
   \begin{align}
  \left\| \nabla \left(\frac12 |\ve_\delta|^2 + e_\delta\right) \right\|_{\LL^{1^+}( \Omega \times (0,T))} \leq C. \label{finqui}
  \end{align}

\subsubsection{Convergence as $k\to\infty$} \label{subsec:convk}
For fixed $\delta>0$ and $m\in \N$, the estimates above suffice to pass to the limit as $k\to\infty$. In particular, we can infer the existence of a nonrelabeled solution $(\ve_\delta, p_\delta, \varphi_\delta, \theta_\delta)$ satisfying \eqref{eq:weakbalmomapprox}-\eqref{eq:weakheatapprox} (with the term $\frac1k\Delta\widetilde \mu_\delta$ going to zero in $L^2(0,T;(H^1(\Omega))')$ due to \eqref{restimate}). 
Additionally, all estimates \eqref{E1}-\eqref{finqui} are still valid in the limit $k\to\infty$ by lower semicontinuity of the norms involved. In particular, the energy inequality \eqref{energineqA} is preserved in the limit, namely 
\begin{align}
 &  \int_\Omega e_\delta(T)+\frac12\norma{\ve_\delta(T)}^2\dx\leq \int_\Omega e_\delta(0)+\frac12\norma{\ve_\delta(0)}^2\dx, \label{energineqAb}
\end{align}
for any $T>0$. On the other hand, we can also pass to the limit in the localized total energy inequality \eqref{eq:weaktoten2bf}, obtaining
\begin{align}   \label{eq:weaktoten2bfnew}
		&\int_{0}^T \int_{\Omega}  \Big(e_\delta + \frac12 |\ve_\delta|^2\Big) \partial_t \eta \, \dx \dt
		+	\int_{0}^T \int_{\Omega}  \Big(e_\delta \ve_\delta + \frac12 |\ve_\delta|^2\mathbf P_\sigma (J_m\ast\ve_\delta) \Big)\cdot \nabla \eta  \, \dx \dt 
        \\& \notag 
		+ \int_{0}^T \int_{\Omega} (\hat\kappa_\delta(\theta_\delta  )+\frac1m q_\delta) \Delta \eta \, \dx \dt	
		+ 	\int_{0}^T \int_{\Omega} p_\delta \ve_\delta  \cdot \nabla \eta \, \dx \dt
         \\& \notag 
		- \int_{0}^T \int_{\Omega} (\mathbb{S}(\theta_\delta)\ve_\delta ):  (\ve_\delta \otimes \nabla \eta ) \, \dx \dt
          +\int_{\Omega} \Big(e_0 + \frac12 |\ve_0|^2\Big) \eta(\cdot, 0)\, \dx \geq 0,\notag
	\end{align}
 for all $\eta \in C^\infty_c([0,T); C^2(\overline{\Omega}))$ such that $\eta \geq 0$ and
    $\nabla \eta  \cdot \no=0$ on $\partial \Omega \times [0,T)$.
Hence, we obtain that this inequality also holds unchanged in the limit as $k\to\infty$. Since the most relevant passage to the limit is that as $\delta\to0$, we omit further details for the limit $k\to\infty$, which proceeds in a similar way, and focus instead on $\delta\to0$.

\subsubsection{Convergence properties as $\delta\to 0$} Collecting the bounds above, which still hold as $r\to\infty$, up to the extraction of (not relabeled) subsequences, we have the following convergences  to some limit quadruple $(\ve_m, p_m, \varphi_m, \theta_m)$, globally defined on $(0,\infty)$, for any $m\in \N$,
 \begin{align}
 	\ve_\delta \rightarrow \ve_m \quad &\text{weakly star in } L^\infty(0,T; \LL^2(\Omega))	\cap L^2(0,T; {\mathbf H}^1(\Omega)), \label{eq:convv}\\
 	\varphi_\delta \rightarrow \varphi_m \quad &\text{weakly star in } L^\infty(0,T; H^1(\Omega))	\cap L^2(0,T; H^2(\Omega)), \label{eq:convu} \\
 	\theta_\delta \rightarrow \theta _m\quad &\text{weakly star in } L^\infty(0,T; L^{1+\alpha}(\Omega))	\cap  L^\beta(0,T; L^{3\beta}(\Omega)) \,, \label{eq:convtheta0} \\
 	\theta_\delta \rightarrow \theta_m \quad &\text{weakly in } L^2(0,T; H^1(\Omega))\label{eq:convtheta1},\\
    p_\delta\rightarrow p_m,\quad &\text{weakly in }L^{\frac53}(0,T;L^\frac53(\Omega)),\label{eq:press}
 \end{align}
 for any $T>0$, as $\delta \rightarrow 0$.  Using the Aubin-Lions lemma, we can infer
  \begin{align}
 \ve_\delta\rightarrow \ve_m \quad \text{strongly in } L^{2}(0,T; \mathbf L^2(\Omega)) , \label{eq:strongv1}
 \end{align}
   \begin{align}
  \varphi_\delta \rightarrow \varphi_m \quad \text{strongly in } C([0,T]; L^2(\Omega)) \cap L^{2}(0,T; H^{2-r }(\Omega)).\label{eq:strongconvu2}
  \end{align}
  for all $ r \in(0,2)$, as well as 
    \begin{align}
  Q_\delta(\theta_\delta) \rightarrow Q_m \quad \text{strongly in }  L^{2}(0,T; L^2(\Omega)),\label{eq:strongconvu3}
  \end{align}
  where we note that $Q_m\geq0$ since $Q_\delta(\theta_\delta)\geq0$ for any $\delta>0$. By continuity of the function $Q^{-1}$, since $Q_m\geq0$ and $Q_\delta^{-1}\to Q^{-1}$ uniformly on $[0,+\infty)$, this also entails, up to a further subsequence, that
  \begin{align}
  \theta_\delta=Q_{\delta}^{-1}(Q_\delta(\theta_\delta)+\delta)\to Q^{-1}(Q_m)\geq 0,
  \label{thetad}
\end{align}
  almost everywhere in $\Omega\times(0,T)$, for any $T>0$, entailing, by uniqueness of the limit, that $\theta_m=Q^{-1}(Q_m)\geq 0$.
  
 \subsubsection{Bound of $\vphi_\delta$ from below.}
Our aim is now to prove the lower bound for $\vphi$ stated in \eqref{eqn:boundsACl}, first at the $\delta$-level of approximation of the system. Then, this property is preserved in the limit by the strong convergence property \eqref{eq:strongconvu2} as $\delta \to 0$. 
We remark that this property can be established only after letting $k\to\infty$, as it relies on a maximum principle, which would not hold if the term $\frac1k\Delta\widetilde\mu_\delta$ was present in the advective Allen-Cahn equation. 

Testing the Allen-Cahn equation \eqref{eq:weakACapprox} with {$-(\vphi_\delta+1)^-$}, we have 
 \begin{align*}
 	\frac 12\frac {\d}{\dt} \norm{(\vphi_\delta  + 1)^-}^2_{L^2(\Omega)}+ \norm{\nabla(\vphi_\delta + 1)^-}^2_{L^2(\Omega)}-\frac1{\varepsilon^2}\int_\Omega W'(\vphi_\delta)(\vphi_\delta +1)^-\dx=-\frac{1}{\varepsilon}\int_\Omega \ell_\delta(\theta_\delta) (\vphi_\delta + 1)^-\dx\leq 0,
 \end{align*}
 where we have exploited the fact that $\opdiv \v_\delta=0$ together with the boundary conditions and the fact that $\ell_\delta(\theta_\delta)\geq 0$. 
Moreover, we also notice that 
$$
-\frac1{\varepsilon^2}\int_\Omega W'(\vphi_\delta)(\vphi_\delta + 1)^- \dx\geq 0,
$$ 
 since $W'(u)=u^3-u\leq 0$ for {$u+1\leq 0$}. 
  Assuming then that $\vert \vphi_0\vert\leq 1$, as it is natural and physical to require, we have 
 $$
 (\vphi_\delta(0) + 1)^-=(\vphi_0 + 1)^-\equiv0,
 $$
and thus Gronwall's Lemma entails 
\begin{align}
\label{b1}
	\vphi_\delta\geq -1\quad\text{almost everywhere in }\Omega\times (0,T).
\end{align}

 \subsubsection{Bound of $\vphi_\delta$ from above (in case of bounded latent heat).}
We prove the upper bound for $\vphi_\delta$ stated in \eqref{eqn:boundsACl}, first at the $\delta$-level of approximation of the system. As one can read below, its proof relies on the boundedness of the latent heat $\ell$, and it is the reason why this property of $\ell$ is taken as an assumption. Once again, this property is preserved in the limit by the strong convergence property \eqref{eq:strongconvu2} as $\delta \to 0$.

For any fixed $\tau>0$, we test the Allen-Cahn equation by $(\vphi_\delta-1-\tau)^+$, to obtain, after some integration by parts exploiting the boundary conditions,
\begin{align*}
	&\frac 12\frac {\d}{\dt}\norm{(\vphi_\delta-1-\tau)^+}^2_{L^2(\Omega)}+\norm{\nabla(\vphi_\delta-1-\tau)^+}^2_{L^2(\Omega)}+\frac1{\varepsilon^2}\int_\Omega W'(\vphi_\delta)(\vphi_\delta-1-\tau)^+\dx\\ &=\frac{1}{\varepsilon}\int_\Omega \ell_\delta(\theta_\delta) (\vphi_\delta-1-\tau)^+ \dx,
\end{align*}
whence, since $W'(s)$ is monotone for $s>1$, and since, crucially, $\ell_\delta(s)\leq L$, 
\begin{align*}
	\frac 12\frac {\d}{\dt}\norm{(\vphi_\delta-1-\tau)^+}^2_{L^2(\Omega)}+\norm{\nabla(\vphi_\delta-1-\tau)^+}^2_{L^2(\Omega)}+\frac1{\varepsilon}\int_\Omega \Big(\frac1{\varepsilon} W'(1+\tau) -L\Big)(\vphi_\delta-1-\tau)^+\dx  \leq 0,
\end{align*}
so that, if $0<\varepsilon<\frac{W'(1+\tau)}{L}$, we conclude that 

\begin{align*}
	\frac 12\frac {\d}{\dt}\norm{(\vphi_\delta-1-\tau)^+}^2_{L^2(\Omega)}+\norm{\nabla(\vphi_\delta-1-\tau)^+}^2_{L^2(\Omega)} \leq 0,
\end{align*}
and since $\vphi_\delta(0)\leq 1$, we always have, after integrating in time, $\norm{(\vphi_\delta-1-\tau)^+}_{L^2(\Omega)}=0$ for any $t\geq0$, giving 
\begin{align}
\label{b2}
	\vphi_\delta\leq 1+\tau\quad\text{almost everywhere in }\Omega\times (0,T).
\end{align}

 \subsubsection{{Strict positivity of $\theta_m$}.} \label{subsubsec:strictpostheta}
For $\alpha\in(\frac12,1]$, we now show that we can further prove that $\theta_m>0$ almost everywhere in $\Omega\times(0,T)$. So far we only have $\theta_\delta\in (Q_\delta^{-1}(\delta),+\infty)$, where $Q_\delta^{-1}(\delta)>0$ by construction, for $\delta$ sufficiently small, so that we infer $\theta_m\geq 0$.
In order to prove the strict positivity, we need to make a further assumption on the initial data, namely \eqref{assumption}. This is due to the lack of mass conservation in the convective Allen-Cahn equation (see also Remark \ref{essa}). We can start by recalling that, from \eqref{logdelta3},
\begin{align}
	\norm{\nabla \log_\delta(\theta_\delta)}_{L^2(0,T;\LL^2(\Omega))}\leq C.
	\label{nablalog}
\end{align} 
Also, note that, from the entropy inequality \eqref{eq:weakent2}, recalling \eqref{b1}-\eqref{b2} and $\vphi_0\geq-1$, we deduce, fixing $\tau_0>0$ so that $\frac{W'(1+\tau_0)}{L}>\eps$, so that \eqref{b2} holds with $\tau=\tau_0$, that 
\begin{align*}
	-\vert \Omega\vert +\int_\Omega \Lambda_\delta(\theta_\delta(0))\, \dx\leq \int_\Omega(\vphi+\Lambda_\delta(\theta_\delta))\, \dx\leq (1+\tau_0)\vert \Omega\vert +\int_\Omega \Lambda_\delta(\theta_\delta)\, \dx,\quad \forall \varepsilon\in(0,\eps_{\tau_0}). 
\end{align*}
This entails
\begin{align}
	\text{ess inf}_{t\in(0,T)}\int_\Omega \Lambda_\delta(\theta_\delta)\, \dx \geq 	-(2+\tau_0)\vert \Omega\vert +\int_\Omega \Lambda_\delta(\theta_\delta(0))\, \dx=c_{0,\delta}>0, \quad \forall \varepsilon\in(0,\eps_{\tau_0}), \label{below}
\end{align}
for $\delta$ sufficiently small, by exploiting assumption \eqref{assumption}, since $\Lambda_\delta\to \Lambda$ on $(0,\infty)$ and $\theta_\delta(0)=Q_\delta^{-1}(Q_\delta(\theta_0)+\delta)\to \theta_0$ as $\delta\to0$ almost everywhere in $\Omega$, thus $c_{0,\delta}\to c_0=-(2+\tau_0)\vert \Omega\vert +\int_\Omega \Lambda(\theta_0)\, \dx>0$ as $\delta\to0$. 
Also, note that, from \eqref{eqn:boundSQ01}-\eqref{eqn:boundSQ}, we infer that
\begin{align}
	\text{ess sup}_{t\in(0,T)}\int_\Omega \Lambda_\delta(\theta_ \delta(\cdot, t))\, \dx\leq C.
	\label{boundLambda}
\end{align}
Let us consider now a subsequence $\delta_l\to 0$ and $M_l\to\infty$ as $l    \to\infty$. Observe that, for almost any $t\in(0,T)$,
\begin{align*}
	\vert\{x\in\Omega: \theta_{\delta_l}(x,t)\geq M_l\} \vert\leq \frac{1}{\Lambda_{\delta_l}(M_l)}\int_\Omega \Lambda_{\delta_l}(\theta_{\delta_l}(\cdot, t))\, \, \dx\leq \frac{C}{\Lambda_{\delta_l}(M_l)}\to 0\text{ as }l\to \infty,
\end{align*}
since $\Lambda_{\delta_l}\to \Lambda$ as $l\to\infty$, and the convergence is uniform.
Additionally, from \eqref{eq:strongconvu3} and since $\Lambda$ is continuous, it is immediate to see that $\Lambda_{\delta_l}(\theta_{\delta_l})=\Lambda_{\delta_l}(Q_{\delta_l}^{-1}(Q_{\delta_l}(\theta_{\delta_l}))$ converges almost everywhere as $l\to\infty$, to $\Lambda(\theta_m)$. Therefore, 
\begin{align}
	\Lambda_{\delta_l}(\theta_{\delta_l}(\cdot, t))\chi_{\{\theta_{\delta_l}(\cdot, t)\geq M_l\}}\to 0,
\end{align}
almost everyhwhere in $\Omega$, as $l\to\infty$. Then, thanks to \eqref{boundLambda}, by dominated convergence we also deduce that, for almost any $t\in(0,T)$,
\begin{align*}
		\int_{\theta_{\delta_l}(\cdot, t)\geq M_l}\Lambda_{\delta_l}(\theta_{\delta_l}(\cdot, t))\, \dx \to 0,\quad \text{ as }l\to \infty.
\end{align*}
Now by the Severini-Egorov theorem, for a fixed $k\in\mathbb N$ there exists $A_k\subset [0,T]$ such that $\vert A_k\vert\leq \frac1{2^k}$ and 
\begin{align}
	\sup_{t\in[0,T]\setminus A_k}	\int_{\theta_{\delta_l}(\cdot, t)\geq M_l}\Lambda_{\delta_l}(\theta_{\delta_l}(\cdot, t))\, \dx \to 0,\quad \text{ as }l  \to \infty.
	\label{uniform1}
\end{align}

Our aim is now to show a control on the entire $\log \theta$ function. The first step is to prove that there exist three positive constants, $N>0$, $\underline{\theta}>0$ and $V_0>0$, such that 
\begin{align}
	\inf_{t\in[0,T]\setminus A_k}\vert\{x\in{\Omega}:\ \underline{\theta}\leq \theta_\delta(x,t)\leq N\} \vert\geq V_0,\label{poinc}
\end{align}
for any $\delta>0$ sufficently small.
Assume by contradiction that this does not hold. This means that there exist a sequence $\{\delta_l\}_l$, $\delta_l\to0$, and a sequence $\{t_l\}_l\subset [0,T]\setminus A_k$ such that 
	 \begin{align}
		\vert\{x\in{\Omega}:\ \frac1l\leq \theta_{\delta_l}(x,t_l)\leq l\} \vert=:\vert E_l\vert< \frac 1{l\Lambda_{\delta_l}(l)}.
		\label{to_zero}
	\end{align}
	This immediately entails that $\Lambda_{\delta_l}(\theta_{\delta_l}(\cdot, t_l))\to 0$ strongly in $L^1(\Omega)$. Indeed, we have
	\begin{align*}
		&\int_\Omega \Lambda_{\delta_l}(\theta_{\delta_l}(\cdot, t_l))\, \dx=	\int_{E_l} \Lambda_{\delta_l}(\theta_{\delta_l}(\cdot, t_l))\, \dx+	\int_{\theta_{\delta_l}(t_l)<\frac1l} \Lambda_{\delta_l}(\theta_{\delta_l}(\cdot, t_l))\, \dx+\int_{\theta_{\delta_l}(\cdot, t_l)> l} \Lambda_{\delta_l}(\theta_{\delta_l}(\cdot, t_l))\, \dx
		\\& \leq \Lambda_{\delta_l}(l)\vert E_l\vert+\Lambda_{\delta_l}(\tfrac1l)\vert \Omega\vert+\int_{\theta_{\delta_l}(\cdot, t_l)> l} \Lambda_{\delta_l}(\theta_{\delta_l}(\cdot, t_l)) \, \dx\\&
		\leq \frac1l+\Lambda_{\delta_l}(\tfrac1l)\vert \Omega\vert+\sup_{t\in [0,T]\setminus A_k}\int_{\theta_{\delta_l}(\cdot, t)> l} \Lambda_{\delta_l}(\theta_{\delta_l}(\cdot, t))\, \dx\to 0\quad\text{as }l\to\infty,
	\end{align*}
where we exploited the uniform convergence \eqref{uniform1}. Clearly, this is in contradiction with \eqref{below}, concluding the proof of \eqref{poinc}. Notice that the constants $N,\underline\theta,V_0$ might depend on $k$. 

From the generalized Poincaré's inequality (e.g., given in \cite[Theorem 10.14]{NovoFe}), it holds, for almost any $t\in[0,T]\setminus A_k$,
\begin{align*}	&\norm{\log_\delta(\theta_\delta)}_{H^1(\Omega)}\leq C(V_0)\left(\norm{\nabla \log_\delta(\theta_\delta)}_{\LL^2(\Omega)} +\int_{\underline\theta \leq \theta\leq N}\vert\log_\delta(\theta_\delta)\vert\, \dx\right)\\&
	\leq C(V_0)\left(\norm{\nabla \log_\delta(\theta_\delta)}_{\LL^2(\Omega)} +\vert \Omega\vert\max\{\log_\delta(\underline{\theta}), \log_\delta(N)\}\right),
\end{align*}
so that, integrating over $[0,T]\setminus A_k$, this gives from \eqref{nablalog}, that 
\begin{align}
\int_{[0,T]\setminus A_k}\norm{\log_\delta(\theta_\delta)}_{H^1(\Omega)}dt\leq C.\label{ff}
	\end{align}
Notice that the constant $C>0$ is uniform in $\delta$, since $\max\{\log_\delta(\underline{\theta}), \log_\delta(N)\}\to \max\{\log(\underline{\theta}), \log(N)\}$ as $\delta\to0$.
Let us now fix $\eta\in(0,1)$, then we have from \eqref{ff}, for $0<\delta<\eta$,
	\begin{align*}
		&
        \vert \{(x,t)\in \Omega\times [0,T]\setminus A_k:\ \theta_\delta(x,t)\leq \eta\}\vert\\&\leq \frac{\int_{\Omega} \vert \log_\delta(\theta_\delta) \vert \, \dx}{\vert \log \eta\vert}\leq \frac{C}{\vert \log \eta\vert},
	\end{align*}	
	entailing, by \eqref{thetad}, letting $\delta\to0$,
    	\begin{align*}
		&
        \vert \{(x,t)\in \Omega\times [0,T]\setminus A_k:\ \theta_m(x,t)\leq \eta\}\vert
        \leq \frac{C}{\vert \log \eta\vert},
	\end{align*}
Then, by Fatou's Lemma letting $\eta\to0^+$, 
	\begin{align}
			\vert \{(x,t)\in \Omega\times [0,T]\setminus A_k:\ \theta_m(x,t)\leq 0\}\vert=0.
		\label{finale1}
	\end{align}
	Therefore, since this can be obtained for any $k\in \mathbb N$, we have, for any $j\in \mathbb N$, 
		\begin{align}
		\vert \{(x,t)\in \Omega\times [0,T]\setminus \bigcup_{k\geq j}A_k:\ \theta_m(x,t)\leq 0\}\vert=0,
		\label{finale}
	\end{align}
	and since by construction
	$$
	\Bigg\vert \bigcup_{k\geq j} A_k\Bigg\vert\leq \sum_{k\geq j}\frac 1 {2^k}\to 0\quad\text{as }j\to \infty,
	$$
	we easily deduce also 
		\begin{align}
		\vert \{(x,t)\in \Omega\times (0,T):\ \theta_m(x,t)\leq 0\}\vert=0,
		\label{finale2}
	\end{align}
	entailing that $\theta_m>0$ almost everywhere in $\Omega\times(0,T)$, for any $T>0$.

 \subsubsection{Weak sequential stability}By the convergences \eqref{eq:convu}-\eqref{eq:strongconvu3}, we can deduce that the quadruple $(\ve_m, p_m,  \varphi_m, \theta_m)$  solves
 \begin{align} \label{eq:weakbalmomA}
		&\int_{0}^T \int_{\Omega} (\ve_m \cdot \partial_t \xe + (\ve_m  \otimes \P_\sigma ( J_m\ast\ve_m)  ): \nabla \xe  
		) \, \dx \dt \notag \\
		&=	\int_{0}^T \int_{\Omega} \mathbb{S}(\theta_m)\ve_m: \nabla \xe \, \dx \dt+\int_{0}^T \int_{\Omega} p_m\opdiv\xe \, \dx \dt\\&\nonumber\quad - 	\int_{0}^T \int_{\Omega} \eps (\nabla \varphi_m  \otimes \nabla \varphi_m ): \nabla \xe \, \dx \dt- \int_\Omega \ve_0 \cdot \xe(\cdot, 0) \, \dx ,
	\end{align}
    for all $\xe \in C^\infty_c([0,T); C^1(\overline{\Omega}; \mathbb{R}^3 ))$ such that $\xe \cdot \no=0$ on $\partial \Omega \times [0,T)$, where $\mathbf v_m$ satisfies the incompressibility condition $\opdiv \mathbf v_m =0$ a.e. in $(0,T)\times \Omega$, the Allen-Cahn system
	\begin{align}
	&	\partial_t \varphi_m +  \ve_m \cdot \nabla \varphi_m  = - \frac1{\varepsilon} \mu_m,
      \notag   \\
        & \mu_m = -\varepsilon \Delta \varphi_m  +  \frac{1}{\varepsilon} \partial_\varphi W(\varphi_m) - \ell(\theta_m),  
        \label{eq:weakACA}
	\end{align}
   a.e. in $\Omega \times (0,T)$, and the heat equation 
    \begin{align} 
  		&\int_{0}^T \int_\Omega Q(\theta_m) \partial_t \zeta \, \dx \dt
  		+ \int_{0}^T \int_\Omega Q(\theta_m)   \ve_m \cdot \nabla \zeta \, \dx \dt
  		\notag\\&\notag
        - \int_{0}^T \int_\Omega \kappa(\theta_m)\nabla \theta_m\cdot \nabla \zeta \, \dx \dt-\frac1m\int_0^T\int_\Omega \nabla Q(\theta_m)\cdot\nabla \zeta \, \dx\dt
  		 \\ 
  		&= - \int_{\Omega}Q(\theta_0) \zeta( \cdot, 0) \, \dx + \int_{0}^T \int_\Omega \ell(\theta_m)\frac{D\varphi_m}{Dt}\zeta \, \dx \dt\notag\\&\quad - \int_{0}^T \int_\Omega \bigg( ( {\nu(\theta_m)}|\nabla \ve_m + \nabla \ve_m^\mathsf{T}|^2) + 
  	\frac{1}{\varepsilon} \big((| \mu_m |^2 \big))\bigg) \zeta\, \dx \dt \label{eq:weakentA},
  	\end{align} 
  	for all $\zeta \in C^\infty_c( [0,T); C^2(\overline \Omega))$ such that $\nabla \zeta \cdot \no = {0}$ along ${\partial \Omega} \times [0,T)$, for any $ T >0$,
 and for any $m\in \N$. 
 
Note that the pressure $p_m$ satisfies 
 \begin{align}\label{eq:weakp2}
		\int_{\Omega} p_m \Delta \eta \, \dx = \int_{\Omega} (\nu(\theta_m)D\bv_m- \eps \nabla \varphi_m \otimes \nabla \varphi _m - \ve_m \otimes \P_\sigma(\ve_m)) : \nabla \nabla \eta \, \dx   , 
	\end{align}
	for any $\eta \in \mathbf H^2(\Omega)$ with $\nabla \eta \cdot \no = {0}$ along ${\partial \Omega}$, for almost any $t>0$. We also recall that $\theta_m>0$ almost everywhere in $\Omega\times(0,\infty)$.

Additionally, from \eqref{energineqAb} and \eqref{eq:weaktoten2bfnew}, we infer the total energy inequality
\begin{align}
 &  \int_\Omega e_m(T')+\frac12\norma{\ve_m(T')}^2\dx\leq \int_\Omega e(0) +\frac12\norma{\ve(0)}^2\dx,\label{energineqAbcd}
\end{align}
for almost every $T' \in (0,T)$, together with the localized total energy inequality
\begin{align}   \label{eq:weaktoten2bf1}
		&\int_{0}^T \int_{\Omega}  \Big(e_m + \frac12 |\ve_m|^2\Big) \partial_t \eta \, \dx \dt
		+	\int_{0}^T \int_{\Omega}  \Big(e_m \ve_m + \frac12 |\ve_m|^2\mathbf P_\sigma (J_m\ast\ve_m) \Big)\cdot \nabla \eta  \, \dx \dt 
        \\& \notag 
		+ \int_{0}^T \int_{\Omega} (\hat\kappa_m(\theta_m  )+\frac1m q_m) \Delta \eta \, \dx \dt	
		+ 	\int_{0}^T \int_{\Omega} p_m \ve_m  \cdot \nabla \eta \, \dx \dt
         \\& \notag 
		- \int_{0}^T \int_{\Omega} (\mathbb{S}(\theta_m)\ve_m ):  (\ve_m \otimes \nabla \eta ) \, \dx \dt
		 +\int_{\Omega} \Big(e_0 + \frac12 |\ve_0|^2\Big) \eta(\cdot, 0)\, \dx \geq 0,\notag
	\end{align}
 for all $\eta \in C^\infty_c([0,T); C^2(\overline{\Omega}))$ such that $\eta \geq 0$ and $\nabla \eta  \cdot \no=0$ on $\partial \Omega \times [0,T)$.
    
    Since $\theta_m>0$, we can divide \eqref{eq:weakentA} by $\ell(\theta_m)$, to obtain the entropy production equality:
  	\begin{align} \notag
  		&\int_{0}^T \int_\Omega(\Lambda(\theta_m) + \varphi_m) (\partial_t \zeta +  \ve_m \cdot \nabla \zeta )\, \dx \dt
  + \int_{0}^T \int_\Omega (h(\theta_\delta )+h_{1}(\theta_m)) \Delta \zeta \, \dx \dt
  		\\ 
        \label{eq:weakent2m}
  	&= - \int_{0}^T \int_\Omega \bigg( \frac{1}{\ell(\theta_m)} ({\nu(\theta_m)}|\nabla \ve_m + \nabla \ve_m^\mathsf{T}|^2) + \frac{\ell^\prime (\theta_m)}{\ell^2(\theta_m)}  {\Big(\kappa(\theta_m)+\frac1m Q'(\theta_m) \Big)}  |\nabla \theta_m|^2
  		 \bigg) \zeta\, \dx \dt   \\ \notag
        & \quad - \int_{0}^T \int_\Omega \bigg(  \frac{1}{\varepsilon\ell(\theta_m)} (|\mu_m |^2) \bigg) \zeta\, \dx \dt \notag   - \int_{\Omega}(\Lambda(\theta_0) +  \varphi_{\eps,0}) \zeta( \cdot, 0) \, \dx ,
  	\end{align}
  	where $\Lambda(s):=\int_0^s \frac{Q(\tau)}{\ell(\tau)}\d\tau$, and $h(\theta_m)= \int_{1}^{\theta_m} \frac{\kappa(s)}{\ell(s)} \, \d s $ 
  	for all { $\zeta \in C^\infty_c( [0,T); C^2(\overline \Omega))$ } such that {$\nabla \zeta \cdot \no =0$ along $\partial \Omega$}, for any $ T >0$, where we set $h_{1}(s):=\frac1m\int_1^s\frac{Q(\tau)}{\ell(\tau)}\d\tau$. 
    
\subsection{Limit as $m\to\infty$.} \label{subsec:limitm}
The aim of this section is to finally pass to the limit $m\to \infty$, obtaining a quadruple $(\bv,p, \vphi,\theta)$ solving \eqref{eq:weakbalmom}-\eqref{ineq:enconsl} and \eqref{eq:weakent}.

\subsubsection{Estimates and convergence properties} 
Thanks to \eqref{energineqAbcd} and to \eqref{eq:weakent2m}, by slightly adapting the arguments, one can deduce that estimates \eqref{E1}-\eqref{finqui} still hold uniformly in $m$. This allows to infer the convergences \eqref{eq:convv}-\eqref{eq:press} to a suitable globally defined quadruple $(\bv,p, \vphi,\theta)$.
Additionally, one can mimic step by step the proofs to obtain \eqref{b1}, \eqref{b2}, and \eqref{ff} (recalling that $\ell$ is bounded), as long as one replace the index $\delta$ with the index $m$. This leads to 
\begin{align*}
-1\leq \vphi_m\leq 1+\tau_0, \text{ almost everywhere in }\Omega\times(0,\infty),
\end{align*}
for some $\tau_0=\tau_0(\eps)$, as well as to 
\begin{align}
\int_0^T\norm{\log\theta_m}_{H^1(\Omega)}^2\dt\leq C.\label{logcontrol}
\end{align}
Recall now that
\begin{equation} \label{eq:defh}
	h(\theta_m) =   \int_{1}^{\theta_m} \frac{\kappa(s)}{\ell(s)} \, \mathrm{d} s = \begin{cases}- \int_{\theta_m}^{1} \frac{\kappa(s)}{\ell(s)} \, \mathrm{d} s  
		&\quad \text{ for } \theta_m < 1,
	\\
	\int_{1}^{\theta_m} \frac{\kappa(s)}{\ell(s)} \, \mathrm{d} s 
&	\quad \text{ for } \theta_m \geq  1
	.
	\end{cases}
\end{equation}
Hence, using the assumptions (d) and (e) on $\ell$ listed in Section \ref{subsec:hp}, we obtain
\begin{itemize}
	\item for $\theta < 1$, 
	\begin{equation} \label{eq:boundh0}
	|	h(\theta_m) | \leq \int_{\theta}^{1} \frac{\kappa(s)}{\ell(s)} \, \mathrm{d} s \leq 
			\int_{\theta}^{1}\frac{\kappa(s)}{\delta_0 \ell^\prime(0) s} \, \mathrm{d} s  = \frac{1}{\delta_0 \ell^\prime(0) } (-\kappa _1 \log \theta_m + \frac{\kappa_2}{\beta}(1- \theta_m^\beta))
	\end{equation}
	\item for $\theta \geq  1$, 
		\begin{equation} \label{eq:boundh1}
		|	h(\theta_m) | \leq 	\int_{1}^{\theta_m} \frac{\kappa(s)}{\ell(s)} \, \mathrm{d} s \leq 	\int_{1}^{\theta_m} \frac{\kappa(s)}{\ell(1)} \, \mathrm{d} s  =  \frac{1}{ \ell(1) } (\kappa _1 (\theta_m -1)+ \frac{\kappa_2}{\beta+1}(\theta_m^{\beta+1}-1)),
	\end{equation}
\end{itemize}
which implies
\begin{align}
\label{htheta1}
\norm{h(\theta_m)}_{L^{1^+}(\Omega\times(0,T))}\leq C(T),\quad \forall T>0,
\end{align}
due to \eqref{logcontrol}. Note that estimate \eqref{logcontrol}, which is here crucial to close the estimate on $h(\theta_m)$, heavily relies on the assumption that $\ell$ is a sublinear latent heat.

Next, we need to show similar estimates to \eqref{dt12}, which allow to obtain the strong convergences analogous to \eqref{eq:strongv1}, \eqref{eq:strongconvu2}, and \eqref{eq:strongconvu3}.
 \paragraph{\it Strong convergence of $\ve_m$.}
 By comparision in \eqref{eq:weakbalmomA}, we deduce
 \begin{align}
 	\| \partial_t\ve_m\|_{L^{\frac{5}{3}}(0,T; { (\mathbf{W}_\sigma ^{1,\frac52}(\Omega))'})} \leq C,
 \end{align}
{We can further prove that
 \begin{align}
 	\norm{\partial_t\ve_m}_{L^\frac43(0,T;\mathbf V_\sigma')} \leq C.\label{regv}
 \end{align}
 Indeed, for any $\xe\in \Vs$, by Gagliardo-Nirenberg's inequalities,
 \begin{align*}
 & \bigg|	\int_\Omega (\ve_m\cdot \nabla)\P_\sigma(\ve_m)\cdot \xe \, \dx\bigg|  = \bigg|-	\int_\Omega (\ve_m\otimes \P_\sigma(\ve_m)): \nabla\xi \, \dx \bigg|\\
 	&\leq C\norm{\ve_m}_{\LL^3(\Omega)}\norm{\ve_m}_{\LL^6(\Omega)}\norm{\nabla\xi}_{\LL^2(\Omega)}\leq C\norm{\ve_m}^\frac12_{\LL^2(\Omega)}\norm{\nabla\ve_m }^\frac32_{\LL^2(\Omega)}\norm{\xi}_{\Vs},
 \end{align*}
 so that, recalling the estimates above,
 \begin{align*}
 	\norm{(\ve_m\cdot \nabla)\P_\sigma(\ve_m)}_{L^\frac43(0,T;\Vs')}\leq  C\norm{\ve_m}_{L^\infty(0,T;\LL^2(\Omega))}\norm{\ve_m}_{L^2(0,T;\Vs)}\leq C.
 \end{align*}
 Also, we have, for any $\xi\in \Vs$,
 \begin{align*}
	\bigg|\int_\Omega \eps(\nabla \vphi_m\otimes  \nabla \vphi_m):\nabla \xe \, \dx\bigg|
	&\leq  C \eps\norm{\vphi_m}_{W^{1,3}(\Omega)} \norm{\vphi_m}_{H^2(\Omega)}\norm{\xe}_{\Vs} 
\\&	\leq C\eps\norm{\vphi_m}_{H^1(\Omega)}^\frac12\norm{\vphi_m}_{H^2(\Omega)}^\frac 32\norm{\xe}_{\Vs}, \end{align*}
entailing 
\begin{align*}
\norm{\eps\text{div}(\nabla \vphi_m\otimes \nabla \vphi_m)}_{L^\frac43(0,,T;\Vs')}	\leq  C\eps\norm{\vphi_m}_{L^\infty(0,T;H^1(\Omega))}\norm{\vphi_m}_{L^2(0,T;H^2(\Omega))}\leq C_\eps,
\end{align*}
thanks to \eqref{eq:estvarphiH2}. Therefore \eqref{regv} follows by comparison in \eqref{eq:weakbalmomA}.}
 
\noindent  From \eqref{evelox} (of course still valid as $r\to\infty$ and $\delta\to0$) and \eqref{regv}, using the Aubin-Lions lemma, we infer
  \begin{align}
 \ve_m \rightarrow \ve \quad \text{strongly in } L^{\frac{4}{3}}(0,T; {\mathbf H}^s(\Omega)) \quad \text{for } s\in(0,1) , \label{eq:strongv}
 \end{align}
  but also
   \begin{align}
  	\ve_m \rightarrow \ve \quad \text{strongly in } C([0,T]; \Vs'),\label{eq:strongvC}
  \end{align}
 as $m \to \infty$.
 
 \paragraph{\it Strong convergence of $\varphi_m$.}
  Using \eqref{eq:estfur2}, by comparison in \eqref{eq:weakACA}, we obtain 
  \begin{align}
  \|\partial_t\varphi_m\|_{L^\frac{5}{3}(0,T; L^{\frac{5}{3}}(\Omega)) }\leq C,
  \end{align}
 Using the Aubin-Lions lemma, from \eqref{eq:convu} we deduce
   \begin{align}
  \varphi_m \rightarrow \varphi \quad \text{strongly in } C([0,T]; L^2(\Omega)) \cap L^{2}(0,T; H^{2-r }(\Omega)).\label{eq:strongconvu}
  \end{align}
  for all $ r \in(0,2)$, as $m \to \infty$.

 \paragraph{\it Strong convergence of the total energy.} 
 Using \eqref{evelox}, \eqref{eq:estfur1},\eqref{eq:estfur2}, \eqref{eq:estfur4}, and \eqref{pressure1}, we obtain
 \begin{align}
 			\| (\mathbb{S}(\theta_m) \ve_m ) \ve_m\|_{\LL^\frac54(\Omega \times (0,T))} &\leq C,\label{eq:convSv}\\
	\| e_m \ve_m\|_{\LL^\frac{10}9(\Omega \times (0,T))} &\leq C, \label{eq:convev}\\
  		\| p_m \ve_m\|_{\LL^\frac{10}9( \Omega \times (0,T))} &\leq C. \label{eq:convpv}
  \end{align}
Observe that 
  \begin{align}
  	 L^\infty(0,T; L^{1+\alpha}(\Omega))	\cap L^\beta(0,T; L^{3\beta}(\Omega))  \hookrightarrow
  	L^q 
  ( \Omega \times (0,T)) 
  	\label{eq:embtheta1},
  \end{align}
  where 
  $$q >  \beta  + 1 \iff \alpha > 1/2.$$
  Hence, we deduce
  \begin{align} 
  	 \| \theta_m \|_{L^{(1+\beta)^+}( \Omega \times (0,T))} \leq C, \label{eq:boundthetabeta}\\
  \| \hat{\kappa}(\theta_m) \|_{L^{1^+}( \Omega \times (0,T))} \leq C, \label{eq:estkappahat}
  \end{align}
  where the notations $(1+\beta)^+$ and $1^+$ indicate an exponent strictly greater than $(1+\beta)$  and $1$, respectively.
  
  Recall that, by lower semicontinuity in \eqref{finquimenouno}, it holds
  \begin{align}
  	\left\| \partial_t \left(\frac12 |\ve_m|^2 + e_m\right) \right\|_{L^1(0,T;{ ({W}_\sigma^{1,10}(\Omega))'})} \leq C, 
  \end{align}
 as well as
   \begin{align}
  \left\| \nabla \left(\frac12 |\ve_m|^2 + e_m\right) \right\|_{\LL^{1^+}( \Omega \times (0,T))} \leq C. 
  \end{align}
  As a consequence, an application of the Aubin-Lions-Simon lemma yields
  \begin{align}\label{eq:stronge}
  	\frac12 |\ve_m|^2 + e_m  \rightarrow	\frac12 |\ve|^2 + e \quad \text{strongly in }  L^{1^+}( \Omega \times (0,T)),
  \end{align}
  as $m \to \infty$.
  Note that the limit $\frac12 |\ve|^2$ can be identified due to \eqref{eq:strongv}, whereas the first two terms in $e=\frac{\eps}2|\nabla \varphi|^2+ \frac1{\eps}W(\varphi)+ Q(\theta)$
due to \eqref{eq:strongconvu}. 
    For the identification of the limit $Q(\theta)$, we need the strong convergence of $\theta_m$ in $L^p( \Omega \times (0,T))$ for some $p\geq1$, whose proof follows below by exploiting the same monotonicity argument as in \cite{ERS3D}.

  \paragraph{\it Strong convergence of $\theta_m$.} 
  For a fixed $\eps>0$, from \eqref{eq:strongconvu} we can infer
   \begin{align}\label{eq:strongCH}
  	\frac{\eps}2|\nabla \varphi_m|^2+ \frac1{\eps}W(\varphi_m) \rightarrow		\frac{\eps}2|\nabla \varphi|^2+ \frac1{\eps}W(\varphi)  \quad \text{strongly in }  L^{1^+}( \Omega \times (0,T)),
  \end{align}
  as $m \to \infty$.
  Using the strong convergences \eqref{eq:stronge} and \eqref{eq:strongCH}, we obtain
  \begin{align*}
  &	\int_{0}^{T} \int_{\Omega} (Q(\theta_{m_1})- Q(\theta_{m_2}))\operatorname{sign}(\theta_{m_1} - \theta_{m_2}  ) \, \dx \dt \\
  &	= 	\int_{0}^{T} \int_{\Omega} (e_{m_1}- e_{m_2})\operatorname{sign}(\theta_{m_1} - \theta_{m_2}  ) \, \dx \dt\\
  & \quad 	- \int_{0}^{T} \int_{\Omega} \Big(\frac{\eps}2|\nabla \varphi_{m_1}|^2+ \frac1{\eps}W(\varphi_{m_1})- \frac{\eps}2|\nabla \varphi_{m_2}|^2- \frac1{\eps}W(\varphi_{m_2}) \Big)\operatorname{sign}(\theta_{m_1} - \theta_{m_2}  ) \, \dx \dt\\
  & \rightarrow 0 \quad \text{ as } {m_1},{m_2} \rightarrow \infty. 
  \end{align*}
  From the monotonicity of $Q$ it follows that $Q(\theta_m)$ is a Cauchy sequence in $L^{1}( \Omega \times (0,T))$, hence converges strongly and, up to a subsequence, almost
  everywhere to some limit $\tilde Q$. Moreover, having \eqref{eq:convtheta0}, using a generalized form of the Lebesgue theorem, we can deduce that $\tilde Q = Q(\theta)$. 
 Additionally, from \eqref{eq:boundthetabeta} we can infer 
  \begin{align} \label{eq:convthetap}
 \theta_m \rightarrow \theta \quad \text{strongly in }  L^{p}( \Omega \times (0,T)) \text{ for } p \in [1, (\beta+1)^+),
  \end{align}
  in particular $	\theta_m \rightarrow \theta $ almost everywhere in $\Omega \times (0,T)$, as $m \to \infty$. {Note that this, together with \eqref{eq:strongconvu}, allows to pass to the limit in \eqref{eq:weakent2m} for any $(\theta_m,\vphi_m)$, obtaining \eqref{extra}.}
  Observe that, since $\theta_m >0$ almost everywhere for any $m\in \N$, then it follows that $\theta\geq0$ almost everywhere in $\Omega\times(0,\infty)$.

\subsubsection{Weak sequential stability} 
   \paragraph{\it Limit of the balance of momentum and of the Allen-Cahn equation.}  
By letting $m \rightarrow \infty$, having the estimates and the convergences above, one can check that both the balance of momentum \eqref{eq:weakbalmomA} and the Allen-Cahn system \eqref{eq:weakACA} pass to the desired limits, up to the presence of the mollifier $J_m$. Indeed, we recover \eqref{eq:weakbalmom} by exploiting \eqref{eq:convv}, \eqref{eq:strongv}, 
and \eqref{eq:strongconvu}, whereas we obtain \eqref{eq:weakAC} by using  \eqref{eq:convu}, \eqref{eq:strongv}, \eqref{eq:strongconvu}, and \eqref{eq:convtheta1}, assuming $\ell$ to be bounded. 
Moreover, note that the initial datum in \eqref{eq:weakbalmom} can be obtained thanks to \eqref{eq:strongvC}.

   \paragraph{\it Limit of the total energy inequality.}  
   By letting $m\rightarrow \infty$ in \eqref{energineqAbcd}, using \eqref{eq:stronge}, we can also obtain the total energy inequality
   \begin{align}
 & \int_\Omega e(T')+\frac12\norma{\ve(T')}^2\dx \leq \int_\Omega e(0)+\frac12\norma{\ve(0)}^2\dx.\label{energineqA2}
\end{align}
for almost every $T' \in (0,T)$. Analogously, we can obtain the localized energy inequality by passing to the limit in \eqref{eq:weaktoten2bf1}, namely we get
\begin{align}   \label{eq:weaktoten2bf2}
		&\int_{0}^T \int_{\Omega}  \Big(e + \frac12 |\ve|^2\Big) \partial_t \eta \, \dx \dt
		+	\int_{0}^T \int_{\Omega}  \Big(e \ve + \frac12 |\ve|^2\ve \Big)\cdot \nabla \eta  \, \dx \dt 
        \\& \notag 
		+ \int_{0}^T \int_{\Omega} \hat\kappa(\theta  ) \Delta \eta \, \dx \dt	
		+ 	\int_{0}^T \int_{\Omega} p \ve  \cdot \nabla \eta \, \dx \dt
         \\& \notag 
		- \int_{0}^T \int_{\Omega} (\mathbb{S}(\theta)\ve ):  (\ve \otimes \nabla \eta ) \, \dx \dt
		+\int_{\Omega} \Big(e_0 + \frac12 |\ve_0|^2\Big) \eta(\cdot, 0)\, \dx \geq 0,\notag
	\end{align}
     for all $\eta \in C^\infty_c([0,T); C^2(\overline{\Omega}))$ such that $\eta \geq 0$ and $\nabla \eta  \cdot \no=0$ on $\partial \Omega \times [0,T)$,
   which corresponds to \eqref{eq:weaktoten}.
  \paragraph{\it Entropy production inequality.}
   At last, we need to prove the entropy production inequality \eqref{eq:weakent} by letting $m \rightarrow +\infty$. First, from \eqref{htheta1}, together with \eqref{eq:convthetap}, we can infer
   \begin{align}
   	h(\theta_m) \rightarrow h(\theta)   \quad \text{strongly in } {L^{1^+}( \Omega \times (0,T)) },\label{eq:strongconvh}
   \end{align}
   for any $T>0$.
   We recall that the thermal part of the entropy $\Lambda$ is given by \eqref{eq:defLambda}. Hence, from the uniform control $\sup_{t\in[0,T]}\int_\Omega\Lambda(\theta_m)\dx\leq C(T)$, coming from analogous computations to \eqref{eqn:boundSQ01}-\eqref{eqn:boundSQ}, using \eqref{eq:convthetap}, we deduce
  \begin{align} 
 	\Lambda(\theta_m) \rightarrow \Lambda(\theta) &\quad \text{strongly in }  L^{p}( \Omega \times (0,T))  \text{ for } p \in \Big[1, \frac{(\beta+1)^+}{\alpha+1}\Big),
 \end{align}
 where $\frac{(\beta+1)^+}{\alpha+1} > \frac32 $.
   Moreover, due \eqref{eq:strongv}, we have
   \begin{align}
   	& \Lambda(\theta_m) \ve_m  \rightarrow \Lambda(\theta)\ve  \quad \text{strongly in }  \LL^{1^+}((0,T) \times \Omega).
   \end{align}
   It follows that we can pass to the limit $m\rightarrow + \infty $ in the non quadratic terms of  \eqref{eq:weakent2m}.
   
   In order to deal with the quadratic terms in the right hand side of \eqref{eq:weakent2m}, first note that 
   \begin{align*}
   \frac1m \int_0^T\int_\Omega Q'(\theta_m))  |\nabla \theta_m|^2\eta \, \dx\dt\geq0,
   \end{align*}
hence we do not need to pass to the limit in this term.
For the other quadratic terms in the right hand side of \eqref{eq:weakent2m}, we need to exploit the lower semicontinuity result of Theorem \ref{ThIoffe}.
   In particular, we first recall \eqref{eq:convv} and \eqref{eq:convtheta1}, as well as the convergences $\nu(\theta_m) \rightarrow \nu(\theta) $, $\ell(\theta_m) \rightarrow \ell(\theta) $, $\ell^\prime(\theta_m) \rightarrow \ell(\theta) $, $\kappa(\theta_m) \rightarrow \kappa(\theta) $ almost everywhere in $\Omega \times (0,T)$, which follow from \eqref{eq:convthetap} and the continuity of the functions $\nu, \ell, \ell',\kappa$. Then, by comparison in \eqref{eq:weakACA}, using\eqref{eq:convv}, \eqref{eq:convu}, and \eqref{eq:convtheta0}, we can infer
   \begin{align}
   	\frac{D\varphi_m}{Dt}= \partial_t \varphi_m + \ve_m \cdot \nabla \varphi_m \rightarrow \frac{D\varphi}{Dt}= \partial_t \varphi +\ve \cdot \nabla \varphi \quad  \text{weakly in }  L^{2}((0,T) \times \Omega).
   \end{align}
   Note that the identification of the limit above follows from \eqref{eq:strongv} and \eqref{eq:strongconvu}. 
   At last, we can apply Theorem \ref{ThIoffe}, obtaining
   \begin{align*}
   	&\liminf_{m \rightarrow +\infty } \int_{0}^{T}  \int_{\Omega }	\frac{\ell^\prime (\theta_m)}{\ell^2(\theta_m)} \kappa(\theta_m) |\nabla \theta_{m}|^2\eta \, \dx \mathrm{d}t \geq 	\int_{0}^{T}  \int_{\Omega }	\frac{\ell^\prime (\theta)}{\ell^2(\theta)} \kappa(\theta) |\nabla \theta|^2 \eta \, \dx \mathrm{d}t,
   	\\&
   	\liminf_{m \rightarrow +\infty } 
   	 \int_{0}^{T}  \int_{\Omega } \left(\frac\eta{\ell(\theta_m)}\right){\nu(\theta_m)}|\nabla \ve_m + \nabla \ve_m^\mathsf{T}|^2 \, \dx \mathrm{d}t
     \geq  \int_{0}^{T}  \int_{\Omega } \frac{\nu(\theta )}{\ell(\theta)} |\nabla \ve + \nabla \ve^\mathsf{T}|^2\eta \, \dx \mathrm{d}t ,
   \\ & 
   	 	\liminf_{m\rightarrow +\infty } \int_{0}^{T}  \int_{\Omega }	 \left(\frac{\eta }{\varepsilon\ell(\theta_m)}\right) | \mu_m |^2 \, \dx \mathrm{d}t
   	 	\geq  \int_{0}^{T}  \int_{\Omega }	 \frac{1 }{\varepsilon\ell(\theta)} | \mu |^2\eta  \, \dx \mathrm{d}t,
   \end{align*}
   where $-\frac{1 }{\varepsilon}\mu = 	\partial_t \varphi +  \ve \cdot \nabla \varphi$ a.e. in $\Omega \times (0,T)$.
   Therefore, we recover the entropy production inequality  \eqref{eq:weakent} in the limit $m \rightarrow +\infty$. The proof of existence of the weak solutions is thus concluded.

   \subsection{Strict separation from zero of $\theta$}\label{Degiorgi}
    We know that the temperature $\theta$ is (almost everywhere) strictly positive.
    Here, we show that, if the initial datum $\theta_0$ is strictly separated from zero, one can deduce the strict separation property stated in Lemma \ref{lemma:DeGiorgi} (see \eqref{thetamin}). Clearly this proof is only formal, as it should be performed in the approximating scheme described in the previous section, so that the heat equation for $\theta$ is still valid. Then, by almost everywhere convergence, since the separation is uniform in the approximating parameters, one retrieves the same result for the limit function.
     To this aim, we adopt a De Giorgi iteration scheme, relying on the assumption $\alpha <1$ in \eqref{eq:defQ}.

 We consider the (heat) equation \eqref{eq:strongheat}.
   Here $\kappa(\theta)\geq \kappa_1$ and $Q(\theta):=\frac{\theta^{1+\alpha}}{1+\alpha}$, with $\alpha\in[0,1)$.
   We introduce a parameter $\Theta>0$ (to be determined below) and define the sequence $k_{m}$
   as $K_r=Q(\Theta)+\frac{Q(\Theta)}{2^m}$. Then it holds
   \begin{equation}
   	Q(\Theta) <k_{m+1}<k_{m}<2Q(\Theta) ,\qquad \forall m\geq 1,\qquad
   	k_{m}\searrow Q(\Theta) \qquad \text{as }m\rightarrow \infty.
   	\label{kn1bis}
   \end{equation}
   We now define $Q(\th)_{m}$ as $Q(\th)_m=(Q(\th)-K_r)^-$.
   For any $m\geq 0$, we also introduce the set: 
   \begin{equation*}
   	A_{m}(t):=\{x\in \Omega :Q(\th(x,t))-k_{m}\leq 0\},\quad \forall t\in [0,T].
   \end{equation*}
   It is evident that: 
   \begin{equation*}
   	A_{m+1}(t)\subseteq A_{m}(t),\qquad \forall m\geq 0,\qquad \forall t\in
   	[0,T].
   \end{equation*}
   
   We define the sequence $y_{m}$ as follows: 
   \begin{equation*}
   	y_{m}=\int_{0}^T\int_{A_{m}(s)}1\, \dx \d s,\qquad \forall m\geq 0.
   \end{equation*}
   For each $m\geq 0$, we consider the test function $v=-Q(\th)_{m}$
   and integrate it over the interval $I=[0,t]$, with $t\leq T$.
   This yields: 
   \begin{align}
   	& -\int_{I}\langle\partial _{t}Q(\th) ,Q(\th)_{m}\rangle_{H^1(\Omega)^{\prime
   		},H^1(\Omega)}\d s-\int_{I}\int_{A_{m}(s)}\v\cdot \nabla Q(\th)Q(\th)_m \, \dx \d s\notag\\&-\int_{I}\int_{A_{m}(s)}\ell(\th)\frac{D\vphi}{Dt}Q(\th)_m \, \dx \d s \notag 
   	-\int_{I}\int_{A_{m}(s)}\kappa(\th)\nabla\th\cdot\nabla Q_m(\th)_m \, \dx \d s
   	\\& = -\int_{I}\int_{A_{m}(s)}\nu(\th)\vert\nabla\v+\nabla\v^T\vert^2Q(\th)_n \, \dx \d s \label{phini}-\varepsilon\int_{I}\int_{A_{m}(s)}\left\vert\frac{D\vphi}{Dt}\right\vert^2Q(\th)_m \, \dx \d s
   \end{align}
   Now it is immediate to infer 
   $$
   -\int_{I}\langle\partial _{t}Q(\th) ,Q(\th)_{m}\rangle_{H^1(\Omega)^{\prime
   	},H^1(\Omega)} \, \d s=\frac12\norm{Q(\th(t))_m}^2-\frac12\norm{Q(\th(0))_m}^2,
   $$
   and, since $\theta_0\geq C$ by assumption, we can choose $\Theta>0$ sufficiently small so that also $Q(\th_0)>2Q(\Theta)$, and thus $\norm{Q(\th(0))_m}=0$ for any $m\in \N$.
   
   Moreover, since $\nabla Q(\th)\chi_{A_m}=-\nabla Q(\th)_m$, and since $\opdiv \v=0$, we have
   \begin{align*}
   	-\int_{I}\int_{A_{m}(s)}\v\cdot \nabla Q(\th)Q(\th)_m \dx \d s=\int_{I}\int_{\Omega}\frac12\v\cdot \nabla (Q(\th)_m)^2 \, \dx \d s=0.
   \end{align*}
   Then, exploiting the definition of $Q(\th)$, we observe 
   \begin{align*}
   	-\int_{I}\int_{A_{m}(s)}\kappa(\th)\nabla\th\cdot\nabla Q(\th)_m\, \dx \d s=\int_{I}\int_{A_{m}(s)}\kappa(\th)\theta^\alpha\vert\nabla \th\vert^2 \, \dx \d s\geq \kappa_1\int_{I}\int_{A_{m}(s)}\theta^\alpha\vert\nabla \th\vert^2 \, \dx \d s
   \end{align*} 
   Moreover, by Cauchy-Schwartz and Young's inequalities, we have
   \begin{align*}
   	\left\vert\int_{I}\int_{A_{m}(s)}\ell(\th)\frac{D\vphi}{Dt}Q(\th)_m \, \dx \d s\right\vert\leq \frac{1}{2{\rho} \varepsilon}\int_{I}\int_{A_{m}(s)} \ell(\th)^2Q(\th)_m \, \dx \d s +\frac{\rho \eps} 2\int_{I}\int_{A_{m}(s)}\left\vert\frac{D\vphi}{Dt}\right\vert^2Q(\th)_m \, \dx \d s
   \end{align*}
   for any $\rho >0$.
   Notice that, since $\th\geq 0$ almost everywhere in $\Omega\times(0,T)$, we deduce
   $$
   0\leq Q(\theta)_m\leq 2Q(\Theta),
   $$
   so that, recalling that, by definition of $Q(\cdot)$, on $A_m(t)$, $0\leq \theta\leq 2^{\frac{1}{1+\alpha}}\Theta$, and that $\ell$ is sublinear (i.e., $\ell(s)\leq \lambda s$),
   $$
   \frac{1}{2\rho}\int_{I}\int_{A_{m}(s)} \ell(\th)^2Q(\th)_m \, \dx \d s \leq \frac{\lambda^22^{\frac{2}{1+\alpha}}}{\rho} Q(\Theta)\Theta^2y_m.
   $$
   Moreover, considering the right-hand side in \eqref{phini}, we have 
   $$
   \frac{\rho \eps }2\int_{I}\int_{A_{m}(s)}\left\vert\frac{D\vphi}{Dt}\right\vert^2Q(\th)_m \, \dx \d s-{\rho \eps }\int_{I}\int_{A_{m}(s)}\left\vert\frac{D\vphi}{Dt}\right\vert^2Q(\th)_m \, \dx \d s\leq 0,
   $$
   so that the terms related to $\left\vert\frac{D\vphi}{Dt}\right\vert$ do not appear anymore in the estimates. In conclusion, we clearly also have that 
   $$
   -\int_{I}\int_{A_{n}(s)}\nu(\th)\vert\nabla\v+\nabla\v^T\vert^2Q(\th)_m \, \dx \d s  \leq 0
   $$
   To sum up, we have obtained  
   \begin{align}
   	\frac12\sup_{t\in[0,T]}\norm{Q(\th(\cdot, t))_m}^2+\kappa_1\int_0^T\int_{A_m(t)} \theta^\alpha\vert \nabla\th\vert^2 \, \dx \dt  \leq \frac{2^{\frac{2}{1+\alpha}}\lambda^2}{\rho} Q(\Theta)\Theta^2y_m:=X_m.
   	\label{rrt}
   \end{align}
   On the other hand, for any $t\in [0,T]$ and for almost any $x\in
   A_{m+1}(t) $, we have 
   \begin{equation*}
   	Q(\th(x,t))_{m}=\underbrace{Q(\Theta)+\frac{Q(\Theta)}{2^{m+1}}-Q(\th(x,t))}_{=Q(\th)_{m+1}\geq0}+Q(\Theta)\left(\frac1{2^m}-\frac{1}{2^{m+1}}\right)\geq \frac{Q(\Theta) }{2^{m+1}},
   \end{equation*}
   implying 
   \begin{equation*}
   	\int_{0}^T\int_{\Omega }|Q(\th)_{m}|^{3}\, \dx \d s\geq
   	\int_{0}^T\int_{A_{m+1}(s)}|Q(\th)_{m}|^{3}\, \dx \d s\geq \left( \frac{Q(\Theta) }{%
   		2^{m+1}}\right) ^{3}\int_{0}^T\int_{A_{m+1}(s)}1\, \dx \d s=\left( \frac{Q(\Theta) 
   	}{2^{m+1}}\right) ^{3}y_{m+1}.
   \end{equation*}
   Thus, by applying H\"{o}lder's inequality, we get
   \begin{align}
   	& \left( \frac{Q(\Theta) }{2^{m+1}}\right) ^{3}y_{m+1}\leq \left(
   	\int_{0}^T\int_{\Omega }|Q(\th)_{m}|^{\frac{10}3}\, \dx \d s\right) ^{\frac{9}{10}}\left(
   	\int_{0}^T\int_{A_{m}(s)}1 \, \dx \d s \right) ^{\frac{1}{10}}.
   	\label{est2pp}
   \end{align}
   Now by Gagliardo-Nirenberg's inequalities (see also \cite{GalPoia,P}) together with \eqref{rrt} and $0<Q(\th)_m\leq 2Q(\Theta)$,
   \begin{align*}
   	&	\int_{0}^T\int_{\Omega }|Q(\th)_{m}|^{\frac{10}3}\, \dx \d s\leq \sup_{t\in[0,T]}\norm{Q(\th(\cdot, t))_{m}}_{L^2(\Omega)}^\frac{4}3\int_0^T\left(\norm{\nabla Q(\th)_m}_{L^2(\Omega)}^2+\norm{Q(\th)_m}^2_{L^2(\Omega)}\right) \, \d s \\&
   	\leq 2^\frac23X_m^\frac23\left(\int_0^T\norm{\nabla Q(\th)_m}^2_{L^2(\Omega)}  \d s +4Q(\Theta )^2y_m\right),
   \end{align*}
   where we considered an equivalent norm on $H^1(\Omega)$ given by $\norm{\cdot}_{H^1(\Omega)}^2=\norm{\nabla\cdot}^2_{L^2(\Omega)}+\norm{\cdot}^2_{L^2(\Omega)}$.
   Recalling the definition of $Q_m(\th)$ and \eqref{rrt}, since, on $A_m(t)$, it holds $0<\th\leq 2^{\frac{1}{1+\alpha}}\Theta$, we thus have
   \begin{align*}
   	\int_0^T\int_\Omega\norma{\nabla Q(\th)_m}^2\, \dx \d s=\int_0^T\int_{A_m(s)}\th^{2\alpha}\norma{\nabla\th}^2\, \dx \d s\leq 2^{\frac{\alpha}{1+\alpha}}\Theta^\alpha \int_0^T\int_{A_m(s)}\th^{\alpha}\norma{\nabla\th}^2\, \dx \d s \leq \frac{2^{\frac{\alpha}{1+\alpha}}}{\kappa_1}\Theta^\alpha X_m,
   \end{align*}
   so that we have
   \begin{align*}
   	&	\int_{0}^T\int_{\Omega }|Q(\th)_{m}|^{\frac{10}3}\, \dx \d s\\&
   	\leq 2^\frac23X_m^\frac23\left(\frac{2^{\frac{\alpha}{1+\alpha}}}{\kappa_1}\Theta^\alpha X_m+4Q(\Theta)^2y_m\right)\\&
   	\leq C(\lambda,\rho)\Theta^\frac43Q(\Theta)^\frac23(\Theta^\alpha Q(\Theta)\Theta^2+Q(\Theta)^2)y_m^\frac{5}{3}\leq C(\lambda,\rho)\Theta^{4+\frac83\alpha}\left(\frac{\Theta}{1+\alpha}+\frac{1}{(1+\alpha)^2}\right)y_m^\frac53.
   \end{align*}
   In conclusion, setting
   $$
   K(\Theta):=\left(\frac{\Theta}{1+\alpha}+\frac{1}{(1+\alpha)^2}\right)^\frac{9}{10}{\Theta^{\frac{3}{5}(1-\alpha)}},
   $$
   we can write, from \eqref{est2pp}
   \begin{align*}
   	y_{m+1}\leq C(\lambda,\rho){2^{3m+3}}\left(\frac{\Theta}{1+\alpha}+\frac{1}{(1+\alpha)^2}\right)^\frac{9}{10}{\Theta^{\frac{3}{5}(1-\alpha)}}y_m^\frac85=C(\lambda,\rho)K(\Theta)2^{3m+3}y_m^\frac85,
   \end{align*}
   where, since $\alpha\in[0,1)$, $K(\Theta)\to 0$ as $\Theta\to 0$.
   
   Therefore, we can apply a well known geometric Lemma \cite[Lemma A.1]{GalPoia} (see also \cite[Lemma 3.8]{P} for a proof). In particular, in the notations of the lemma, we have $b=2^{3}>1$, $C=C(\lambda, \rho)K(\Theta)2^3>0$, $%
   \epsilon =\frac{3}{5}$, and we conclude that ${y}_{m}\rightarrow 0$, as
   long as 
   \begin{equation}
   	y_{0}\leq 2^{-\frac{40}{3}}C(\lambda,\rho)^{-\frac53}K(\Theta)^{-\frac{5}3},
   	\label{last5}
   \end{equation}
   but since
   \begin{align*}
   	y_{0} & = \int_{0}^T\int_{A_{0}(s)}1\, \dx \d s  \\ &= \int_{0}^T\int_{\{x\in \Omega : \ Q(\th (x,s))\leq 2Q(\Theta) \}}1\, \dx \d s\leq T\vert\Omega\vert,
   \end{align*}
   it is enough to ensure
   $$
   T\vert\Omega\vert\leq 2^{-\frac{40}{3}}C(\lambda,\rho)^{-\frac53}K(\Theta)^{-\frac{5}3},
   $$
   which can be satisfied as long as $\Theta$ is sufficiently small, since $K(\Theta)\to 0$ as $\Theta \to 0$.
   
   \noindent In summary, with this choice of $\Theta$ and by taking the limit as $m\rightarrow \infty$, we conclude that
   \begin{equation*}
   	\Vert (Q(\th) -Q(\Theta))^{-}\Vert _{L^{\infty }(\Omega \times (0,T))}=0.
   \end{equation*}
   This concludes the proof of \eqref{thetamin}, since it entails that there exists $c=c(T,\lambda,\rho)>0$ such that
   $$
   \theta\geq c \quad \text{almost everywhere in }\Omega\times(0,T).
   $$.

\section{Sharp interface limit of a non-isothermal Navier-Stokes/Allen-Cahn system} \label{sec:sharp}
This section concerns the proof of Theorem \ref{theo:convergenceBV}. The proof is divided into serveral steps, which are explained in Section \ref{subsec:proofmain}.
In some steps, we need to exploit additional intermediate results, whose statement are first collected in Section \ref{subsec:interesults}, then proved in Section \ref{subsec:proofinter} for the sake of readibility.

\subsection{Intermediate results}\label{subsec:interesults}
{Defining the localized interface energies 
	\begin{align}\label{eq:deflocEdiff}
		E_{\text{int},\eps}(\vphi_\eps; \zeta)
		&:= \int_{\Omega} \zeta
		\left(\frac{\eps}{2} |\nabla \vphi_\eps|^2+\frac{1}{\eps} W(\vphi_\eps)  
		\right) \, \dx
		,
		\\ \label{eq:deflocEsharp}
		E_{\text{int}}(\chi; \zeta)
		&:=  \sigma \int_{\partial^\ast A} \zeta \,  \mathrm{d}\mathcal{H}^{d-1} .
	\end{align}
	where $\zeta \in C(\overline\Omega)$, we have the following result, which is proven in Setion \ref{I1}.
	\begin{lemm} \label{lemma:convlocen}
		Let $\zeta \in C(\overline\Omega \times [0,T])$. Under the energy convergence assumption \eqref{eqn:hpenergyconv} , it holds
		\begin{align}\label{eq:convlocenint}
			\lim_{\eps \rightarrow 0} \int_{0}^T 
			E_{\operatorname{int},\eps}(\vphi_\eps; \zeta) \, \mathrm{d}t = \int_{0}^T 
			E_{\operatorname{int}}(\chi_A; \zeta) \, \mathrm{d}t . 
		\end{align} 
\end{lemm}
}
{We introduce the localized relative (interface) energy functionals
	\begin{align} \label{eq:locreleneps}
		\mathcal{E}_{\text{int},\eps}(\vphi_\eps; \zeta, \bm{\xi}):=& \, E_{\text{int},\eps}(\vphi_\eps; \zeta) - \int_{\Omega} \zeta (\bm{\xi} \cdot \nabla )\psi_\eps \, \dx  ,\\
		=& \,   E_{\text{int},\eps}(\vphi_\eps; \zeta) + \int_{\Omega} \psi_\eps (\zeta (\text{div} \bm{\xi}) + (\bm{\xi} \cdot \nabla )\zeta) \, \dx ,\notag 
		\\ \label{eq:locrelen}
		\mathcal{E}_{\text{int}} (\chi; \zeta, \bm{\xi} ):=&\, E_{\text{int}}(\chi; \zeta) 
		- \sigma \int_{\partial^\ast A} \zeta \bm{\xi} \cdot \n_A \frac{\nabla \chi}{|\nabla \chi|} 
		\, \mathrm{d}\mathcal{H}^{d-1} ,\\
		=&\, E_{\text{int}}(\chi; \zeta) 
		+ \sigma \int_{A} (\zeta (\text{div} \bm{\xi}) + (\bm{\xi} \cdot \nabla )\zeta) \, \dx ,\notag 
	\end{align}
	where $\zeta \in C^1(\overline \Omega ; [0,\infty))$ and $\bm{\xi} \in C^1(\overline \Omega ; \R^d)$ is such that $\bm{\xi}\cdot \no =0$ along $\partial \Omega$.
	The following result will be proven in Section \ref{I2}.
	\begin{lemm} \label{lemma:convlocrelen}
		Under the energy convergence assumption \eqref{eqn:hpenergyconv}, 
		then it holds 
		\begin{align}\label{eq:convlocrelen}
			\lim_{\eps \rightarrow 0} \int_{0}^T 
			\mathcal{E}_{\textnormal{int}, \eps}(\vphi_\eps ; \zeta, \bm{\xi} ) \, \mathrm{d}t= \int_{0}^T 
			\mathcal{E}_{\textnormal{int}}(\chi ; \zeta, \bm{\xi} ) \, \mathrm{d}t
		\end{align}
		for all 
		$\zeta \in C([0, T]; C^1(\overline\Omega; [0, \infty ))) $ and for all $\xe \in C([0, T]; C^1(\overline\Omega; \R^d)) $ 
	such that $\bm{\xi}\cdot \no=0$ along $\partial \Omega \times (0,T)$.
	\end{lemm}
Let ${\mathbf s} \in \mathbb{S}^{d-1}$ be a fixed but arbitrary unit vector. Defining $\bm{\nu}_\eps:= \frac{\nabla \vphi_\eps}{|\nabla \vphi_\eps|}$ if $\nabla \vphi_{\eps} \neq 0$, $\bm{\nu}_\eps:= {\mathbf s}$ if $\nabla \vphi_{\eps} = 0$, and $\n_A:= \frac{\nabla \chi_A}{|\nabla \chi_A|}$. Provided that  $$|\xe| \leq  1  \text{ on } \overline \Omega,$$ we can prove coercivity properties of the functionals \eqref{eq:locreleneps} and \eqref{eq:locrelen}, which give a tilt-excess type penalization of the difference between the normals, namely $\xe$ and $\bm{\nu}_\eps$, or $\xe$ and $\n_A$.
First, we notice that
\begin{align}
	\sigma  \int_{\partial^\ast A} \zeta \left|\n_A- \xe\right|^2 \, \mathrm{d}\mathcal{H}^{d-1} 	
	\leq \sigma  \int_{\partial^\ast A} \zeta \left(1- \xe \cdot \n_A\right) \, \d \mathcal{H}^{d-1}
	=	\mathcal{E}_{\text{int}}(\chi_A; \zeta, \xe ).
\end{align}
Second, we compute 
\begin{align}\label{eqn:coersharp}
	0 \leq 	\int_{\Omega} 
	\zeta \frac12	\left(\sqrt{\eps} |\nabla \vphi_\eps| - \frac{1}{\sqrt{\eps} } \sqrt{W(\vphi_\eps)}   \right)^2 \, \dx  
	+ \int_{\Omega} \zeta \Big(1- \xe\cdot \bm{\nu}_\eps\Big) |\nabla \psi_\eps  | \, \dx
	=	\mathcal{E}_{\text{int}, \eps} (\vphi_\eps; \zeta, \xe ) ,
\end{align}
which implies 
\begin{align} \label{eqn:coerequip}
	\int_{\Omega} 
	\zeta \frac12	\left(\sqrt{\eps} |\nabla \vphi_\eps| - \frac{1}{\sqrt{\eps} } \sqrt{W(\vphi_\eps)}   \right)^2 \, \dx  &\leq 	\mathcal{E}_{\text{int}, \eps} (\vphi_\eps; \zeta, 	\xe),\\ \label{eqn:coernormgrad}
	\int_{\Omega} \zeta \frac12 | \bm{\nu}_\eps - \xe|^2 |\nabla \psi_\eps|	\, \dx &\leq \mathcal{E}_{\text{int}, \eps}(\vphi_\eps; \zeta, \xe ) .
\end{align}
Furthermore, following the argument in \cite[Lemma 4]{FischerLauxSimon}, one can additionally prove that
\begin{align} \label{eqn:coernormpsi}
	\int_{\Omega} 	 \zeta \frac12 | \bm{\nu}_\eps - \xe|^2 \eps |\nabla \vphi_\eps|^2	\, \dx &\leq 	\mathcal{E}_{\text{int}, \eps}(\vphi_\eps; \zeta, \xe). 
\end{align}
}

\subsection{Proof of Theorem \ref{theo:convergenceBV}} \label{subsec:proofmain}
We recall that weak solutions $(\v_\eps , \varphi_{\eps}, \theta_\eps)$  to the non-isothermal Navier-Stokes/Allen-Cahn system satisfies further properties, which are collected in Lemma \ref{lemma:furtherproperties}. 

\subsubsection{Compactness}  	\label{subsec:compactnesssharp}
The first step consists in the extraction of an accumulation point for each sequence of $(\v_\eps , \varphi_{\eps}, \theta_\eps)_{\eps>0}$.

\textit{Proof of  \eqref{eqn:convergencepsi}.} We compute  $\partial_t \psi_\eps = \sqrt{2W(\varphi_\eps)} \partial_t \varphi_\eps $ and $\nabla \psi_\eps = \sqrt{2W(\varphi_\eps)} \nabla \varphi_\eps $. 
From the boundeness of the latent heat function $\ell$ (i.e., assumption (a) in Section \ref{subsec:hp}), it follows that 
\begin{align}
	\| (\partial_t + \ve_\eps \cdot \nabla) \psi_\eps  \|_{L^1(\Omega \times (0,T))} 
	& \leq \sqrt{2} \Big( \int_{0}^{T} \int_{\Omega} \frac{1}{\varepsilon} W(\varphi_{\varepsilon}) \, \dx \dt \Big)^\frac12 \Big( \int_{0}^{T} \int_{\Omega} {\varepsilon} |(\partial_t + \ve_\eps \cdot \nabla) \varphi_{\varepsilon}|^2 \, \dx \dt \Big)^\frac12
	\notag 	\\
	& \leq \sqrt{2L} \Big( \int_{0}^{T} \int_{\Omega} \frac{1}{\varepsilon} W(\varphi_{\varepsilon}) \, \dx \dt \Big)^\frac12 \Big( \int_{0}^{T} \int_{\Omega} \frac{\varepsilon}{\ell(\theta_{\eps})} {|(\partial_t + \ve_\eps  \cdot \nabla) \varphi_{\varepsilon}|^2 }\, \dx \dt \Big)^\frac12
	\notag 	\\
	&\leq C_T  \left(\sup_{\eps}  E_{\operatorname{int}, \eps }[\varphi_{\varepsilon,0}]+ { \sup_{\eps} S [\varphi_{\varepsilon} (\cdot, T), \theta_{\varepsilon}(\cdot, T)] } \right) \notag\\
	&\leq C_T  2 \sup_{\eps}  E_{\operatorname{tot}, \eps}[\varphi_{\varepsilon,0}, \ve_0, \theta_{0}]+ C(T, \Omega)< \infty,
\end{align}
where in the last inequality we used the bounds \eqref{eqn:boundSQ0}, \eqref{eqn:boundSQ}, as well as the uniform bound on $\frac{D\varphi}{Dt}$.

Similarly, we can show that
$\| \nabla  \psi_\eps  \|_{L^1(\Omega \times (0,T))}  \leq C_T  \sup_{\eps} E_{\operatorname{tot}, \eps}[\varphi_{\varepsilon,0}, \ve_0, \theta_{0}]  < \infty $ and
$$\textnormal{ess sup}_{t \in (0,T)} \| \nabla  \psi_\eps(\cdot,t)  \|_{L^1(\Omega )}  \leq C  \sup_{\eps} E_{\operatorname{tot}, \eps}[\varphi_{\varepsilon,0}, \ve_0, \theta_{0}]  < \infty ,$$
whence, together with the boundness of $\psi_\eps$, we can deduce 
\begin{align} \label{eq:esttotalvar}
	\textnormal{ess sup}_{t \in (0,T)} \|  \psi_\eps(\cdot,t)  \|_{W^{1,1}(\Omega )}  \leq C   .
\end{align}

From \eqref{ineq:enconsl} it follows that $\v_\eps $ is uniformly bounded in  $L^\infty( 0,T; \LL^2(\Omega))$, hence we have that $\v_\eps \psi_\eps \in L^\infty( 0,T; \LL^2(\Omega))$ due to \eqref{eqn:boundsACl}. 
This implies $\v_\eps \cdot   \nabla \psi_\eps  \in L^1(0,T;(H^{1}(\Omega))^\prime)$. Then, we observe that $\partial_ t \psi_\eps = \frac{D}{Dt} \psi_{\eps} - \v_\eps  (\nabla \cdot \psi_\eps)$, whence we deduce
that  $\partial_ t \psi_{\eps} \in L^1(0,T;(H^{1}(\Omega))^\prime)$. Since $W^{1,1}(\Omega ) \hookrightarrow  L^1(\Omega ) $ compactly and $ L^1(\Omega ) \hookrightarrow (H^{1}(\Omega))^\prime$, by an application of the Aubin-Lions-Simon Lemma, we obtain the strong convergence $\psi_{\eps} \to w$ in $L^1(  \Omega \times (0,T))$ as $\eps \to 0$, in particular (up to subsequence) $\psi_{\eps} \to w$ a.e. in $\Omega \times (0,T)$ and
$\psi_{\eps}(\cdot,t) \to w(\cdot,t)$ in $L^1(  \Omega)$  for almost all $t \in (0,T)$. Moreover, from the characterization \cite[Proposition 3.13]{AFP} of weak-$*$ convergence in $BV$, we can deduce $\psi_{\eps}(\cdot,t) \to w(\cdot,t)$ weakly-$*$ in $BV(\Omega)$ for almost all $t \in (0,T)$, hence by lower seminicontinuity of the total variation we can infer from \eqref{eq:esttotalvar} that $w \in L^\infty(0,T; BV(\Omega))$.

Finally, we need to identify the limit $w \in H^1(0,T; (H^1(\Omega))') \cap L^\infty(0,T; BV(\Omega))$. We first observe that the function $\psi$ as defined in \eqref{eqn:convergencepsi} is  continuous and invertible entailing that $\varphi_\eps \to \psi^{-1}(w)$ a.e. in $\Omega \times (0,T)$. Then, since for a.e. $t \in (0,T)$
$$\int_{\Omega} W(\psi^{-1}(w(\cdot, t))) \, \dx \leq \liminf_{\eps } \int_{\Omega} W(\phi_\eps(\cdot, t)) \, \dx   \leq \liminf_{\eps } \eps  E[\varphi_{\varepsilon,0}, \ve_0, \theta_{0}] = 0, $$
then we deduce $W(\psi^{-1}(w)) =0$ a.e. in $\Omega \times (0,T)$, implying $\psi^{-1}(w) \in \{-1,1\}$ a.e. in $\Omega \times (0,T)$. Setting $A(t)= \{x \in \Omega: \psi^{-1}(w(x,t)) =1\}$ for a.a. $t \in (0,T)$, we obtain $\psi^{-1}(w)= 2\chi_A -1$. From \eqref{eq:esttotalvar} it easily follows that $\chi_A(\cdot, t) \in BV(\Omega)$ for a.e. $t \in (0,T)$ and $w= \sigma \chi_A $.

\noindent \textit{Proof of  \eqref{eqn:convergencevel}.} We have
\begin{align*} 
	\| \ve_\eps \|_{L^\infty(0,T;\LL^2(\Omega))} < \infty,  \quad \| \ve_\eps\|_{{L^2(0,T;\Vs)} } < \infty, \quad 
		\| \partial_t \ve_\eps\|_{L^{1}(0,T; \mathbf H^{m}(\Omega )')} \leq C,
\end{align*}
for some $m \geq 3 $ by comparison in \eqref{eq:strongmom}, using the boundeness of $\sqrt{\eps}\nabla \varphi_\eps$ in $L^\infty (0,T;\LL^2(\Omega))$ due to \eqref{ineq:enconsl}. 
Since $\mathbf H^1(\Omega ) \hookrightarrow  \mathbf \LL^q(\Omega ) $ with $q < 6 $ compactly and $ \LL^q(\Omega ) \hookrightarrow \mathbf H^{m}(\Omega )'$, by an application of the Aubin-Lions-Simon Lemma, we obtain that
\begin{align}\label{eq:strongconvveps}
	\ve_\eps  \rightarrow \ve \quad \text{ strongly in } L^2([0,T];\LL^q(\Omega)), \quad q < 6,
\end{align}
and, in particular,
$$
 \lim_{\eps\to 0} \int_{0}^{T} \int_{\Omega} \frac12 |\ve_\eps|^2
\, \mathrm{d}x  \mathrm{d}t  =  \int_{0}^{T} \int_{\Omega} \frac12 |\ve|^2
\, \mathrm{d}x
 \mathrm{d}t. 
$$

\noindent \textit{Proof of  \eqref{eqn:convergenceQ}.} Having \eqref{ineq:enconsl} and \eqref{eq:entropyldissip}, proceeding as in Section \ref{subsec:estl}, we can infer
\begin{align} \label{eq:boundsQeps}
\| Q(\theta_\eps)\|_{L^\infty(0,T; L^1(\Omega))} < \infty \quad \text{ and } \quad \| \nabla \theta_\eps\|_{L^2(\Omega \times (0,T))}  < \infty.
\end{align}
recalling that $\alpha \in (1/2,1]$ in the definition of $Q$ given by \eqref{eq:defQ}.
Hence, from the assumption \eqref{eqn:hpQ} we can deduce
$$
\lim_{\eps\to 0} \int_{0}^{T} \int_{\Omega} Q(\theta_\eps)
\, \mathrm{d}x\mathrm{d}t  =  \int_{0}^{T} \int_{\Omega} Q(\theta)
\, \mathrm{d}x \mathrm{d}t ,
$$
hence the strong convergence \eqref{eqn:convergenceQ}.

\subsubsection{Total energy inequality}	
The total energy inequality \eqref{eqn:energyconservation} for almost every time $T' \in (0,T)$ follows from \eqref{ineq:enconsl}, \eqref{eqn:convE0}, \eqref{eqn:convergencevel}, \eqref{eqn:convergenceQ}, and \eqref{eqn:hpenergyconv}.
In particular, we have
\begin{align}\label{eq:oss1}
&\lim_{\eps\to 0} \int_{0}^{T} E_{ \operatorname{tot}, \eps} [\vphi_\eps(\cdot, t), \ve_\eps(\cdot, t) ,\theta_{\eps}(\cdot, t)]\, \mathrm{d}t\\ \notag 
& =  
\int_{0}^{T} E_{\operatorname{int}}[\chi_A(\cdot, t)]\, \mathrm{d}t  + \lim_{\eps\to 0} \int_{0}^{T} \int_{\Omega} \frac12 |\ve_\eps|^2
\, \mathrm{d}x  \mathrm{d}t  +  	\lim_{\eps\to 0} \int_{0}^{T} \int_{\Omega}Q(\theta_\eps)
\,\mathrm{d}x \mathrm{d}t \\ \notag 
&  =  	\int_{0}^{T} E_{\operatorname{int}}[\chi_A(\cdot, t)]\, \mathrm{d}t + \int_{0}^{T} \int_{\Omega} \frac12 |\ve|^2
\, \mathrm{d}x  \mathrm{d}t  +  	 \int_{0}^{T} \int_{\Omega}Q(\theta) \, \mathrm{d}x  \mathrm{d}t =  \int_{0}^{T} E [\chi_A(\cdot, t), \ve (\cdot, t), \theta (\cdot, t)]\, \mathrm{d}t
\end{align}
and 
\begin{align}\label{eq:oss2}
&\lim_{\eps\to 0} \int_{0}^{T} E_{ \operatorname{tot}, \eps} [\vphi_\eps(\cdot, t),\ve_\eps(\cdot, t) , \theta_{\eps}(\cdot, t)]\, \mathrm{d}t 
\\
&\leq  
\lim_{\eps\to 0} {T} E_{ \operatorname{tot}, \eps} [\vphi_{\eps,0}, \ve_0, \theta_{0}] =  T E_{\operatorname{int}}[\chi_{A(0)}]+ T \int_{\Omega} \frac12 |\ve_0|^2
\, \mathrm{d}x 
+ T \int_{\Omega}Q(\theta_0  ) \, \mathrm{d}x= T E [\chi_{A(0)}, \ve_0, \theta _0] .  \notag 
\end{align}

\subsubsection{Existence of a normal velocity}
We now show that the measure $\sigma (\partial_t \chi_A+ \ve \cdot \nabla \chi_A)$ is absolutely
continuous with respect to the product measure $\mathcal{L}^1 \llcorner(0,T) \otimes (\sigma \mathcal{H}^{d-1} \llcorner (\partial^\ast A(t)) _{t \in (0,T)})$.  Then, we denote by $V:\Omega \times (0,T)\rightarrow \R$ the associated Radon–Nikodym derivative, which therefore satisfies \eqref{eqn:transport}. Then,  we have the $L^2$-estimate
\begin{align}\label{eqn:boundV}
\int_0^T \int_{\partial^*A(t)} \sigma | V|^2 \, \mathrm{d}\mathcal{H}^{d-1}\mathrm{d}t& \leq \liminf_{\eps \to 0}  	\int_0^T \int_{\Omega}\eps |(\partial_{t}+ \ve_\eps \cdot 	\nabla) \varphi_\eps |^2 \, \dx \mathrm{d}t\\
&\leq L  \liminf_{\eps \to 0}  	\int_0^T \int_{\Omega}\eps \frac{|(\partial_{t}+ \ve_\eps \cdot 	\nabla) \varphi_\eps |^2}{\ell(\theta_\eps)} \, \dx \mathrm{d}t  < \infty, \notag 
\end{align}
which can be proved by adapting the proof of {\cite[Lemma 3.8]{HenselLaux}} as it follows.
In particular,  let $U \subset \Omega $ and $I  \subset  (0, T)$ be two open sets, and $ \zeta \in C^\infty_c (U \times I)$
such that $|\zeta| \leq 1$ and $\text{supp} \, \zeta \subset  U \times I $ . We then estimate 
\begin{align}
\nonumber&\sigma \langle \partial_t \chi_A + \ve \cdot \nabla \chi_A, \zeta \rangle  \\
&= \lim_{\eps \rightarrow 0} \Big( - \int_0^T \int_{\Omega} \psi_{\eps} (\partial_{t}+ \ve_\eps \cdot 	\nabla) \zeta \, \dx \dt \Big)=\lim_{\eps \rightarrow 0} \Big( \int_0^T \int_{\Omega} (\partial_{t}+ \ve_\eps \cdot 	\nabla)\psi_{\eps}  \zeta \, \dx \dt \Big) \notag \\
& \nonumber\leq \Big( \liminf_{\eps \to 0} \int_0^T \int_{\Omega}   \varepsilon |(\partial_{t}+ \ve_\eps \cdot 	\nabla) \vphi_{\eps}|^2 \, \dx \dt \Big)^\frac12
\Big( \limsup_{\eps \to 0} \int_0^T \int_{\Omega} |\zeta|^2  \frac1\varepsilon 2W(\vphi_{\eps}) \, \dx \dt \Big)^\frac12\\
&\nonumber \leq \sqrt{2L}  \Big( \liminf_{\eps \to 0} S[\vphi_\eps(\cdot,T), \theta_\eps(\cdot,T)]  \Big)^\frac12
\Big( \limsup_{\eps \to 0} \int_0^T \int_{\Omega} |\zeta|^2  \frac1\varepsilon W(\vphi_{\eps}) \, \dx \dt \Big)^\frac12 \\
&\nonumber \leq \sqrt{2L} C \Big( 1+ \liminf_{\eps \to 0} E_{ \operatorname{tot}, \eps} [\vphi_\eps(\cdot,T),\ve_\eps(\cdot,T), \theta_\eps(\cdot,T)]  \Big)^\frac12
\Big( \limsup_{\eps \to 0} \int_0^T \int_{\Omega} |\zeta|^2  \frac1\varepsilon W(\vphi_{\eps}) \, \dx \dt \Big)^\frac12\\
& \leq \sqrt{2L}  C\Big( 1 + E_{ \operatorname{tot}, \eps} [\vphi_{\eps,0}, \ve_0 , \theta_{0}]  \Big)^\frac12
\Big( \sigma \int_0^T \int_{\partial^\ast A} |\zeta|^2   \, \d \mathcal{H}^{d-1} \dt  \Big)^\frac12,\label{radon}
\end{align}
where we used the compactness \eqref{eqn:convergencepsi}, the entropy production \eqref{eq:entropyldissip} together with the boundeness of $\ell$ (see assumption (a) in Section \ref{subsec:hp}), then the bounds \eqref{eqn:boundSQ0},\eqref{eqn:boundSQ} and \eqref{eqn:boundsACl}, the total energy inequality \eqref{ineq:enconsl}, and at last the convergence of the local interface energy \eqref{eq:convlocenint} (cf. Lemma \ref{lemma:convlocen}). 
Proceeding as in the proof of {\cite[Lemma 3.8]{HenselLaux}}, one can improve the above estimate (suboptimal by a factor of $\sqrt{2}$) to obtain the $L^2$-estimate \eqref{eqn:boundV} for the associated Radon–Nikodym derivative $V$.

\subsubsection{Derivation of the motion law}
First, we test the non-isothermal Allen-Cahn equation \eqref{eq:strongAC} with $\eps (\Psi \cdot \nabla ) \varphi_\eps$, where $\Psi \in C^1_c(\overline\Omega \times (0,T); \R^d)$ with ${\Psi} \cdot  \no= 0$ on $\partial \Omega \times (0, T) $. Note that this test function is admissible due to the regularity stated in Definition \ref{def:weaksol}. 
We then obtain
\begin{align*}
& \int_0^T \int_{\Omega} (\eps (\Psi \cdot \nabla) \varphi_\eps ) (\partial_t \varphi_\eps + \ve_\eps \cdot \nabla\varphi_\eps)  \, \dx \mathrm{d}t
\\ \notag
&= - 	\int_0^T \int_{\Omega} \nabla (\eps (\Psi \cdot \nabla) \varphi_\eps ) \cdot \nabla \varphi_\eps 
+ (\eps (\Psi \cdot \nabla) \varphi_\eps ) \frac{1}{\eps^2} W'(\varphi_\eps) - (\eps (\Psi \cdot \nabla) \varphi_\eps ) \frac{1}{\eps} \ell (\theta_\eps ) 
\, \dx \mathrm{d}t.
\end{align*}
Since $W'(\varphi_\eps)  (\Psi \cdot \nabla) \varphi_\eps = (\Psi \cdot \nabla) W(\varphi_\eps)$, then integrating by parts we get
\[
-	\int_0^T \int_{\Omega} 
(\eps (\Psi \cdot \nabla) \vphi_\eps ) \frac{1}{\eps^2} W'(\vphi_\eps) \, \dx \mathrm{d}t
= 
\int_0^T \int_{\Omega} 
\frac{1}{\eps} W(\varphi_\eps)   (\text{div} \Psi) \, \dx \mathrm{d}t
\]
and using the symmetry $\nabla^2 \vphi_\eps$
\begin{align*}
&- 	\int_0^T \int_{\Omega} \nabla (\eps (\Psi \cdot \nabla) \vphi_\eps ) \cdot \nabla \vphi_\eps \, \dx \mathrm{d}t \\&=
- 	\int_0^T \int_{\Omega} \eps \nabla  \Psi  : \nabla \vphi_\eps  \otimes \nabla \vphi_\eps \, \dx \mathrm{d}t
+ \int_0^T \int_{\Omega} \frac{\eps}{2} |\nabla \vphi_\eps|^2 (\text{div} \Psi) \, \dx \mathrm{d}t .
\end{align*}
As a consequence, we obtain
\begin{align}\label{eq:motioneps}
& \int_0^T \int_{\Omega} (\eps (\Psi\cdot \nabla) \vphi_\eps ) \Big(\partial_t \vphi_\eps -  \frac{1}{\eps}\ell(\theta_{\eps})\Big) \, \dx \mathrm{d}t
\\ \notag
&= 
\int_0^T \int_{\Omega} 
\left(\left(\frac{\eps}{2} |\nabla \vphi_\eps|^2+\frac{1}{\eps} W(\vphi_\eps)  \right) \text{Id} - \eps\nabla \vphi_\eps  \otimes \nabla \vphi_\eps  \right) : \nabla  \Psi \, \dx \mathrm{d}t.
\end{align}

\paragraph{\textit{Step 1 - The curvature term.}}
In order to take the limit $\varepsilon \rightarrow 0$ in the first term on the right hand side of \eqref{eq:motioneps}, we consider the localized interface energies {defined in \eqref{eq:deflocEdiff} and \eqref{eq:deflocEsharp}.}
We take $\zeta = \text{div} \Psi$ as test function in \eqref{eq:convlocenint} and we obtain 
\begin{align}
\lim_{\eps \rightarrow 0} 	\int_0^T \int_{\Omega} 
\left(\frac{\eps}{2} |\nabla \vphi_\eps|^2+\frac{1}{\eps} W(\vphi_\eps)  \right)  \text{Id} : \nabla  \Psi \, \dx \mathrm{d}t 
= \sigma 	\int_0^T  \int_{\partial^\ast A}  \text{Id} : \nabla \Psi \,  \mathrm{d}\mathcal{H}^{d-1}  \mathrm{d}t. 
\end{align}
In order to take the limit $\varepsilon \rightarrow 0$ in the second term on the right hand side of \eqref{eq:motioneps}, we exploit 
{the ingredients given by Section \ref{subsec:interesults}}.

{More precisely, we can prove the following result.}
\begin{prop} \label{prop:convcurvature}
Under the energy convergence assumption \eqref{eqn:hpenergyconv} and the convergences \eqref{eqn:convergencepsi} and \eqref{eqn:convergenceQ}, it holds
\begin{align}\label{eq:convcurvature}
	\lim_{\eps \rightarrow 0} 	\int_0^T \int_{\Omega} 
	\eps\nabla \vphi_\eps  \otimes \nabla \vphi_\eps   : \nabla  \Psi \, \dx  \mathrm{d}t 
	= \sigma 	\int_0^T  \int_{\partial^\ast A} \n_A \otimes \n_A  : \nabla \Psi \,  \mathrm{d}\mathcal{H}^{d-1} \mathrm{d}t
\end{align}
for all $\Psi \in C^1_c(\overline\Omega \times (0,T); \R^d)$ with $\Psi \cdot  \no= 0$ on $\partial \Omega \times (0, T) $.
\end{prop}

{
A localization argument (by means of a suitable partition of unity) together with a local optimization of the choice of the vector fields $\xe$ (cf. \cite[Proof of Theorem 1.2]{LauxSimon}) 
shows that the convergence \eqref{eq:convcurvature} is implied by the following claim
\begin{align*}
	&\lim_{\eps \rightarrow 0} 	\Big| \int_0^T \int_{\Omega} 
	\eps\nabla \vphi_\eps  \otimes \nabla \vphi_\eps   : \nabla  \Psi \, \dx  \mathrm{d}t 
	- \sigma 	\int_0^T  \int_{\partial^\ast A} \nu_A \otimes \nu_A  : \nabla \Psi \,  \mathrm{d}\mathcal{H}^{d-1} \mathrm{d}t\Big|\\
	& \leq C_T \frac1\rho\int_0^T \mathcal{E}_{\text{int}, \eps}(\chi; \zeta, \xe)\, \dt  + \rho E[\chi_0, \theta_0]
\end{align*}
for all $\zeta \in C([0, T]; C^1(\overline\Omega; [0, \infty ))) $, for all $\Psi \in C^1_c(\overline\Omega \times (0,T); \R^d)$ with $\Psi \cdot  \no= 0$ on $\partial \Omega \times (0, T) $, for all $\xe \in C([0, T]; C^1(\overline\Omega; \R^d)) $ with $|\xe| \leq  1$ such that $\xe \cdot \no =0$ along $\partial \Omega \times (0,T)$, and any $\rho>0$.
The proof of the claim reads exactly as in \cite[Proof of Proposition 3.7]{HenselLaux} (neglecting the additional boundary contribution to the energy there considered). 
}

\paragraph{\textit{Step 2 - The velocity term.}}We now pass to the limit in the first term on left hand side of \eqref{eq:motioneps}.
We claim thet there exists a subsequence such that
\begin{equation}\label{eqn:claimconvvel}
\lim_{\eps\to 0} \left( \int_{0}^{T} \int_{\Omega}\eps ( \nabla \vphi_\eps \cdot \Psi) 
(\partial_t + \ve_\eps  \cdot \nabla )\vphi_\eps  
\, \dx \mathrm{d}t \right) = \int_{0}^{T} \int_{\partial^\ast A(t)} \Psi \cdot \n_A \sigma 
V
\, \mathcal{H}^{d-1} \dt 
\end{equation}
for all $\Psi \in C^1_c(\overline\Omega \times (0,T); \R^d)$ with $\Psi \cdot  \no = 0$ on $\partial \Omega \times (0, T) $.

For all any $\eta \in C^1([0, T]; C(\overline\Omega; [0, \infty ))) $ and $\xe \in C^1([0, T]; C(\overline\Omega; \R^d)) $ with $|\xe| \leq  1$ such that $\xe \cdot \no =0$ along $\partial \Omega \times (0,T)$
, adding zeros twice, we can write
\begin{align*}
&\int_{0}^{T} \int_{\Omega} \eta (\eps ( \nabla \vphi_\eps \cdot \Psi) (\partial_t + \ve_\eps  \cdot \nabla )\vphi_\eps )\, \dx \dt \\
&=  \int_{0}^{T} \int_{\Omega} \eta \Psi \cdot (\bm{\nu}_\eps - \xe ) \sqrt{\eps} |\nabla \vphi_\eps|  \sqrt{\eps} (\partial_t + \ve_\eps  \cdot \nabla ) \vphi_\eps \, \dx \dt \\
&\quad + \int_{0}^{T} \int_{\Omega} \eta (\Psi \cdot \xe) \Big(\sqrt{\eps} |\nabla \vphi_\eps|  - \frac{1}{\sqrt{\eps}} \sqrt{2W(\vphi_\eps)}\Big) \sqrt{\eps}  (\partial_t + \ve_\eps  \cdot \nabla ) \vphi_\eps \, \dx \dt  \\
&\quad + \int_{0}^{T} \int_{\Omega} \eta   (\Psi \cdot \xe) (\partial_t + \ve_\eps  \cdot \nabla )\psi_\eps\, \dx \dt .
\end{align*}
The first two terms on the right hand side can be bounded by means of the Young inequality combined together with \eqref{eqn:coerequip}, \eqref{eqn:coernormgrad}, and the entropy production \eqref{eqn:boundV}. 
Integrating by parts, using the evolution equation \eqref{eqn:transport}, and taking limits based on \eqref{eqn:convergencepsi}, we obtain
\begin{align*}
& \lim_{\eps\to 0}	\int_{0}^{T} \int_{\Omega} \eta  (\Psi \cdot \xe) (\partial_t + \ve_\eps  \cdot \nabla )\psi_\eps \, \dx \dt \\
&= 	- \lim_{\eps\to 0} \bigg(\int_{0}^{T} \int_{\Omega} \psi_\eps   (\partial_t + \ve_\eps  \cdot \nabla ) (\eta  \Psi \cdot \xe ) \, \dx \dt +\int_{\Omega}   \eta (\Psi \cdot \xe)(\cdot, 0) \psi_{\eps,0} \, \dx \dt \bigg) \\
&=   - \int_{0}^{T} \int_{\Omega} \sigma \chi_A    (\partial_t + \ve  \cdot \nabla )( \eta \Psi \cdot \xe ) \, \dx \dt - \int_{\Omega}  \eta  (\Psi \cdot \xe)(\cdot, 0) \sigma \chi_A(\cdot, 0) \, \dx \dt \\
&=  \int_{0}^{T} \int_{\partial^\ast A(t)} \eta  (\Psi \cdot \xe ) \sigma V  \, \d\mathcal{H}^{d-1} \dt. 
\end{align*}
As a consequence, using the convergences \eqref{eqn:convE0} and \eqref{eq:convlocrelen}, as well as the bound \eqref{eqn:VL2}, for a small constant $\delta>0$ we have
\begin{align*}
&\lim_{\eps\to 0} \Big| 	\int_{0}^{T} \int_{\Omega} \eta (\eps ( \nabla \vphi_\eps \cdot \Psi) (\partial_t + \ve_\eps  \cdot \nabla ) \vphi_\eps )\, \dx \dt  - \int_{0}^{T} \int_{\partial^\ast A(t)} \eta  (\Psi \cdot \n_A ) \sigma V  \, \d\mathcal{H}^{d-1} \dt \Big|\\
&\leq \lim_{\eps\to 0} \Big| 	\int_{0}^{T} \int_{\Omega} \eta (\eps ( \nabla \vphi_\eps \cdot \Psi) (\partial_t + \ve_\eps  \cdot \nabla ) \vphi_\eps )\, \dx \dt  - \int_{0}^{T} \int_{\partial^\ast A(t)} \eta  (\Psi \cdot \xi ) \sigma V  \, \d\mathcal{H}^{d-1} \dt \Big| \\
& \quad + \lim_{\eps\to 0}  \Big| 	\int_{0}^{T} \int_{\partial^\ast A(t)} \eta  (\Psi \cdot (\xe - \n_A) ) \sigma V  \, \d\mathcal{H}^{d-1} \dt \Big|\\
& \leq C_T \frac1\delta \int_0^T \mathcal{E}_{\text{int}, \eps}(\chi_A; \eta, \xe )\, \dt  + \delta (1+E[\chi_{A(0)}, \theta_0])
\end{align*}
We can conclude the proof of the claim \eqref{eqn:claimconvvel} by exploiting a localization argument (by means of a suitable partition of unity) together with a local optimization of the choice of the vector fields $\xe$ (cf. \cite[Proof of
Theorem 1.2]{LauxSimon}).

\paragraph{\textit{Step 3 - The temperature term.}}
At last, we consider the limit of the second term on the left hand side of \eqref{eq:motioneps}.
We compute 
\begin{align*}
- 	\int_{0}^{T} \int_{\R^d} (\eps ( \nabla \vphi_\eps \cdot \Psi) \frac{1}{\eps}\ell(\theta_{\eps} )\, \dx \dt 
&=   	\int_{0}^{T} \int_{\Omega}  \vphi_\eps (\text{div} \Psi) \ell(\theta_{\eps} )\, \dx \dt + 	\int_{0}^{T} \int_{\Omega}  \vphi_\eps \ell^\prime(\theta_{\eps} )\Psi \cdot \nabla\theta_{\eps} \, \dx \dt .
\end{align*}
Note that {$\varphi_{\varepsilon} \rightarrow 2\chi_A-1$} almost everywhere in $\Omega \times (0,T)$. Hence, having the bound \eqref{eqn:boundsACl}, we can deduce that $
	\varphi_{\varepsilon} \rightarrow \chi_A$ strongly in $L^p (\Omega \times (0,T))$ for any $p\in [1, \infty)$. 
Since $\ell$ is bounded and sublinear, the convergence \eqref{eqn:convergenceQ} yields
\begin{align*}
& \lim_{\eps\to 0}	\int_{0}^{T} \int_{\Omega}  \vphi_\eps (\text{div} \Psi) \ell(\theta_{\eps} )\, \dx \dt + 	\int_{0}^{T} \int_{\Omega}  \vphi_\eps \ell^\prime(\theta_{\eps} )\Psi \cdot \nabla\theta_{\eps} \, \dx \dt \\
&=   \int_{0}^{T} \int_{\Omega}  \chi_A (\text{div} \Psi) \ell(\theta )\, \dx \dt + 	\int_{0}^{T} \int_{\Omega}   \chi_A \ell^\prime(\theta )\Psi \cdot \nabla\theta \, \dx \dt \\
&= \int_{0}^{T} \int_{\Omega}  \chi_A  \text{div} (\Psi \ell(\theta ))\, \dx \dt, 
\end{align*}
which concludes the proof of the motion law \eqref{eqn:motionlaw}.

\subsubsection{Balance of momentum}  	
We aim to pass to the limit $\eps \rightarrow 0 $ in the balance of momentum \eqref{eq:weakbalmom}, i.e.,
\begin{align}
	&\int_{0}^T \int_{\Omega} (\ve_\eps \cdot \partial_t \xe + (\ve_\eps \otimes \ve_\eps ): \nabla \xe  + p_\eps \opdiv \xe ) \, \dx \dt \notag \\
	&=	\int_{0}^T \int_{\Omega} \mathbb{S}_\eps : \nabla \xe \,  \dx \dt  -	\int_{0}^T \int_{\Omega} \eps (\nabla \varphi_{\varepsilon} \otimes \nabla \varphi_{\varepsilon} ): \nabla \xe  \, \dx \dt - \int_\Omega \ve_0 \cdot \xe(\cdot, 0) \, \dx ,
\end{align}
where $\xe \in C^1_c([0,T); C^\infty(\overline{\Omega}; \mathbb{R}^d ))$ such that $\operatorname{div}\xe  =0$ and $\xe \cdot \no =0$ on $\partial \Omega \times (0, T) $.
Note that the last term in the first line is zero due to the choice of test functions.
From \eqref{eq:entropyldissip} and the boundenesss of $\ell$, we can deduce that 
$\ve_\eps \rightarrow \ve $ weakly in $L^2(0,T;H^1(\Omega))$. This fact combined with the strong convergence \eqref{eqn:convergencevel} and an application of Proposition \ref{prop:convcurvature} yield the balance of momentum \eqref{eqn:momentumsharp} in the limit  $\eps \rightarrow 0$.

\subsubsection{Entropy production inequality}	 	
We consider the sharp interface limit $\varepsilon \rightarrow 0$ of the entropy production inequality \eqref{eq:weakent}, namely
	\begin{align} \label{eq:entprodineqeps}
	&	
	-\int_{0}^{T}  \int_{\Omega }	(\Lambda (\theta_\eps) + \varphi_\varepsilon)(t) (\partial_t \zeta + \ve_\eps \cdot \nabla \zeta)\,  \dx \mathrm{d}t
	-	\int_{0}^{T}  \int_{\Omega }	  h(\theta_\eps) \Delta\zeta \,  \dx \mathrm{d}t
	\notag \\
	&	\geq \int_{\Omega }	(\Lambda (\theta_{0}) + \varphi_\varepsilon(\cdot, 0))\zeta(\cdot, 0) \, \dx +
	\int_{0}^{T}  \int_{\Omega }\frac{\ell^\prime (\theta_\eps)}{\ell^2(\theta_\eps)} \kappa(\theta_\eps) |\nabla \theta_\eps|^2 \zeta \, \dx \mathrm{d}t\nonumber \\
	& \quad + \int_{0}^{T}  \int_{\Omega }	\frac{\nu(\theta_\varepsilon)}{\ell(\theta_\varepsilon)} |\nabla \ve_\varepsilon + \nabla \ve_\varepsilon^\mathsf{T}|^2
	\zeta \, \dx \mathrm{d}t 
	+ \int_{0}^{T}  \int_{\Omega }\frac{\eps }{\ell(\theta_\eps)} |(\partial_t + \ve \cdot \nabla)\vphi_\eps|^2 \zeta \, \dx \mathrm{d}t,
\end{align}
where $h(\theta_\eps):= \int_{1}^{\theta_\eps} \frac{\kappa(s)}{\ell(s)} \, \d s $, for 
{$\zeta \in C^\infty_c( [0,T); C^2(\overline \Omega))$ }
such that $\zeta \geq 0 $ and {$\nabla \zeta \cdot \no =0$ along $\partial \Omega \times [0,T)$}, so that the first term on the right hand side is zero.

We recall that  $\varphi_{\varepsilon} \rightarrow 2\chi_A-1$ strongly in $L^p (\Omega \times (0,T))$ for any $p\in [1, \infty)$, the bounds \eqref{eqn:boundSQ0}-\eqref{eqn:boundSQ}, and the strong convergences \eqref{eqn:convergenceQ} and \eqref{eq:strongconvveps}. Moreover, from \eqref{eq:boundsQeps} we can deduce that $Q(\theta_\eps) 
$ is uniformly bounded in 
\begin{align*}
	L^\infty(0,T; L^1(\Omega)) \cap L^\frac{2}{1+\alpha}(0,T; L^\frac{6}{1+\alpha}(\Omega)) \hookrightarrow L^{4^-}(0,T; L^{{\frac{6}{5}}^+}(\Omega)) .
\end{align*}
As a consequence, we can conclude that the first term on the left side of \eqref{eq:entprodineqeps} passes to the desired limit in \eqref{eqn:entropydissipBV}. {Furthermore, we observe that the bound \eqref{eq:boundthetabeta} is uniform in $\eps>0$.
Therefore, recalling \eqref{eq:boundh0} and \eqref{eq:boundh1}, {it holds $h(\theta_\eps)\to h(\theta)$ strongly in $L^{1+}(\Omega\times(0,T))$ (see also \eqref{eq:strongconvh}), so that} the second term in \eqref{eq:entprodineqeps} easily passes to the limit.}

In order to pass to the limit $\varepsilon \rightarrow 0$ in the first two quadratic terms, we exploit the lower semicontinuity result stated in Theorem \ref{ThIoffe}
Recalling the strong convergence \eqref{eqn:convergenceQ}, we have that $\nu(\theta_\eps) \rightarrow \nu(\theta) $, $\ell(\theta_\eps) \rightarrow \ell(\theta) $, $\ell^\prime(\theta_\eps) \rightarrow \ell(\theta) $, $\kappa(\theta_\eps) \rightarrow \kappa(\theta) $ almost everywhere in $\Omega \times (0,T)$, due to the continuity of the functions $\nu, \ell, \ell',\kappa$. Moreover, we recall the properties of the bounded and sublinear $\ell$ (as listed in \ref{subsec:hp}) and the entropy production equality \eqref{eq:entropyldissip}, whence we can deduce the weak convergence in $L^1(\Omega \times (0,T))$ of the quadratic terms (thanks to the bounds derived in Section \ref{subsec:estl}).
Then, we can apply Theorem \ref{ThIoffe}, obtaining 
\begin{align}
	\label{basic1}\liminf_{\eps \rightarrow 0} \int_{0}^{T}  \int_{\Omega }	\frac{\ell^\prime (\theta_\eps)}{\ell^2(\theta_\eps)} \kappa(\theta_\eps) |\nabla \theta_{\eps}|^2 \zeta \, \dx \mathrm{d}t &\geq 	\int_{0}^{T}  \int_{\Omega }	\frac{\ell^\prime (\theta)}{\ell^2(\theta)} \kappa(\theta) |\nabla \theta|^2 \zeta\, \dx \mathrm{d}t,
	\\
	\liminf_{\eps \rightarrow 0 } 
	\int_{0}^{T}  \int_{\Omega } \frac{\nu(\theta_\eps )}{\ell(\theta_\eps)} |\nabla \ve_\eps + \nabla \ve_\eps^\mathsf{T}|^2 \zeta \, \dx \mathrm{d}t &\geq  \int_{0}^{T}  \int_{\Omega } \frac{\nu(\theta )}{\ell(\theta)} |\nabla \ve  + \nabla \ve^\mathsf{T}|^2 \zeta \, \dx \mathrm{d}t .
\end{align}
At last, we recall \eqref{eqn:boundV}, namely
\begin{align*}
\int_0^T \int_{\partial^*A(t)} \sigma | V|^2 \, \mathrm{d}\mathcal{H}^{d-1}\mathrm{d}t &\leq \liminf_{\eps \to 0}  	\int_0^T \int_{\Omega}\eps |(\partial_t + \ve_\eps   \cdot \nabla) \varphi_\eps |^2 \, \dx \mathrm{d}t \\
&\leq L  \liminf_{\eps \to 0}  	\int_0^T \int_{\Omega}\eps \frac{|(\partial_t + \ve_\eps   \cdot \nabla)\varphi_\eps |^2}{\ell(\theta_\eps)} \, \dx \mathrm{d}t  < \infty.
\end{align*}
We observe that $\frac{\varepsilon }{\ell(\theta_\varepsilon)} |(\partial_t + \ve_\eps   \cdot \nabla) \varphi_\varepsilon|^2  \dx \mathrm{d}t \rightarrow \mu_\ell  \text{ in }  \mathcal{M}_+(\Omega \times (0,T))$ as well as ${\varepsilon }|(\partial_t + \ve_\eps   \cdot \nabla)\varphi_\varepsilon|^2  \dx \mathrm{d}t \rightarrow \mu:= \mu_\ell \ell(\theta)  \text{ in }  \mathcal{M}_+(\Omega \times (0,T))$, due to the boundeness and sublinearity of $\ell$ together with the convergence \eqref{eqn:convergenceQ}. 
Hence, we have (up to a subsequence)
\begin{align*}
L  \langle \mu , 1 \rangle  \geq 	\langle \ell(\theta )\mu_\ell, 1 \rangle = \lim_{\eps \to 0}  	\int_0^T \int_{\Omega}\eps {|(\partial_t + \ve \cdot \nabla) \varphi_\eps |^2} \, \dx \mathrm{d}t \geq 
\int_0^T \int_{\partial^*A(t)} \sigma | V|^2 \, \mathrm{d}\mathcal{H}^{d-1}\mathrm{d}t 
\end{align*}	
Therefore, we obtain the entropy production inequality  \eqref{eq:weakent} in the sharp interface limit $\eps\rightarrow 0$.

\subsubsection{{Strict positivity of $\theta$}}
We need to show that, as $\eps\to 0$, we still have $\theta>0$ almost everywhere in $\Omega\times(0,T)$. This part is very similar to the proof of \eqref{finale2}. We only highlight how to pass the limit as $\eps \to 0$. First, we consider \eqref{extra}, holding for any $\eps>0$, and pass to the limit, recalling that $\theta_\eps\to\theta$ and $\vphi_\eps\to 2\chi_A-1$ almost everywhere in $\Omega\times(0,T)$, so that we have
\begin{align} \notag 
	-\vert \Omega\vert +\int_\Omega \Lambda(\theta_0)\, \dx &\leq \int_\Omega (\Lambda(\theta_0))\, \dx +\int_\Omega (2\chi_{A(0)}-1)\, \dx\\
	&\leq \int_\Omega (\Lambda(\theta)+2\chi_A-1)\, \dx\leq \vert \Omega\vert +\int_\Omega\Lambda(\theta)\, \dx.\label{below1}
\end{align} 
By assumption \eqref{assumption} (note that here also $\tau_0=0$ is already sufficient), the inequality \eqref{below1} entails 
$$\int_\Omega \Lambda(\theta)\, \dx\geq c_0>0,\quad \forall t\in(0,T).$$
{Note that also in this case we have the bound \eqref{nablalog}, as it follows from \eqref{eq:entropyldissip} together with an application of Theorem \ref{ThIoffe} (see \eqref{basic1}),
 \eqref{eqn:boundSQ0}-\eqref{eqn:boundSQ} and \eqref{ineq:enconsl}. Moreover, we have \eqref{boundLambda} due to \eqref{eqn:energyconservation} and once again \eqref{eqn:boundSQ0}-\eqref{eqn:boundSQ}.}
We recall that, since $\theta_\eps>0$ for any $\eps>0$, and $\theta_\eps\to \theta$ almost everywhere in $\Omega\times(0,T)$, it also holds $\theta\geq0$. Therefore, we can repeat the proof of \eqref{finale2} (without the dependence on $\delta$)  to infer that $\theta>0$ almost everywhere in $\Omega\times(0,T)$.

\subsection{Proof of intermediate results} \label{subsec:proofinter}
In this subsection we collect the proof of the intermediate results, {which we stated in Section \ref{subsec:interesults} and then exploited in the previous subsection to prove Theorem \ref{theo:convergenceBV}.}

\subsubsection{Proof of  Lemma \ref{lemma:convlocen}}\label{I1}
By linearity, it is enough to prove the statement for $\zeta \in [0, 1]$.
Having the convergence of the total energy \eqref{eqn:hpenergyconv} (corresponding to $\zeta = 1$), upon replacing $\zeta$ by $1-\zeta$ it is sufficient to prove the lower
bound. To this aim, we use the Young inequality and compute 
\begin{align*}
	&	 \int_{0}^T \int_{\Omega} \zeta 
	\left(\frac{\eps}{2} |\nabla \vphi_\eps|^2+\frac{1}{\eps} W(\vphi_\eps) 
	\right) \, \dx \dt \\
	& \geq 	 \int_{0}^T \int_{\Omega} \zeta 
	|\nabla \vphi_\eps| \sqrt{W(\vphi_\eps)  } 
	\, \dx =  \int_{0}^T \int_{\Omega} \zeta 
	|\nabla \psi_\eps|   
	 \, \dx \dt . 
\end{align*}
By (localized) lower semi-continuity statement, we have
\begin{align*}
	&	\lim_{\eps \to 0} \int_{0}^T \int_{\Omega} \zeta |\nabla \psi_\eps|  
	\, \dx	\geq \sigma \int_{\partial^\ast A} \zeta \,  \mathrm{d}\mathcal{H}^{d-1}, 
\end{align*}
whence it follows that
\begin{align*}
	&	\liminf_{\eps \to 0} \int_{0}^T \int_{\Omega} \zeta 
	\left(\frac{\eps}{2} |\nabla \vphi_\eps|^2+\frac{1}{\eps} W(\vphi_\eps)  
	\right) \, \dx	\geq \sigma \int_{\partial^\ast A} \zeta \,  \mathrm{d}\mathcal{H}^{d-1}.  
\end{align*}

\subsubsection{Proof of  Lemma \ref{lemma:convlocrelen}}\label{I2}
The convergence \eqref{eq:convlocrelen} of the localized relative entropies is a direct
consequence of the convergence \eqref{eq:convlocenint} of the localized interface energies, the representations \eqref{eq:locreleneps} and \eqref{eq:locrelen},
and the compactness \eqref{eqn:convergencepsi}. 
\\ \

\textbf{Acknowledgments.} The authors thank Elisabetta Rocca for the helpful initial discussions on the problem. Part of this contribution was completed while A.P. was visiting A.M. at the University of Bonn, whose hospitality is kindly acknowledged.
A.M. is funded by the Deutsche Forschungsgemeinschaft (DFG, German Research Foundation) - CRC 1720 - 539309657.
A.M. and A.P. research is sponsored by the Alexander von Humboldt Foundation, which is acknowledged.
A.P. research was funded in part by the Austrian Science Fund (FWF) \href{https://doi.org/10.55776/ESP552}{10.55776/ESP552}.
A.P. is a member of Gruppo Nazionale per l’Analisi Matematica, la Probabilità e le loro Applicazioni (GNAMPA) of
Istituto Nazionale per l’Alta Matematica (INdAM).  For open access purposes, the authors have applied a CC BY public copyright license to any author accepted manuscript version arising from this submission.


\appendix

\section{A lower semicontinuity result}
For the reader's convenience, we provide the statement of a lower semicontinuity result due to Ioffe \cite{Ioffe}. 
For the purpose of this work, this result is useful in order to pass to the limit in the quadratic terms terms on the right hand side of \eqref{eq:entropybal}, obtaining an entropy production inequality in the limit, namely \eqref{eq:weakent} at the level of weak solutions for the non-isothermal phase-field model and \eqref{eqn:entropydissipBV} for weak $BV$ solutions to the sharp interface model.
\begin{teor}[Ioffe's theorem] \label{ThIoffe}
	Let \(Q\subset \mathbb{R}^d\) be a smooth, bounded, open subset and let \(f: Q \times \mathbb{R}^l \times \mathbb{R}^m \rightarrow [0, + \infty]\), \(d, l, m \in \mathbb{N}, \, d,l, m\geq 1,\) be a measurable non-negative function such that: 
	\[
	f(y, \cdot, \cdot) \text{ is lower semicontinuous on  } \mathbb{R}^l \times \mathbb{R}^m \text{ for every } y \in Q, 
	\]
	\[
	f(y, w, \cdot) \text{ is convex on  } \mathbb{R}^m  \text{ for every } (y, w) \in Q \times \mathbb{R}^l.
	\]
	Let also \((w_n, v_n)\), \((w,v) : Q \rightarrow \mathbb{R}^l \times \mathbb{R}^m\) be measurable functions such that
	\[
	w_n \rightarrow w \text{ almost everywhere in } Q , \; \;  v_n \rightharpoonup v \text{ in }  L^1(Q).
	\]
	Then, 
	\[
	\liminf_{n \rightarrow +\infty } \int_{Q} f(y, w_n(y), v_n(y)) \, {\rm d}y \geq \int_{Q} f(y, w(y), v(y)) \, {\rm d}y .
	\]
\end{teor}

\bibliography{Biblio_nonisoNSAC}

\bibliographystyle{abbrv}

\end{document}